\documentclass[12pt]{article}
\usepackage{amssymb}
\pdfoutput=1
\usepackage{epsfig}
\usepackage{graphicx}
\usepackage{subfigure}
\usepackage{rotating}
\usepackage{longtable}

\usepackage{amssymb}
\usepackage{amsmath}

\usepackage{hyperref}

\usepackage[round]{natbib}
\begin{document}

\newenvironment {proof}{{\noindent\bf Proof.}}{\hfill $\Box$ \medskip}

\newtheorem{theorem}{Theorem}[section]
\newtheorem{lemma}[theorem]{Lemma}
\newtheorem{condition}[theorem]{Condition}
\newtheorem{proposition}[theorem]{Proposition}
\newtheorem{remark}[theorem]{Remark}
\newtheorem{hypothesis}[theorem]{Hypothesis}
\newtheorem{corollary}[theorem]{Corollary}
\newtheorem{example}[theorem]{Example}
\newtheorem{definition}[theorem]{Definition}

\renewcommand {\theequation}{\arabic{section}.\arabic{equation}}
\def \non{{\nonumber}}
\def \hat{\widehat}
\def \bar{\overline}
\def \tilde{\widetilde}
\def \nd{\and}
\def \cprime{$'$}

\title{Separation of time-scales and model reduction
for stochastic reaction networks
 \thanks{Research supported in part by NSF grants
DMS 05-53687 and 08-05793} 
}
                                                  
\author{\begin{tabular}{ll}                               
Hye-Won Kang & Thomas G. Kurtz\\
Department of Mathematics & Departments of Mathematics and Statistics \\
University of Minnesota & University of Wisconsin - Madison\\
206 Church St. S.E. & 480 Lincoln Drive\\
Minneapolis, MN 55455 & Madison, WI  53706-1388\\
hkang@math.umn.edu & kurtz@math.wisc.edu                                     
\end{tabular}}

\date{August 15, 2010}

\maketitle

\begin{abstract}
A stochastic model for a chemical reaction network is 
embedded in a one-parameter family of models with species 
numbers and rate constants scaled by powers of the 
parameter.  A systematic approach is developed for 
determining appropriate choices of the exponents that 
can be applied to large complex networks.  When the 
scaling implies subnetworks have different time-scales, 
the subnetworks can be approximated separately 
providing insight into the behavior of the full network 
through the analysis of these lower dimensional 
approximations.

\vspace{.1in}

\noindent
{\bf MSC 2000 subject classifications:  }
60J27, 60J80, 60F17, 92C45, 80A30\hfill\break

\end{abstract}
                                             
\setcounter{equation}{0}

\section{Introduction}
Chemical reaction networks in biological cells involve 
chemical species with vastly differing numbers of 
molecules and reactions with rate constants that also 
vary over several orders of magnitude.  This wide 
variation in number and rate yield phenomena that 
evolve on very different time-scales.  As in many other 
areas of application, these differing time-scales can be 
exploited to obtain simplifications of complex models.  
Papers by \citet*{RA03} and \citet*{HR02} stimulated 
considerable interest in this approach and notable 
contributions by \citet*{CGP05}, \citet*{Gou05}, 
\citet*{ELV07}, \citet*{MHR07}, \citet*{CDR09}, and others. All of the cited work 
considers
models of chemical reaction networks given by 
continuous time Markov chains where the state of the 
chain is an integer vector whose components give the 
numbers of molecules of each of the chemical species 
involved in the reaction.  Most of the analysis carried 
out in this previous work is based on the chemical 
master equation (the Kolmogorov forward equation) 
determining the one-dimensional distributions of the 
process and is focused on simplifying simulation methods 
for the process.  In contrast, the analysis 
in \citet*{BKPR06}, is based primarily on 
stochastic equations determining the process and focuses 
on the derivation of simplified models obtained as limits 
of rescaled versions of the original model.

The present 
paper gives a systematic development of many of the ideas 
introduced in   \citet{BKPR06}.  
First, recognizing that the variation in 
time-scales is due both to variation in species number 
and to variation in rate constants, we normalize species 
numbers and rate constants by powers of a fixed 
constant $N_0$ which we assume to be ``large.''  

Second, we 
replace $N_0$ by a parameter $N$ to obtain a one-parameter 
family of models and obtain our approximate models as 
rigorous limits as $N\rightarrow\infty$.  It is natural to compare this 
approach to singular perturbation analysis of 
deterministic models (cf.  \citet*{SS89}) and many of the 
same ideas and problems arise.  This kind of 
analysis is implicit in some of the earlier work and is 
the basis for the work in \citet{BKPR06}.  

Third, as in  \citet*{BKPR06}, the different 
time-scales are identified with powers $N_0^{\gamma}$, and making a 
change of time variable (replacing $t$ by $tN^{\gamma}$) we get 
different limiting/approximate models involving different 
subsets of the chemical species.  As observed in 
\citet*{CGP05} and \citet*{ELV07}, the variables in the 
approximate models may correspond to linear 
combinations of species numbers.  
We identify the time-scale of a species or a reaction 
with the exponent $\gamma$  for which the asymptotic behavior 
is nondegenerate, that is, the quantity has a nonconstant, 
well-behaved limit.
The 
time-scale of a reaction is determined by 
the scaling of its rate constant and by the scaling of the 
species numbers of the species that determine the 
intensity/propensity function for the reaction.
The time-scale of a 
species will depend both on the scaling of the 
intensity/propensity functions (the reaction time-scales) 
and on the scaling of the species number. It can 
happen that the scaling of a species number will need to 
be different for different time scales, and a species may 
appear in the limiting model for more than one of the 
time scales.

Fourth, the limiting models may be stochastic, 
deterministic or ``hybrid'' involving stochastically driven 
differential equations, that is, {\em piecewise deterministic }
Markov processes (see \cite{Dav93}). 
\citet*{HR02} obtain hybrid 
models and hybrid models have been used elsewhere in reaction 
network modeling (for example, \citet*{HRY09}, \citet*{ZFL10} ) and are 
a primary focus of \citet*{CDR09}.

Finally, as in \citet{BKPR06}, we carry out our analysis 
using stochastic equations of the form 
\[X(t)=X(0)+\sum_kY_k(\int_0^t\lambda_k(X(s))ds)\zeta_k\]
that determine the continuous time Markov chain model.  
Here the $Y_k$ are independent unit Poisson processes and 
the $\zeta_k$ are vectors in ${\Bbb Z}^d$.  These equations are
rescaled and the analysis carried out exploiting the law 
of large numbers and martingale properties of the $Y_k$.  
(For more information, see \cite{Kur77} and 
\citet*{EK86}, Chapter 11.)  The other critical component of the 
analysis is averaging methods that date back at least 
to \cite{Kha66a,Kha66b}.  (We follow \cite{Kur92}.  See 
that paper for additional references.)

If $N_0$ is large but not large enough,
the limiting model obtained by the procedure outlined 
above may have components that exhibit no fluctuation
but corresponding to components in the original model 
that exhibit substantial fluctuation. This observation 
suggests the possibility of some kind of 
diffusion/Langevin approximation.
Under what we will call the classical scaling (see Section 
\ref{sectmodel}), 
diffusion/Langevin 
approximations can be determined simply by replacing the 
rescaled Poisson processes by their appropriate Brownian 
approximations.  In systems with multiple time-scales 
that involve averaging fast components, fluctuations 
around averaged quantities may also contribute to the 
diffusion terms, and identifying an appropriate diffusion 
approximation becomes more delicate.  
These ``higher order'' corrections will 
be discussed in a later paper, \citet*{KKP10}.

Section \ref{sectmodel}\ introduces the general 
class of models to be considered and defines the scaling 
parameters used in our approach.  For comparison 
purposes, we will also describe 
the ``classical scaling'' that leads to the deterministic 
law of mass action.  Section \ref{sectscl} describes 
systematic approaches to the selection of the scaling 
parameters.  Unfortunately, even with these methods 
there may be as much art as science in their selection, 
although perhaps we should claim that this is a ``feature'' 
(flexibility) rather than a ``bug'' (ambiguity). Section 
\ref{sectlim}\ discusses identification of principal 
time-scales and derivation of the limiting models.  
Section \ref{sectave}\ reviews general averaging methods, 
and Section \ref{sectexamp} gives additional examples.  
We believe that these methods provide tools for the 
systematic reduction of highly complex models.  Further 
evidence for that claim is provided in \cite{Kan10} in
which the methods are applied to obtain a three 
time-scale reduction of a model of the heat shock 
response in {\em E. coli\/} given by \citet*{SPB01}. 

\subsection{Terminology}
This paper relies on work in both the stochastic 
processes and the chemical physics and biochemical 
literature.  Since the two communities use different 
terminology, we offer a brief translation 
table.

\begin{center}
\begin{tabular}{lcl}
{\bf Chemistry }& & {\bf Probability }\\
\\
propensity &\qquad& intensity\\
master equation &\qquad& forward equation\\
Langevin approximation && diffusion 
approximation\\
Van Kampen approximation && central limit theorem\\
quasi steady state/partial equilibrium analysis&& averaging\\
\end{tabular}
\end{center}

The terminology in the last line is less settled on both 
sides, and the methods we will discuss in Section 
\ref{sectave}\ may not yield ``averages'' at all, although 
when they don't they still correspond well to the 
quasi-steady state assumption in the chemical literature.

\subsection{Acknowledgments} The authors thank the 
other members of the NSF sponsored Focused Research 
Group on Intracellular Reaction Networks, David 
Anderson, George Craciun, Lea Popovic, Greg Rempala and 
John Yin, for many helpful conversations during the 
long gestation period of the ideas presented here.  They 
provided many valuable insights and much 
encouragement.  This work was completed while the 
first author held a postdoctoral appointment under Hans 
Othmer at the University of Minnesota and the second 
author was a Visiting Fellow at the Isaac Newton Institute in 
Cambridge, UK.  The hospitality and support provided by 
these appointments is gratefully acknowledged.
                                             
\setcounter{equation}{0}

\section{Equations for the system state}\label{sectmodel}
The standard notation for a chemical reaction
\[A+B\rightharpoonup C\]
is interpreted as ``a molecule of $A$ combines with a 
molecule of $B$ to give a molecule of $C$.''
\[A+B\rightleftharpoons C\]
means that the reaction can go in either direction, that 
is, in addition to the previous reaction, a molecule of $C$ 
can dissociate into a molecule of $A$ and a molecule of $B$.  
We consider a {\em network\/} of reactions involving $s_0$ chemical 
species, 
$S_1,\ldots ,S_{s_0}$, and $r_0$ chemical reactions
\[\sum_{i=1}^{s_0}\nu_{ik}S_i\rightharpoonup\sum_{i=1}^{s_0}\nu_{
ik}'S_i,\quad k=1,\ldots ,r_0,\]
where the $\nu_{ik}$ and $\nu_{ik}'$ are nonnegative integers.  If the 
$k$th reaction occurs, then for $i=1,\ldots ,s_0$, $\nu_{ik}$ molecules 
of $S_i$
are consumed and $\nu'_{ik}$ molecules are produced. We write 
reversible reactions as two separate reactions.

Let $X(t)\in {\Bbb N}^{s_0}$ be the vector whose components give the 
numbers of molecules of each species in the system at 
time $t$.  Let $\nu_k$ be the vector with components $\nu_{ik}$ and $
\nu_k'$ the 
vector with components $\nu_{ik}'$. 
If 
the $k$th reaction occurs at time $t$, then the state 
satisfies
\[X(t)=X(t-)+\nu'_k-\nu_k.\]
If $R_k(t)$ is the number of times that the $k$th reaction 
occurs by time $t$, then
\begin{eqnarray*}
X(t)&=&X(0)+\sum_kR_k(t)(\nu'_k-\nu_k)=X(0)+(\nu'-\nu )R(t),\end{eqnarray*}
where $\nu'$ is the $s_0\times r_0$-matrix with columns given by the $
\nu'_k$,
$\nu$ is the matrix with columns given by the $\nu_k$, and
$R(t)\in {\Bbb N}^{r_0}$ 
is the vector with components $R_k(t)$.

Modeling $X$ as a continuous time Markov chain, we can 
write
\begin{equation}R_k(t)=Y_k(\int_0^t\lambda_k(X(s))ds),\label{cntr}\end{equation}
where the $Y_k$ are independent unit Poisson processes and
$\lambda_k(x)$ is the rate at which the $k$th
reaction occurs if the chain is in state $x$, that is, $\lambda_k
(X(t))$ 
gives the {\em intensity\/} ({\em propensity\/} in the  chemical 
literature) for the $k$th reaction.  Then $X$ is the solution 
of
\begin{equation}\\
X(t)=X(0)+\sum_kY_k(\int_0^t\lambda_k(X(s))ds)(\nu'_k-\nu_k).\label{speq}\end{equation}

Define $\zeta_k=\nu'_k-\nu_k$.  The generator of the process has the form
\[{\Bbb B}f(x)=\sum_k\lambda_k(x)(f(x+\zeta_k)-f(x)).\]
Assuming that the solution of (\ref{speq}) exists for all 
time, that is, $X$ jumps only finitely often in a finite 
time interval,
\begin{equation}f(X(t))-f(X(0))-\int_0^t{\Bbb B}f(X(s))ds\label{mgp}\end{equation}
is at least a local martingale for all functions on the 
state space of the process $X$.

 If (\ref{mgp}) is a martingale, 
then its expectation is zero and 
\begin{equation}\sum_xf(x)p(x,t)=\sum_xf(x)p(x,0)+\int_0^t{\Bbb B}
f(x)p(x,s)ds,\label{wkfwd}\end{equation}
where $p(x,t)=P\{X(t)=x\}$.  Taking $f(x)={\bf 1}_{\{y\}}(x)$, 
(\ref{wkfwd}) gives the Kolmogorov forward equations (or 
{\em master equation\/} in the chemical literature) 
\begin{equation}\dot {p}(y,t)=\sum_k\lambda_k(y-\zeta_k)p(y-\zeta_
k,t)-\sum_k\lambda_k(y)p(y,t).\label{fwd}\end{equation}

The stochastic equation (\ref{speq}), the martingales 
(\ref{mgp}), and the forward equation (\ref{fwd}) provide 
three different ways of specifying the same model.  
This paper focuses primarily on the stochastic equation 
which seems to be the simplest approach to identifying 
and analyzing 
the rescaled families of models that we will introduce.

In what follows, we will focus on reactions that are at most 
binary (that is, consume at most two molecules), so $\lambda_k(x)$ 
must have one of the following forms:

\begin{center}
\begin{tabular}{lcc}
$\lambda_k$ & Reaction & $\nu_k$ \\ \hline
$\kappa_k'$ & $\emptyset\rightarrow\mbox{\rm stuff}$ & 0 \\
$\kappa_k'x_i$ & $S_i\rightarrow\mbox{\rm stuff}$ & $e_i$ \\
$\kappa_k'V^{-1}x_i(x_i-1)$ & $2S_i\rightarrow\mbox{\rm stuff}$ & $2e_
i$ \\
$\kappa_k'V^{-1}x_ix_j$ & $S_i+S_j\rightarrow\mbox{\rm stuff}$ & $e_
i+e_j$
\end{tabular}
\end{center}
Here $V$ denotes some measure of the volume of the 
system, and the form of the rates reflects the fact that 
the rate of a binary reaction in a well-stirred system
should vary inversely with 
the volume of the system.
Note that if $\zeta_{ik}<0$, then $\lambda_k(x)$ must have $x_i$ as a factor.  
Higher order reactions can be included at the cost of 
more complicated expressions for the $\lambda_k$.

Our intent is to embed the model of primary interest $X$ into a 
family of models $X^N$ indexed by a large parameter $N$.  The 
model $X$ corresponds to a particular value of 
the parameter $N=N_0$, that is $X=X^{N_0}$.

For each species 
$i$, let $\alpha_i\geq 0$ and define the {\em normalized abundance\/} (or 
simply, the abundance) for the $N$th model by
\[Z_i^N(t)=N^{-\alpha_i}X_i^N(t).\]
Note that the abundance may be the species number 
($\alpha_i=0$),
the species concentration, or something else.
The exponent $\alpha_i$ should be selected so that $Z_i^N=O(1)$.  
To be precise, we want 
 $\{Z_i^N(t)\}$ to be
stochastically bounded,
that is, for each $\epsilon >0$, there exists 
 $K_{\epsilon ,t}<\infty$  such that
\[\inf_NP\{\sup_{s\leq t}Z_i^N(s)\leq K_{\epsilon ,t}\}\geq 1-\epsilon 
.\]
In other words, we want $\alpha_i$ to be ``large enough.''  On 
the other hand, we do not want $\alpha_i$ to be so large that 
$Z_i^N$ converges to zero as $N\rightarrow\infty$.  For example, the existence of 
$\delta_{\epsilon}$ such that
\[\inf_NP\{\inf_{s\leq t}Z_i^N(s)\geq\delta_{\epsilon ,t}\}\geq 1
-\epsilon\]
would suffice; however, there are natural situations in 
which   $\alpha_i=0$ and $Z_i^N$ is occasionally or even 
frequently zero, so this requirement would in general be 
too restrictive.  For the moment, we just keep in mind 
that $\alpha_i$ cannot be ``too big.''

The rate constants may also vary over several orders of 
magnitude,  
so we define $\kappa_k$ by setting
$\kappa_k'=\kappa_kN_0^{\beta_k}$ for unary reactions and $\kappa_
k'V^{-1}=\kappa_kN_0^{\beta_k}$ for 
binary reactions.
 The $\beta_k$ should be selected so 
that the  $\kappa_k$  are of order one, although we again avoid 
being too precise regarding the meaning of ``order one.'' 
For a unary 
reaction, the intensity for the model of primary interest 
becomes
\[\kappa_k'x_i=N_0^{\beta_k+\alpha_i}z_i=N_0^{\beta_k+\nu_k\cdot\alpha}
z_i,\]
and for binary reactions, 
\[\kappa_k'V^{-1}x_ix_j=N_0^{\beta_k+\alpha_i+\alpha_j}\kappa_kz_
iz_j=N_0^{\beta_k+\nu_k\cdot\alpha}\kappa_kz_iz_j\]
and 
\begin{equation}\kappa_k'V^{-1}x_i(x_i-1)=N_0^{\beta_k+2\alpha_i}
\kappa_kz_i(z_i-N_0^{-\alpha_i})=N_0^{\beta_k+\nu_k\cdot\alpha}\kappa_
kz_i(z_i-N_0^{-\alpha_i}).\label{redi}\end{equation}

The $N$th model in the scaled family is given by the system
\[Z^N_i(t)=Z^N_i(0)+\sum_kN^{-\alpha_i}Y_k(\int_0^tN^{\beta_k+\nu_
k\cdot\alpha}\lambda_k(Z^N(s))ds)(\nu_{ik}'-\nu_{ik}).\]
For binary reactions of the form $2S_i\rightarrow\mbox{\rm stuff}$ with $
\alpha_i>0$, 
$\lambda_k(z)=\kappa_kz_i(z_i-N^{-\alpha_i})$ 
depends on $N$, but to simplify notation we still write
$\lambda_k$ rather than $\lambda_k^N$.

Let $\Lambda_N=\mbox{\rm diag}(N^{-\alpha_1},\ldots ,N^{-\alpha_{
s_0}})$, $\rho_k=\beta_k+\nu_k\cdot\alpha$, and 
$\zeta_k=\nu_k'-\nu_k$.  The generator for $Z^N$ is 
\[{\Bbb B}_Nf(z)=\sum_kN^{\rho_k}\lambda_k(z)(f(z+\Lambda_N\zeta_
k)-f(z)).\]

Even after the $\beta_k$ and $\alpha_i$ are selected, we still have the 
choice of time-scale on which to study the model, that 
is, we can consider
\begin{equation}Z^{N,\gamma}_i(t)=Z^N_i(tN^{\gamma})=Z^N_i(0)+\sum_
kN^{-\alpha_i}Y_k(\int_0^tN^{\gamma +\beta_k+\nu_k\cdot\alpha}\lambda_
k(Z^{N,\gamma}(s))ds)(\nu_{ik}'-\nu_{ik})\label{zgmn}\end{equation}
for any $\gamma\in {\Bbb R}$.  Different choices of $\gamma$ may give 
interesting approximations for different subsets of species.  
To identify that approximation, note that if 
$\lim_{N\rightarrow\infty}Z^{N,\gamma}_i=Z^{\gamma}_i$ and $N_0$ is ``large'', then we should have 
\[X_i(t)\equiv X^{N_0}_i(t)\approx N_0^{\alpha_i}Z^{\gamma}_i(tN^{
-\gamma}_0).\]

In what we will call the {\em classical scaling\/} (see, for example, 
\cite{Kur72,Kur77}) $N_0$ has the interpretation of 
volume times Avogadro's number and $\alpha_i=1$, for all $i$,
so  $Z_i^{N_0}$ is the 
concentration of $S_i$.  
Taking 
$\beta_k=0$ for a 
unary reaction and $\beta_k=-1$ for a binary reaction, the 
intensities are all of the form $N\lambda_k(z)$, and hence taking 
$\gamma =0$, $Z^N=Z^{N,0}$ converges to the solution of
\begin{equation}Z_i(t)=Z_i(0)+\sum_k\int_0^tZ(s)^{\nu_k}ds(\nu_{i
k}'-\nu_{ik}),\label{lma}\end{equation}
where $z^{\nu_k}=\prod_iz_i^{\nu_{ik}}$.  Note that (\ref{lma}) is just the 
usual {\em law of mass action\/} model for the network.

\setcounter{equation}{0}

\section{Determining the scaling exponents}\label{sectscl}
For systems with a diversity of scales because 
of wide variations in species numbers or rate constants 
or both, 
the challenge is to select the $\alpha_i$ and the $\beta_k$ in ways that 
capture  this variation and produce interesting 
approximate models.
Once the exponents and $N_0$ are selected,
\[X_i^N(0)=\lfloor\left(\frac N{N_0}\right)^{\alpha_i}X_i(0)\rfloor 
,\]
and the family of models to be studied is determined.

Suppose
\[\kappa_1'\geq\kappa_2'\geq\cdots\geq\kappa'_{r_0}.\]
Then it is reasonable to select the $\beta_i$ so that 
$\beta_1\geq\cdots\geq\beta_{r_0}$, although it may be natural to impose this 
order separately for unary and binary reactions.  (See 
the ``classical'' scaling.)

Typically, we want to select the $\alpha_i$ so that 
$Z_i^N(t)=N^{-\alpha_i}X_i^N(t)=O(1)$, or more precisely, assuming 
$\lim_{N\rightarrow\infty}Z_i^N(0)=Z_i(0)>0$, for all $i$, we want to 
avoid $\alpha$, $\beta$, and $\gamma$ for which
$\lim_{N\rightarrow\infty}Z^N_i(tN^{\gamma})=0,$
for all $t>0$ or $\lim_{N\rightarrow\infty}Z^N_i(tN^{\gamma})=\infty$, for all $
t>0$.  This 
goal places constraints on $\alpha$, $\beta$, and possibly $\gamma$.

\subsection{Species balance}
Consider the reaction system
\begin{eqnarray*}
S_1+S_2&\rightharpoonup&S_3+S_4\\
S_3+S_5&\rightharpoonup&S_6.\end{eqnarray*}
Then the equation for $Z_3^{N,\gamma}$ is
\begin{eqnarray*}
Z_3^{N,\gamma}(t)=Z_3^N(0)+N^{-\alpha_3}Y_1(N^{\gamma +\beta_1+\alpha_
1+\alpha_2}\int_0^t\kappa_1Z_1^{N,\gamma}(s)Z_2^{N,\gamma}(s)ds)\\
-N^{-\alpha_3}Y_2(N^{\gamma +\beta_2+\alpha_3+\alpha_5}\int_0^t\kappa_
2Z_3^{N,\gamma}(s)Z_5^{N,\gamma}(s)ds)\;.\end{eqnarray*}
Assuming that $Z_i^{N,\gamma}=O(1)$ for $i\neq 3$ and $Z_3^N(0)=O
(1)$,
 $Z_3^{N,\gamma}=O(1)$ if 
\[(\beta_1+\alpha_1+\alpha_2+\gamma )\vee (\beta_2+\alpha_3+\alpha_
5+\gamma )\leq\alpha_3\]
(the power of $N$ outside the Poisson processes dominates 
the power inside) or if
\begin{equation}\beta_1+\alpha_1+\alpha_2=\beta_2+\alpha_3+\alpha_
5.\label{fbal}\end{equation}
Assuming (\ref{fbal}), if $Z_3^{N,\gamma}>\frac {\kappa_1Z_1^{N,\gamma}
(s)Z_2^{N,\gamma}(s)ds}{\kappa_2Z_5^{N,\gamma}(s)}$, the rate of 
consumption of $S_3$ exceeds the rate of production, and if 
the inequality is reversed, the rate of production 
exceeds the rate of consumption ensuring that $Z_3^N$ 
neither explodes nor is driven to zero. 

In general, let $\Gamma_i^{+}=\{k:\nu_{ik}'>\nu_{ik}\}$, that is, $
\Gamma_i^{+}$ gives the set of 
reactions that result in an increase in the $i$th species, 
and let $\Gamma_i^{-}=\{k:\nu_{ik}'<\nu_{ik}\}$.  Then for each $
i$, we want 
either
\begin{equation}\max_{k\in\Gamma_i^{-}}(\beta_k+\nu_k\cdot\alpha 
)=\max_{k\in\Gamma^{+}_i}(\beta_k+\nu_k\cdot\alpha ).\label{sbal1}\end{equation}
or
\begin{equation}\max_{k\in\Gamma^{+}_i\cup\Gamma_i^{-}}(\beta_k+\nu_
k\cdot\alpha )+\gamma\leq\alpha_i.\label{sbal2}\end{equation}
We will refer to (\ref{sbal1}) as the {\em balance equation\/} for species $
i$ 
and to (\ref{sbal2}) as a {\em time-scale constraint\/} since it is 
equivalent to
\[\gamma\leq\alpha_i-\max_{k\in\Gamma^{+}_i\cup\Gamma_i^{-}}(\beta_
k+\nu_k\cdot\alpha ).\]
The requirement that either a
species be balanced or the time-scale 
constraint be satisfied will be called the
{\em species balance condition}.

Equation (\ref{sbal1}) is the requirement that 
the maximum rate at which a species is produced is of the 
same order of magnitude as the rate at which it is 
consumed.  Since consumption rates are proportional to 
the normalized species state $Z_i$, $Z_i$ should 
remain $O(1)$ provided the same is true for the other $Z_j$ 
even if the normalized reaction numbers blow up.
If (\ref{sbal1}) fails to hold, then
(\ref{sbal2}) ensures that $Z_i(t)=O(1)$, again provided the 
other $Z_j$ remain $O(1)$.

Note that if $\zeta_{ik}\neq 0$, then 
\begin{equation}\gamma =\gamma_{ik}=\alpha_i-(\beta_k+\nu_k\cdot\alpha 
)\label{natr}\end{equation}
is in 
some sense the {\em natural time-scale\/} for the normalized 
reaction number 
\[N^{-\alpha_i}R_k^{N,\gamma}(t)=N^{-\alpha_i}Y_k(N^{\gamma +\beta_
k+\nu_k\cdot k}\int_0^t\lambda_k(Z^{N,\gamma}(s))ds).\]
Then, regardless of whether (\ref{sbal1}) or (\ref{sbal2}) holds,
\begin{equation}\gamma_i=\min_{k\in\Gamma^{+}_i\cup\Gamma_i^{-}}\gamma_{
ik}=\alpha_i-\max_{k\in\Gamma^{+}_i\cup\Gamma_i^{-}}(\beta_k+\nu_
k\cdot\alpha )\label{nati}\end{equation}
is the natural time-scale for species $S_i$.  
With reference to (\ref{zgmn}), if $\gamma <\gamma_i$, we expect 
$Z_i^{N,\gamma}(t)$ to converge to $\lim_{N\rightarrow\infty}Z_i^
N(0)$.  If $\gamma =\gamma_i$ and 
$\alpha_i>0$, then 
we expect
\[\lim_{N\rightarrow\infty}Z_i^{N,\gamma_i}(t)=\lim_{N\rightarrow
\infty}(Z_i^N(0)+\sum_{k\in\Gamma_{i,0}}\int_0^t\lambda_k(Z^{N,\gamma_
i}(s))ds(\nu_{ik}'-\nu_{ik})),\]
where 
\[\Gamma_{i,0}=\{l:\beta_l+\nu_l\cdot\alpha =\max_{k\in\Gamma^{+}_
i\cup\Gamma_i^{-}}(\beta_k+\nu_k\cdot\alpha )\}\]
and each integral on the right side is nonconstant but 
well behaved.  If $\alpha_i=0$, we expect
\[\lim_{N\rightarrow\infty}Z_i^{N,\gamma_i}(t)=\lim_{N\rightarrow
\infty}(Z_i^N(0)+\sum_{k\in\Gamma_{i,0}}Y_k(\int_0^t\lambda_k(Z^{
N,\gamma_i}(s))ds)(\nu_{ik}'-\nu_{ik})).\]

It is important to notice that we associate ``time-scales'' 
with species (and as we will see below, with collections 
of species) and that one
reaction may  determine different time-scales associated 
with different species.

\subsection{Collective species balance}
The species balance condition, however, does 
not by itself ensure that the normalized species numbers 
are asymptotically all $O(1)$.  There may also be subsets of 
species such that the collective rate of production is 
of a different order of magnitude than the collective 
rate of consumption.
Consider the following simple 
network:
\[\emptyset\mathop{\rightharpoonup}^{\kappa_1}S_1\mathop{\rightleftharpoons}^{
\kappa_2}_{\kappa_3}S_2\mathop{\rightharpoonup}^{\kappa_4}\emptyset 
.\]
If $0<\beta_4<\beta_1<\beta_2=\beta_3$ and $\alpha_1=\alpha_2=0$, then 
\begin{eqnarray}
Z_1^N(t)&=&Z_1^N(0)+Y_1(\kappa_1N^{\beta_1}t)+Y_3(\kappa_3N^{\beta_
3}\int_0^tZ_2^N(s)ds)-Y_2(\kappa_2N^{\beta_2}\int_0^tZ^N_1(s)ds)\nonumber\\
Z_2^N(t)&=&Z^N_2(0)+Y_2(\kappa_2N^{\beta_2}\int_0^tZ^N_1(s)ds)-Y_
3(\kappa_3N^{\beta_3}\int_0^tZ_2^N(s)ds)\label{unbal}\\
&&\qquad\qquad\qquad\qquad\qquad\qquad\qquad -Y_4(\kappa_4N^{\beta_
4}\int_0^tZ^N_2(s)ds)\;.\nonumber\end{eqnarray}
Since $\beta_2=\beta_3\vee\beta_1$ and $\beta_2=\beta_3\vee\beta_
4$,
the species balance condition is 
satisfied for all species, 
but noting that 
\[Z_1^N(t)+Z_2^N(t)=Z_1^N(0)+Z_2^N(0)+Y_1(\kappa_1N^{\beta_1}t)-Y_
4(\kappa_4N^{\beta_4}\int_0^tZ^N_2(s)ds),\]
the species numbers still go to infinity as 
$N\rightarrow\infty$.  This example suggests the need to consider linear 
combinations of species. These linear combinations may, 
in fact, play the role of ``virtual'' species or 
auxiliary variables needed in the specification of the 
reduced models (cf.  
\citet*{CGP05} and \citet*{ELV05,ELV07}).

To simplify notation, define 
\[\rho_k=\beta_k+\nu_k\cdot\alpha ,\]
so the scaled model 
satisfies
\begin{eqnarray*}
Z^{N,\gamma}(t)=Z^{N,\gamma}(0)+\Lambda_N\sum_kY_k(N^{\beta_k+\nu_
k\cdot\alpha +\gamma}\int_0^t\lambda_k(Z^{N,\gamma}(s))ds)\zeta_k\\
=Z^{N,\gamma}(0)+\Lambda_N\sum_kY_k(N^{\rho_k+\gamma}\int_0^t\lambda_
k(Z^{N,\gamma}(s))ds)\zeta_k,\end{eqnarray*}
where $\Lambda_N$ is the diagonal matrix with entries $N^{-\alpha_
i}$.

\begin{definition}
For $\theta\in [0,\infty )^{s_0}$, define $\Gamma_{\theta}^{+}=\{
k:\theta\cdot\zeta_k>0\}$ and 
$\Gamma^{-}_{\theta}=\{k:\theta\cdot\zeta_k<0\}$.  
\end{definition}

Then, noting that 
\[\theta^T\Lambda_N^{-1}Z^{N,\gamma}(t)=\sum_{i=1}^{s_0}\theta_iN^{
\alpha_i}Z_i^{N,\gamma}(t)=\sum_{i=1}^{s_0}\theta_iX^N_i(N^{\gamma}
t),\]
\begin{eqnarray*}
\theta^T\Lambda_N^{-1}Z^{N,\gamma}(t)&=&\theta^T\Lambda_N^{-1}Z^{
N,\gamma}(0)+\sum_k(\theta\cdot\zeta_k)Y_k(N^{\rho_k+\gamma}\int_
0^t\lambda_k(Z^{N,\gamma}(s))ds)\\
&=&\theta^T\Lambda_N^{-1}Z^{N,\gamma}(0)+\sum_{k\in\Gamma^{+}_{\theta}}
(\theta\cdot\zeta_k)R_k^{N,\gamma}(t)-\sum_{k\in\Gamma^{-}_{\theta}}
|(\theta\cdot\zeta_k)|R_k^{N,\gamma}(t).\end{eqnarray*}
To avoid some kind of degeneracy in the limit, either 
the positive and negative sums must cancel, or they 
must grow no faster than $N^{\alpha_i}$ for some $i$ with $\theta_
i>0$.  
Consequently, we extend the species balance condition to 
linear combinations of species.
For each $\theta\in [0,\infty )^{s_0}$, the following condition must hold.  

\begin{condition}\label{cndz}
\begin{equation}\max_{k\in\Gamma_{\theta}^{-}}(\beta_k+\nu_k\cdot
\alpha )=\max_{k\in\Gamma^{+}_{\theta}}(\beta_k+\nu_k\cdot\alpha 
)\label{cnd1}\end{equation}
or
\begin{equation}\gamma\leq\gamma_{\theta}\equiv\max_{i:\theta_i>0}
\alpha_i-\max_{k\in\Gamma^{+}_{\theta}\cup\Gamma_{\theta}^{-}}(\beta_
k+\nu_k\cdot\alpha ).\label{cnd2}\end{equation}
\end{condition}

Of course, if $\theta_i>0$ for only a single species, then this 
requirement is just the species balance condition, so 
Condition \ref{cndz} includes that
condition.  Again, we will refer to (\ref{cnd1}) as the 
balance equation for the linear 
combination $\theta\cdot X=\sum_i\theta_iX_i$.  In the special 
case of $\theta =e_i$, the vector with  $i$th component $1$ and 
other components $0$, we say that $X_i$ is balanced
 or that the species 
$S_i$ is balanced.  If (\ref{cnd1}) 
fails for $\theta$, we say that $\theta\cdot X$ is {\em unbalanced}.  
The inequalities given by  (\ref{cnd2}) are again called 
{\em time-scale constraints\/} as they 
imply
\begin{equation}\gamma\leq\min_{\theta\cdot X\mbox{\rm \ unbalanced}}
\gamma_{\theta}.\label{gamcnst}\end{equation}

For example, consider the network
\[\emptyset\mathop{\rightharpoonup}^{\kappa_1}S_1\mathop{\rightleftharpoons}^{
\kappa_2}_{\kappa_3}S_2,\]
and assume that $\kappa_k'=\kappa_kN_0^{\beta_k}$, where $\beta_1
=\beta_2>\beta_3$.  For $S_2$ 
to be balanced, we must have $\beta_2+\alpha_1=\beta_3+\alpha_2$ and for $
S_1$ 
to be balanced, we must have
\[\beta_1\vee (\beta_3+\alpha_2)=\beta_2+\alpha_1.\]
Let  $\alpha_1=0$ and $\alpha_2=\beta_2-\beta_3$ so $S_1$ and $S_
2$ are balanced.  
For $\theta =(1,1)$,  
$\Gamma_{\theta}^{+}=\{1\}$, and $\Gamma^{-}_{\theta}=\emptyset$.  Consequently, (\ref{cnd1}) fails, so
we require
\begin{equation}\gamma\leq\alpha_1\vee\alpha_2-\beta_1=-\beta_3.\label{nbal2a}\end{equation}
There are two time-scales of interest in this model, 
$\gamma =-\beta_1$, the natural time-scale of $S_1$ and $\gamma =
-\beta_3$, the 
natural time-scale 
of $S_2$.  The system of equations is
\begin{eqnarray*}
Z^{N,\gamma}_1(t)&=&Z_1^N(0)+Y_1(\kappa_1N^{\gamma +\beta_1}t)-Y_
2(\kappa_2N^{\gamma +\beta_2}\int_0^tZ_1^{N,\gamma}(s)ds)\\
&&\qquad +Y_3(\kappa_3N^{\gamma +\beta_3+\alpha_2}\int_0^tZ_2^{N,
\gamma}(s)\\
Z_2^{N,\gamma}(t)&=&Z_2^N(0)+N^{-\alpha_2}Y_2(\kappa_2N^{\gamma +
\beta_2}\int_0^tZ_1^{N,\gamma}(s)ds)\\
&&\qquad -N^{-\alpha_2}Y_3(\kappa_3N^{\gamma +\beta_3+\alpha_2}\int_
0^tZ_2^{N,\gamma}(s).\end{eqnarray*}

For $\gamma =-\beta_1$, since $\beta_1=\beta_2=\beta_3+\alpha_2$, the limit of $
Z^{N,\gamma}$ satisfies
\begin{eqnarray*}
Z_1(t)&=&Z_1(0)+Y_1(\kappa_1t)-Y_2(\kappa_2\int_0^tZ_1(s)ds)+Y_3(
\kappa_3\int_0^tZ_2(s))\\
&=&Z_1(0)+Y_1(\kappa_1t)-Y_2(\kappa_2\int_0^tZ_1(s)ds)+Y_3(\kappa_
3Z_2(0)t)\\
Z_2(t)&=&Z_2(0).\end{eqnarray*}

For $\gamma =-\beta_3$, if we divide the equation for $Z_1^{N,\gamma}$ by 
$N^{\alpha_2}=N^{\beta_1-\beta_3}$, we see that 
\begin{eqnarray}
0&=&\lim_{N\rightarrow\infty}N^{-\alpha_2}Z^{N,\gamma}_1(t)
\label{chtave}\\
&=&\lim_{N\rightarrow\infty}N^{-\alpha_2}Z_1^N(0)+N^{-\alpha_2}Y_
1(\kappa_1N^{\gamma +\beta_1}t)-N^{-\alpha_2}Y_2(\kappa_2N^{\gamma 
+\beta_2}\int_0^tZ_1^{N,\gamma}(s)ds)\nonumber\\
&&\qquad\qquad\qquad\qquad\qquad +N^{-\alpha_2}Y_3(\kappa_3N^{\gamma 
+\beta_3+\alpha_2}\int_0^tZ_2^{N,\gamma}(s)\nonumber\\
&=&\lim_{N\rightarrow\infty}\left(\kappa_1t+\kappa_3\int_0^tZ^{N,
\gamma}_2(s)ds-\kappa_2\int_0^tZ_1^{N,\gamma}(s)ds\right)\nonumber\end{eqnarray}
and $Z_2^{N,\gamma}$ converges to 
\[Z_2(t)=Z_2(0)+\kappa_1t.\]
With reference to (\ref{nbal2a}), if $\gamma >-\beta_3$, then 
$Z_2^{N,\gamma}(t)\rightarrow\infty$, for each $t>0$, demonstrating the significance 
of the time-scale constraints.  

For $\gamma =-\beta_3$, $Z_1^{N,\gamma}$ fluctuates rapidly and does not 
converge in a functional sense.  Its behavior is captured, 
at least to some extent, by its occupation measure 
\[V_1^{N,\gamma}(C\times [0,t])=\int_0^t{\bf 1}_C(Z_1^{N,\gamma}(
s))ds.\]
Applying the generator to functions of $z_1$ 
and using the 
fact 
that $\beta_1-\beta_3=\beta_2-\beta_3=\alpha_2$, ${\Bbb B}^{N,\gamma}
f(z_1,z_2)=N^{\alpha_2}{\Bbb C}_{z_2}f(z_1)$, 
where
\[{\Bbb C}_{z_2}f(z_1)=(\kappa_1+\kappa_3z_2)(f(z_1+1)-f(z_1))+\kappa_
2z_1(f(z_1-1)-f(z_1)).\]
Then
\[f(Z_1^{N,\gamma}(t))-f(Z_1^{N,\gamma}(0))-N^{\alpha_2}\int_{{\Bbb N}
\times [0,t]}{\Bbb C}_{Z_2^{N,\gamma}(s)}f(z_1)V_1^{N,\gamma}(dz_
1\times ds)\]
is a martingale, and dividing by $N^{\alpha_2}$ and passing to the 
limit, it is not difficult to see that $V_1^{N,\gamma}$ converges to a 
measure satisfying
\[\]
\[\int_{{\Bbb N}\times [0,t]}{\Bbb C}_{Z_2(s)}f(z_1)V_1(dz_1\times 
ds)=0.\]
(See Section \ref{sectave}.)  Writing 
$V_1(dz_1\times ds)=v_s(dz_1)ds$, it follows that $v_s$ is the Poisson 
distribution with mean $\frac {\kappa_1+\kappa_3Z_2(s)}{\kappa_2}$. We will refer to $
v_s$ as 
the {\em conditional-equilibrium\/} or {\em local-averaging\/} distribution.

\subsection{Auxiliary variables}\label{sectaux}
 While (\ref{nati}) gives the natural time-scale for 
individual species, it is clear 
from examples considered by
\citet*{ELV05},  that the species time-scales may 
not be the only time-scales of interest.  For example, 
they consider the network
\[S_1\mathop{\rightleftharpoons}^{\kappa'_1}_{\kappa'_2}S_2\mathop{
\rightleftharpoons}^{\kappa'_3}_{\kappa'_4}S_3\mathop{\rightleftharpoons}^{
\kappa'_5}_{\kappa'_6}S_4\]
with $\kappa'_1,\kappa'_2,\kappa'_5,\kappa'_6>>\kappa'_3,\kappa'_
4$.  The scaled model is given by
\begin{eqnarray*}
Z_1^N(t)&=&Z_1^N(0)+N^{-\alpha_1}Y_2(\kappa_2N^{\beta_2+\alpha_2}
\int_0^tZ^N_2(s)ds)-N^{-\alpha_1}Y_1(\kappa_1N^{\beta_1+\alpha_1}
\int_0^tZ^N_1(s)ds)\\
Z_2^N(t)&=&Z_2^N(0)+N^{-\alpha_2}Y_1(\kappa_1N^{\beta_1+\alpha_1}
\int_0^tZ^N_1(s)ds)-N^{-\alpha_2}Y_2(\kappa_2N^{\beta_2+\alpha_2}
\int_0^tZ^N_2(s)ds)\\
&&\qquad\qquad +N^{-\alpha_2}Y_4(\kappa_4N^{\beta_4+\alpha_3}\int_
0^tZ^N_3(s)ds)-N^{-\alpha_2}Y_3(\kappa_3N^{\beta_3+\alpha_2}\int_
0^tZ^N_2(s)ds)\\
Z_3^N(t)&=&Z_3^N(0)+N^{-\alpha_3}Y_6(\kappa_6N^{\beta_6+\alpha_4}
\int_0^tZ^N_4(s)ds)-N^{-\alpha_3}Y_5(\kappa_5N^{\beta_5+\alpha_3}
\int_0^tZ^N_3(s)ds)\\
&&\qquad\qquad +N^{-\alpha_3}Y_3(\kappa_3N^{\beta_3+\alpha_2}\int_
0^tZ^N_2(s)ds)-N^{-\alpha_3}Y_4(\kappa_4N^{\beta_4+\alpha_3}\int_
0^tZ^N_3(s)ds)\\
Z_4^N(t)&=&Z_4^N(0)+N^{-\alpha_4}Y_5(\kappa_5N^{_{\beta_5+\alpha_
3}}\int_0^tZ^N_3(s)ds)-N^{-\alpha_4}Y_6(\kappa_6N^{\beta_6+\alpha_
4}\int_0^tZ^N_4(s)ds).\end{eqnarray*}
Assume that 
$\beta_1=\beta_2>\beta_5=\beta_6>\beta_3>\beta_4$.  Then if we look for a scaling 
under which all $\theta\cdot X$ are balanced, 
$\alpha_1=\alpha_2$, $\alpha_3=\alpha_4$, and $\alpha_2+\beta_3=\alpha_
3+\beta_4$, so 
$\alpha_3=\alpha_2+\beta_3-\beta_4$.  For definiteness, take $\alpha_
1=\alpha_2=0$.

The natural time-scale for $S_1$ and $S_2$ 
is $-\beta_1$, and the natural time-scale for $S_3$ and $S_4$ is $
-\beta_5$,
but on either of these time-scales $Z_1+Z_2$ and $Z_3+Z_4$ are 
constant.  In particular,
\begin{eqnarray*}
U_1^{N,\gamma}(t)&\equiv&Z_1^{N,\gamma}(t)+Z_2^{N,\gamma}(t)=Z_1^
N(0)+Z_1^N(0)+Y_4(\kappa_4N^{\gamma +\beta_4+\alpha_3}\int_0^tZ_3^{
N,\gamma}(s)ds)\\
&&\qquad -Y_3(\kappa_3N^{\gamma +\beta_3}\int_0^tZ_2^{N,\gamma}(s
)ds)\\
U_2^{N,\gamma}(t)&\equiv&Z_3^{N,\gamma}(t)+Z_4^{N,\gamma}(t)=Z_3^
N(0)+Z_4^N(0)-N^{-\alpha_3}Y_4(\kappa_4N^{\gamma +\beta_4+\alpha_
3}\int_0^tZ_3^{N,\gamma}(s)ds)\\
&&\qquad +N^{-\alpha_3}Y_3(\kappa_3N^{\gamma +\beta_3}\int_0^tZ_2^{
N,\gamma}(s)ds).\end{eqnarray*}
For $\gamma_1=\gamma_2=-\beta_1=-\beta_2$, $(Z_1^{N,\gamma_1},Z_2^{
N,\gamma_1})$ converges to
\begin{eqnarray*}
Z_1^{\gamma_1}(t)&=&Z_1(0)+Y_2(\kappa_2\int_0^tZ^{\gamma_1}_2(s)d
s)-Y_1(\kappa_1\int_0^tZ^{\gamma_1}_1(s)ds)\\
Z_2^{\gamma_1}(t)&=&Z_2(0)+Y_1(\kappa_1\int_0^tZ^{\gamma_1}_1(s)d
s)-Y_2(\kappa_2\int_0^tZ^{\gamma_1}_2(s)ds),\end{eqnarray*}
and for $\gamma_3=\gamma_4=-\beta_5=-\beta_6$, 
\begin{eqnarray*}
Z_3^{\gamma_3}(t)&=&Z_3(0)+\kappa_6\int_0^tZ^{\gamma_3}_4(s)ds-\kappa_
5\int_0^tZ^{\gamma_3}_3(s)ds\\
Z_4^{\gamma_3}(t)&=&Z_4(0)+\kappa_5\int_0^tZ^{\gamma_3}_3(s)ds-\kappa_
6\int_0^tZ^{\gamma_3}_4(s)ds.\end{eqnarray*}

Taking  $\gamma_{12}=-\beta_3=-(\alpha_3+\beta_4)$ and dividing the equation 
for $Z_4^{N,\gamma_{12}}$ by $N^{\beta_5-\beta_3}$, we see that 
\begin{equation}\kappa_5\int_0^tZ_3^{N,\gamma_{12}}(s)ds-\kappa_6
\int_0^tZ_4^{N,\gamma_{12}}(s)ds\rightarrow 0\label{intave}\end{equation}
and hence 
\begin{equation}\int_0^tZ_3^{N,\gamma_{12}}(s)ds-\frac {\kappa_6}{
\kappa_5+\kappa_6}\int_0^tU_2^{N,\gamma_{12}}(s)ds\rightarrow 0.\label{u2int}\end{equation}
Similarly, dividing the equation for $Z_1^{N,\gamma_{12}}$ by $N^{
\beta_2-\beta_3}$,
\[\int_0^tZ_2^{N,\gamma_{12}}(s)ds-\frac {\kappa_1}{\kappa_1+\kappa_
2}\int_0^tU_1^{N,\gamma_{12}}(s)ds\rightarrow 0.\]
Since $U_2^{N,\gamma_{12}}$ converges to $U_2(0)$ uniformly on bounded 
time intervals, 
$U_1^{N,\gamma_{12}}$ converges 
to the solution of
\[U_1^{\gamma_{12}}(t)=U_1(0)+Y_4(\frac {\kappa_4\kappa_6}{\kappa_
5+\kappa_6}U_2(0)t)-Y_3(\frac {\kappa_3\kappa_1}{\kappa_1+\kappa_
2}\int_0^tU_1^{_{12}}(s)ds).\]

Finally, taking $\gamma_{34}=-\beta_4$, as in (\ref{u2int}),
\[\int_0^tZ_3^{N,\gamma_{34}}(s)ds-\frac {\kappa_6}{\kappa_5+\kappa_
6}\int_0^tU_2^{N,\gamma_{34}}(s)ds\rightarrow 0,\]
and dividing the equation for $U_1^{N,\gamma_{23}}$ by $N^{\beta_
3-\beta_4}$,
\[\int_0^tZ_2^{N,\gamma_{34}}(s)ds-\frac {\kappa_4}{\kappa_3}\int_
0^tZ_3^{N,\gamma_{34}}(s)ds\rightarrow 0.\]
Consequently, even on this faster time-scale,
 $U_2^{N,\gamma_{34}}$ converges to $U_2(0)$ uniformly on 
bounded time intervals.

\subsection{Checking the balance conditions}
Condition \ref{cndz}\ only depends on the support of $\theta$, 
$\mbox{\rm supp}(\theta )=\{i:\theta_i\neq 0\}$, and on the signs of $
\theta\cdot\zeta_k$, so the 
condition needs to be checked for only finitely many 
$\theta$.  For $k\in \{1,\ldots ,r_0\}$, define
\[\Lambda_k^{+}=\{\theta\in [0,\infty )^{s_0}:\theta\cdot\zeta_k>
0\},\quad\Lambda_k^{-}=\{\theta\in [0,\infty )^{s_0}:\theta\cdot\zeta_
k<0\},\quad\Lambda_k^0=\{\theta\in [0,\infty )^{s_0}:\theta\cdot\zeta_
k=0\},\]
and for disjoint $\Gamma_{-}$, $\Gamma_{+}$, $\Gamma_0$ satisfying 
$\Gamma_{-}\cup\Gamma_{+}\cup\Gamma_0=\{1,\cdots ,r_0\}$, define
\[\Lambda_{\Gamma_{-},\Gamma_{+},\Gamma_0}=(\cap_{k\in\Gamma_{-}}
\Lambda_k^{-})\cap (\cap_{k\in\Gamma_{+}}\Lambda_k^{+})\cap (\cap_{
k\in\Gamma_0}\Lambda_k^0).\]
The following lemma is immediate.

\begin{lemma}\label{cndzr}
Fix $\gamma$.  Condition \ref{cndz}\ holds for all $\theta\in [0,
\infty )^{s_0}$ provided
\begin{equation}\max_{k\in\Gamma_{-}^{}}(\beta_k+\nu_k\cdot\alpha 
)=\max_{k\in\Gamma_{+}}(\beta_k+\nu_k\cdot\alpha )\label{cnd1r}\end{equation}
or
\begin{equation}\gamma\leq\min_{\theta\in\Lambda_{\Gamma_{-},\Gamma_{
+},\Gamma_0}}\max_{i:\theta_i>0}\alpha_i-\max_{k\in\Gamma_{+}\cup
\Gamma_{-}}(\beta_k+\nu_k\cdot\alpha )\label{cnd2r}\end{equation}
for all partitions $\{\Gamma_{-},\Gamma_{+},\Gamma_0\}$ for which 
$\Lambda_{\Gamma_{-},\Gamma_{+},\Gamma_0}\neq\emptyset$.
\end{lemma}

Checking the conditions of Lemma \ref{cndzr}  
could still be a formidable task.
 The next lemmas 
significantly
reduce the effort required.
Observe that for
 $\theta^1,\theta^2\in [0,\infty )^{s_0}$ and  $c_1,c_2>0$, $k\in
\Gamma^{+}_{c_1\theta^1+c_2\theta^2}$ implies 
$k\in\Gamma^{+}_{\theta^1}\cup\Gamma^{+}_{\theta^2}$ and similarly for $
\Gamma^{-}_{c_1\theta^1+c_2\theta^2}$, so
\begin{equation}\max_{k\in\Gamma^{+}_{c_1\theta^1+c_2\theta^2}}\rho_
k\leq\max_{k\in\Gamma^{+}_{\theta^1}}\rho_k\vee\max_{k\in\Gamma^{
+}_{\theta^2}}\rho_k\label{sump}\end{equation}
and
\begin{equation}\max_{k\in\Gamma^{-}_{c_1\theta^1+c_2\theta^2}}\rho_
k\leq\max_{k\in\Gamma^{-}_{\theta^1}}\rho_k\vee\max_{k\in\Gamma^{
-}_{\theta^2}}\rho_k.\label{suml}\end{equation}

Let $G$ be  a directed graph in which the nodes are identified 
with the species and a directed edge is drawn from $S_i$  
to $S_j$ if there is a
reaction that consumes $S_i$ and produces $S_j$.  
A subgraph $G_0\subset G$ is {\em strongly connected\/} if and only if for 
each pair $S_i,S_j\in G_0$, there is a directed path in $G_0$ beginning 
at $S_i$ and ending at $S_j$.  Single nodes are understood 
to form strongly connected subgraphs.  Recall that $G$ has a unique 
decomposition $G=\cup_jG_j$ into maximal strongly connected subgraphs.

The following lemma may significantly reduce the work 
needed to verify Condition \ref{cndz}.

\begin{lemma}\label{lembalcnd} 
Let $\theta\in [0,\infty )^{s_0}$, and fix $\gamma$.  Write 
\begin{equation}\theta =\sum_{j=1}^m\theta^j,\label{irrdec}\end{equation}
where $\mbox{\rm supp}(\theta^j)\subset G_j$ 
for some maximal strongly connected subgraph $G_j$ and $G_j\neq G_
i$ for 
$i\neq j$.  
 If 
Condition \ref{cndz} holds for each $\theta^j$, then it holds for 
$\theta$.  More specifically, 
if the balance equation (\ref{cnd1}) holds for 
each $\theta^j$, then the balance equation  holds for $\theta$, and if 
(\ref{cnd2}) holds for 
each $\theta^j$, then (\ref{cnd2})  holds for $\theta$.

Consequently, 
if Condition \ref{cndz} holds 
for each 
 $\theta\in [0,\infty )^{s_0}$ with 
support in some strongly connected subgraph, then Condition 
\ref{cndz}
holds for all $\theta\in [0,\infty )^{s_0}$;
 if (\ref{cnd1}) holds for each $\theta\in [0,\infty )^{s_0}$ with 
support in some strongly connected subgraph, then (\ref{cnd1}) 
holds for all $\theta\in [0,\infty )^{s_0}$; and 
if (\ref{cnd2}) holds for each $\theta\in [0,\infty )^{s_0}$ with 
support in some strongly connected subgraph, then (\ref{cnd2}) 
holds for all $\theta\in [0,\infty )^{s_0}$.
\end{lemma}
\begin{proof}
Assume that Condition \ref{cndz} holds for each $\theta^j$, 
$j=1,\ldots ,m$.
First, assume that $\Gamma^{+}_{\theta}\neq\emptyset$. Select  $l_
1\in\Gamma_{\theta}^{+}$ satisfying  
\begin{equation}\rho_{l_1}=\max_{k\in\Gamma_{\theta}^{+}}\rho_k.\label{lem15}\end{equation}
Since $\Gamma^{+}_{\theta}\subset\cup_j\Gamma^{+}_{\theta^j}$, 
there exists $j_1$ such that $l_1\in\Gamma^{+}_{\theta^{j_1}}$,
and using 
(\ref{lem15}), we have 
\begin{equation}\max_{k\in\Gamma^{+}_{\theta}}\rho_k=\rho_{l_1}\le\max_{
k\in\Gamma^{+}_{\theta^{j_1}}}\rho_k.\label{lem16}\end{equation}
We have three possible cases. 
First, 
if $\max_{k\in\Gamma_{\theta^{j_1}}^{+}}\rho_k\neq\max_{k\in\Gamma_{
\theta^{j_1}}^{-}}\rho_k$, then by (\ref{cnd2}),
there exists $i_1\in\mbox{\rm supp}(\theta^{j_1})$ such that
\begin{equation}\gamma +\max_{k\in\Gamma^{+}_{\theta^{j_1}}\cup\Gamma^{
-}_{\theta^{j_1}}}\rho_k\le\alpha_{i_1},\label{lem19}\end{equation}
and by (\ref{lem16}),
\begin{equation}\gamma +\max_{k\in\Gamma^{+}_{\theta}}\rho_k\le\alpha_{
i_1}\le\max_{i\in\mbox{\rm supp}(\theta )}\alpha_i.\label{lem110}\end{equation}

Second, if $\max_{k\in\Gamma_{\theta^{j_1}}^{+}}\rho_k=\max_{k\in
\Gamma_{\theta^{j_1}}^{-}}\rho_k\leq\max_{k\in\Gamma_{\theta}^{-}}
\rho_k$,
then by (\ref{lem16}), we obtain
\begin{equation}\max_{k\in\Gamma^{+}_{\theta}}\rho_k\le\max_{k\in
\Gamma^{+}_{\theta^{j_1}}}\rho_k=\max_{k\in\Gamma^{-}_{\theta^{j_
1}}}\rho_k\leq\max_{k\in\Gamma_{\theta}^{-}}\rho_k.\label{lem17}\end{equation}

Finally, if 
\begin{equation}\max_{k\in\Gamma_{\theta^{j_1}}^{+}}\rho_k=\max_{
k\in\Gamma_{\theta^{j_1}}^{-}}\rho_k>\max_{k\in\Gamma_{\theta}^{-}}
\rho_k\label{lem3a}\end{equation}
we select $l_2$ in $\Gamma_{\theta^{j_1}}^{-}$ with $\rho_{l_2}=\max_{
k\in\Gamma_{\theta^{j_1}}^{-}}\rho_k$.  The fact that 
$\rho_{l_2}>\max_{k\in\Gamma_{\theta}^{-}}\rho_k$ ensures the existence of
 $j_2$ such that $l_2\in\Gamma_{\theta^{j_2}}^{+}$. Then we have
\begin{equation}\max_{k\in\Gamma_{\theta^{j_1}}^{+}}\rho_k=\max_{
k\in\Gamma_{\theta^{j_1}}^{-}}\rho_k=\rho_{l_2}\le\max_{k\in\Gamma_{
\theta^{j_2}}^{+}}\rho_k.\label{lem111}\end{equation}

We recursively select $l_n$ and $j_n$ with $l_n\in\Gamma^{+}_{\theta^{
j_n}}$ such that
\[\max_{k\in\Gamma_{\theta^{j_{n-1}}}^{+}}\rho_k=\max_{k\in\Gamma_{
\theta^{j_{n-1}}}^{-}}\rho_k=\rho_{l_n}\le\max_{k\in\Gamma_{\theta^{
j_n}}^{+}}\rho_k\]
until we find $ $$l_n$ for which this is no longer 
possible. 
Since the $G_j$ are
maximal 
strongly connected subgraphs, 
there is no possibility that the same $\theta^j$
is selected more than once.  
Thus, 
the process will terminate for some $n$ and when it does
$\max_{k\in\Gamma_{\theta^{j_n}}^{+}}\rho_k\neq\max_{k\in\Gamma_{
\theta^{j_n}}^{-}}\rho_k$ and
\begin{equation}\gamma +\max_{k\in\Gamma^{+}_{\theta}}\rho_k\le\gamma 
+\max_{k\in\Gamma_{\theta^{j_n}}^{+}}\rho_k\leq\max_{i\in\mbox{\rm supp}
(\theta^{j_n})}\alpha_i\le\max_{i\in\mbox{\rm supp}(\theta )}\alpha_
i.\label{lem130a}\end{equation}
Consequently, we always have either
\begin{equation}\gamma +\max_{k\in\Gamma_{\theta}^{+}}\rho_k\le\max_{
i\in\mbox{\rm supp}(\theta )}\alpha_i\label{lem130}\end{equation}
or
\begin{equation}\max_{k\in\Gamma_{\theta}^{+}}\rho_k\le\max_{k\in
\Gamma_{\theta}^{-}}\rho_k.\label{lem131}\end{equation}

If $\Gamma^{-}_{\theta}\neq\emptyset$, interchanging $-$ and $+$, we see that either
\begin{equation}\gamma +\max_{k\in\Gamma_{\theta}^{-}}\rho_k\le\max_{
i\in\mbox{\rm supp}(\theta )}\alpha_i\label{lem120}\end{equation}
or
\begin{equation}\max_{k\in\Gamma_{\theta}^{-}}\rho_k\le\max_{k\in
\Gamma_{\theta}^{+}}\rho_k.\label{lem121}\end{equation}

Assume that both $\Gamma^{+}_{\theta}$ and $\Gamma^{-}_{\theta}$ are nonempty.
If both (\ref{lem131}) and (\ref{lem121}) hold, then 
(\ref{cnd1}) is satisfied.  If (\ref{lem130})  and 
(\ref{lem120}) hold, then taking the maximum of the left 
and right sides, (\ref{cnd2}) holds.  If (\ref{lem130}) and 
(\ref{lem121}) hold, then (\ref{cnd2}) holds and similarly for 
(\ref{lem131}) and (\ref{lem120}).

If (\ref{cnd1}) holds for all $\theta^j$, then the first and third 
cases above cannot hold so (\ref{lem17}) must hold giving 
(\ref{lem131}) and by the same argument (\ref{lem121}).  
Consequently, (\ref{cnd1}) must hold for $\theta$.  If (\ref{cnd2}) 
holds 
for all $\theta^j$, then the first case above holds giving 
(\ref{lem130}) and by the same argument (\ref{lem120}), so 
(\ref{cnd2}) must hold for $\theta$.

If $\Gamma^{+}_{\theta}=\emptyset$ and $\Gamma^{-}_{\theta}\neq\emptyset$, then (\ref{lem120}) must hold and 
(\ref{cnd2}) holds for $\theta$ and similarly with the $+$ and $-$ 
interchanged.  

If both $\Gamma^{+}_{\theta}$ and $\Gamma^{-}_{\theta}$ are empty, then (\ref{cnd1}) holds 
($-\infty =-\infty$).  In particular, $\theta\cdot\zeta_k=0$ for all $
\zeta_k$.
\end{proof}

The remaining lemmas in this section may be useful in 
verifying Condition \ref{cndz} for the cases that remain, 
that is, for $\theta$ with support in some strongly connected subgraph.

\begin{lemma}\label{zlem1} 
Fix $\gamma\in {\Bbb R}$, and suppose (\ref{cnd2}) holds for 
$\theta^1,\ldots ,\theta^m\in [0,\infty )^{s_0}$.  Then for $c_j>
0$, $j=1,\ldots ,m$, 
(\ref{cnd2}) holds for $\theta =\sum_{j=1}^mc_j\theta^j$.  
\end{lemma}

\begin{proof}
Since $\theta\cdot\zeta_k>0$ implies $c_j\theta^j\cdot\zeta_k>0$ for some $
j$ and 
$\theta\cdot\zeta_k<0$ implies $c_j\theta^j\cdot\zeta_k<0$ for some $
j$,
\[\max_{k\in\Gamma^{+}_{\theta}\cup\Gamma_{\theta}^{-}}\rho_k\leq\max_{
1\leq j\leq m}\max_{k\in\Gamma^{+}_{\theta^j}\cup\Gamma_{\theta^j}^{
-}}\rho_k\]
and there exists $j$   such that
\[\gamma\leq\max_{i:\theta_i^j>0}\alpha_i-\max_{k\in\Gamma^{+}_{\theta^
j}\cup\Gamma_{\theta^j}^{-}}\rho_k\leq\max_{i:\theta_i>0}\alpha_i
-\max_{k\in\Gamma^{+}_{\theta}\cup\Gamma_{\theta}^{-}}\rho_k.\]
\end{proof}

\begin{lemma}\label{zlem4}
For $\theta^1,\theta^2\in [0,\infty )^{s_0}$, suppose that
\begin{equation}\max_{k\in\Gamma^{-}_{\theta^1}}\rho_k=\max_{k\in
\Gamma^{+}_{\theta^1}}\rho_k>\max_{k\in\Gamma^{+}_{\theta^2}\cup\Gamma_{
\theta^2}^{-}}\rho_k.\label{z4a}\end{equation}
Then for $c_1,c_2>0$, (\ref{cnd1}) holds for $c_1\theta^1+c_2\theta^
2$.
\end{lemma}

\begin{proof}
If $l\in\Gamma^{+}_{\theta^1}$ and $\rho_l=\max_{k\in\Gamma^{+}_{
\theta^1}}\rho_k$, then by (\ref{z4a}), 
$l\notin\Gamma^{-}_{\theta^2}$.  Consequently, $l\in\Gamma^{+}_{c_
1\theta^1+c_2\theta^2}$ and by (\ref{sump}), we 
must have 
\[\max_{k\in\Gamma^{+}_{c_1\theta^1+c_2\theta^2}}\rho_k=\max_{k\in
\Gamma^{+}_{\theta^1}}\rho_k.\]
By the same argument, 
\[\max_{k\in\Gamma^{-}_{c_1\theta^1+c_2\theta^2}}\rho_k=\max_{k\in
\Gamma^{-}_{\theta^1}}\rho_k,\]
and it follows that (\ref{cnd1}) holds for $c_1\theta^1+c_2\theta^
2$.
\end{proof}

\begin{lemma}
Fix $\gamma$, and suppose that (\ref{cnd1}) holds for $\theta^1$ and (\ref{cnd2}) 
for $\theta^2$.  Then for $c_1,c_2>0$, Condition \ref{cndz} holds for 
$c_1\theta^1+c_2\theta^2$.
\end{lemma}

\begin{proof}
If
\begin{equation}\max_{k\in\Gamma_{\theta^1}^{-}}\rho_k=\max_{k\in
\Gamma_{\theta^1}^{+}}\rho_k>\max_{k\in\Gamma_{\theta^2}^{-}}\rho_
k\vee\max_{k\in\Gamma^{+}_{\theta^2}}\rho_k,\label{zcase1}\end{equation}
then Lemma \ref{zlem4} implies $c_1\theta^1+c_2\theta^2$ 
satisfies (\ref{cnd1}), so assume that
\begin{equation}\max_{k\in\Gamma_{\theta^1}^{-}}\rho_k=\max_{k\in
\Gamma_{\theta^1}^{+}}\rho_k\leq\max_{k\in\Gamma_{\theta^2}^{-}}\rho_
k\vee\max_{k\in\Gamma^{+}_{\theta^2}}\rho_k.\label{zcase2}\end{equation}
Then
\[\max_{k\in\Gamma^{-}_{c_1\theta^1+c_2\theta^2}}\rho_k\leq\max_{
k\in\Gamma_{\theta^1}^{-}}\rho_k\vee\max_{k\in\Gamma_{\theta^2}^{
-}}\rho_k\leq\max_{k\in\Gamma_{\theta^2}^{-}}\rho_k\vee\max_{k\in
\Gamma^{+}_{\theta^2}}\rho_k\]
and
\[\max_{k\in\Gamma^{+}_{c_1\theta^1+c_2\theta^2}}\rho_k\leq\max_{
k\in\Gamma_{\theta^1}^{+}}\rho_k\vee\max_{k\in\Gamma_{\theta^2}^{
+}}\rho_k\leq\max_{k\in\Gamma_{\theta^2}^{-}}\rho_k\vee\max_{k\in
\Gamma^{+}_{\theta^2}}\rho_k,\]
so
\[\max_{k\in\Gamma^{-}_{c_1\theta^1+c_2\theta^2}}\rho_k\vee\max_{
k\in\Gamma^{+}_{c_1\theta^1+c_2\theta^2}}\rho_k\leq\max_{k\in\Gamma_{
\theta^2}^{-}}\rho_k\vee\max_{k\in\Gamma^{+}_{\theta^2}}\rho_k,\]
and since $\mbox{\rm supp}(c_1\theta^1+c_2\theta^2)\supset\mbox{\rm supp}
(\theta^2)$, (\ref{cnd2}) for $\theta^2$ 
implies (\ref{cnd2}) for $c_1\theta^1+c_2\theta^2$.
\end{proof}

If Condition \ref{cndz} holds for $\theta^1$  and $\theta^2$ and $
c_1,c_2>0$, then the 
previous lemmas imply 
Condition 
\ref{cndz} holds for $c_1\theta^1+c_2\theta^2$ 
 except in one possible situation, that is,
\begin{equation}\max_{k\in\Gamma_{\theta^1}^{-}}\rho_k=\max_{k\in
\Gamma_{\theta^1}^{+}}\rho_k=\max_{k\in\Gamma_{\theta^2}^{-}}\rho_
k=\max_{k\in\Gamma_{\theta^2}^{+}}\rho_k.\label{z1c}\end{equation}
Since the species balance condition does not imply 
Condition \ref{cndz} for $\theta =(1,1)$ for the system $(Z_1^N,Z_
2^N)$ 
given by (\ref{unbal}), some additional condition must be 
required to be able to conclude Condition \ref{cndz}\ holds 
for 
$c_1\theta^1+c_2\theta^2$ when (\ref{z1c}) 
holds.  The following lemmas give such conditions.

\begin{lemma}\label{zlem2}
Fix $\gamma\in {\Bbb R}$, and suppose  that Condition \ref{cndz} holds for 
$\theta^1,\theta^2\in [0,\infty )^{s_0}$.  If $\Gamma_{\theta^1}^{
+}\cap\Gamma_{\theta^2}^{-}=\emptyset$ or $\Gamma_{\theta^1}^{-}\cap
\Gamma_{\theta^2}^{+}=\emptyset$ and $c_1,c_2>0$,
then Condition \ref{cndz} holds for 
$c_1\theta^1+c_2\theta^2$.  

If Condition \ref{cnd1}\ holds for $\theta_1$ and $\theta_2$, 
 $\Gamma_{\theta^1}^{+}\cap\Gamma_{\theta^2}^{-}=\emptyset$ or $\Gamma_{
\theta^1}^{-}\cap\Gamma_{\theta^2}^{+}=\emptyset$, and $c_1,c_2>0$,
then Condition \ref{cnd1} holds for 
$c_1\theta^1+c_2\theta^2$.  
\end{lemma}

\begin{remark}
 If no reaction that consumes a 
species in the support of $\theta^1$  produces a 
species in the support of $\theta^2$, then $\Gamma_{\theta^1}^{-}
\cap\Gamma_{\theta^2}^{+}=\emptyset$.
  That condition is, of 
course, equivalent to the requirement that a reaction 
that produces a species in the support of $\theta^2$ does not 
consume a species in the support of $\theta^1$.
\end{remark}

\begin{proof}
As noted, the previous lemmas cover all possible situations except 
in the case that (\ref{z1c}) holds.
Suppose $\Gamma_{\theta^1}^{-}\cap\Gamma_{\theta^2}^{+}=\emptyset$.  If $
\theta^1\cdot\zeta_k<0$, then $\theta^2\cdot\zeta_k\leq 0$ and 
$(c_1\theta^1+c_2\theta^2)\cdot\zeta_k<0$, and if $(c_1\theta^1+c_
2\theta^2)\cdot\zeta_k<0$, then either 
$\theta^1\cdot\zeta_k<0$ or $\theta^2\cdot\zeta_k<0$, so
\begin{equation}\max_{k\in\Gamma_{\theta^1}^{-}}\rho_k\leq\max_{k
\in\Gamma^{-}_{c_1\theta^1+c_2\theta^2}}\rho_k\leq\max_{k\in\Gamma_{
\theta^1}^{-}}\rho_k\vee\max_{k\in\Gamma_{\theta^2}^{-}}\rho_k.\label{z1a}\end{equation}
Similarly, noting that $\theta^2\cdot\zeta_k>0$ implies $\theta^1
\cdot\zeta_k\geq 0$,
\begin{equation}\max_{k\in\Gamma_{\theta^2}^{+}}\rho_k\leq\max_{k
\in\Gamma^{+}_{c_1\theta^1+c_2\theta^2}}\rho_k\leq\max_{k\in\Gamma_{
\theta^1}^{+}}\rho_k\vee\max_{k\in\Gamma_{\theta^2}^{+}}\rho_k.\label{z1b}\end{equation}
But (\ref{z1c}) implies equality holds throughout (\ref{z1a}) and 
(\ref{z1b}) and  (\ref{cnd1}) holds 
for $c_1\theta^1+c_2\theta^2$.  
\end{proof}

\begin{lemma}
Suppose (\ref{cnd1}) holds for $\theta^1$ and $\theta^2$ and for
$\theta^1-\frac {\theta^1\cdot\zeta_k}{\theta^2\cdot\zeta_k}\theta^
2$ for all  $k\in (\Gamma_{\theta^1}^{+}\cap\Gamma_{\theta^2}^{-}
)\cup (\Gamma_{\theta^1}^{-}\cap\Gamma_{\theta^2}^{+})$.  (Note that 
$-\frac {\theta^1\cdot\zeta_k}{\theta^2\cdot\zeta_k}>0$.)  Then (\ref{cnd1}) holds for $
c_1\theta^1+c_2\theta^2$ for all 
$c_1,c_2>0$.
\end{lemma}

\begin{proof}
By Lemma \ref{zlem4}, we can restrict our attention to 
the case (\ref{z1c}), and
it is enough to consider $\theta^1+c\theta^2$ for $c>0$.  Note that for 
$c$ sufficiently small, $\Gamma_{\theta^1+c\theta^2}^{+}\supset\Gamma_{
\theta^1}^{+}$  and $\Gamma_{\theta^1+c\theta^2}^{-}\supset\Gamma_{
\theta^1}^{-}$ and
\[\max_{k\in\Gamma_{\theta^1+c\theta^2}^{+}}\rho_k=\max_{k\in\Gamma_{
\theta^1}^{+}}\rho_k=\max_{k\in\Gamma_{\theta^1+c\theta^2}^{-}}\rho_
k=\max_{k\in\Gamma_{\theta^1}^{-}}\rho_k.\]
Let 
\[c_0=\inf\{c:\max_{k\in\Gamma_{\theta^1+c\theta^2}^{+}}\rho_k\neq\max_{
k\in\Gamma_{\theta^1}^{+}}\rho_k\mbox{\rm \ or }\max_{k\in\Gamma_{
\theta^1+c\theta^2}^{-}}\rho_k\neq\max_{k\in\Gamma_{\theta^1}^{-}}
\rho_k\},\]
and note that for $0<c<c_0$, (\ref{cnd1}) holds.
If $c_0<\infty$, then $c_0=-\frac {\theta^1\cdot\zeta_k}{\theta^2
\cdot\zeta_k}>0$  for some $k$, and by the 
assumptions of the lemma 
\begin{equation}\max_{k\in\Gamma_{\theta^1+c_0\theta^2}^{+}}\rho_
k=\max_{k\in\Gamma_{\theta^1+c_0\theta^2}^{-}}\rho_k<\max_{k\in\Gamma_{
\theta^1}^{+}}\rho_k=\max_{k\in\Gamma_{\theta^1}^{-}}\rho_k.\label{sgch}\end{equation}
But (\ref{sgch}) can hold only if there exists $l^{+}\in\Gamma^{+}_{
\theta_1}$ such 
that $\rho_{l^{+}}=\max_{k\in\Gamma_{\theta^1}^{+}}\rho_k$ and $c_
0=-\frac {\theta^1\cdot\zeta_{l^{+}}}{\theta^2\cdot\zeta_{l^{+}}}$ and  $
l^{-}\in\Gamma^{-}_{\theta_1}$ 
such that
$\rho_{l^{-}}=\max_{k\in\Gamma_{\theta^1}^{-}}\rho_k$ and $c_0=-\frac {
\theta^1\cdot\zeta_{l^{-}}}{\theta^2\cdot\zeta_{l^{-}}}$.  Then, for $
c>c_0$, 
$l^{+}\in\Gamma^{-}_{\theta^1+c\theta^2}$ and $l^{-}\in\Gamma^{-}_{
\theta^1+c\theta^2}$, and the lemma follows. 
\end{proof}

\setcounter{equation}{0}

\section{Derivation of limiting models}\label{sectlim}
As can be seen from the examples, derivation of the 
limiting models can frequently be carried out by 
straightforward analysis of the stochastic equations.  
The results of this section may be harder to apply than 
direct analysis, but they should give added confidence 
that the limits hold in great generality for complex 
models.

We assume throughout this section that 
$\lim_{N\rightarrow\infty}Z^{N,\gamma}_i(0)$ exists and is positive for all $
i$.
If 
\begin{equation}\gamma =r_1\equiv\min_i\gamma_i=\min_i(\alpha_i-\max_{
k\in\Gamma^{+}_i\cup\Gamma_i^{-}}(\beta_k+\nu_k\cdot\alpha )),\label{scl1}\end{equation}
then $\lim_{N\rightarrow\infty}Z^{N,\gamma}$ exists, at least on some interval 
$[0,\tau_{\infty})$ with $\tau_{\infty}>0$, 
and is easy to calculate since on any time interval over 
which $\sup_{t\leq T}|Z^{N,\gamma}(t)|<\infty$,
each term 
\[N^{-\alpha_i}Y_k(\int_0^tN^{\gamma +\rho_k}\lambda_k(Z^{N,\gamma}
(s))ds)\]
either converges to zero (if $\alpha_i>\gamma +\rho_k$), is 
dependent on $N$ only through $Z^{N,\gamma}$ (if 
$\alpha_i=\gamma +\rho_k=0$), or is asymptotic to
\[\int_0^t\lambda_k(Z^{N,\gamma}(s))ds\]
(if $\alpha_i=\gamma +\rho_k>0$), since
\[\lim_{N\rightarrow\infty}\sup_{u\leq u_0}\left|N^{-\alpha_i}Y_k
(N^{\alpha_i}u)-u\right|=0,\quad u_0>0.\]
The caveat regarding the interval $[0,\tau_{\infty})$ reflects the 
fact that we have not ruled out ``reaction'' networks of the form 
$2S_1\rightarrow 3S_1$, $S_1\rightarrow\emptyset$ which would be modeled by
\[X_1(t)=X_1(0)+Y_1(\kappa_1\int_0^tX_1(s)(X_1(s)-1)ds)-Y_2(\kappa_
2\int_0^tX_1(s)ds)\]
and has positive probability of exploding in finite time, if 
$X_1(0)>1$.

\begin{theorem}
Let $\Gamma_i^{\gamma}=\{k:\gamma +\rho_k=\alpha_i\}$.
For $r_1$ defined by (\ref{scl1}), $Z^{N,r_1}\Rightarrow Z^{r_1}$ on $
[0,\tau_{\infty})$, where 
if $\alpha_i>0$,
\[Z^{r_1}_i(t)=Z_i(0)+\sum_{k\in\Gamma_i^{r_1}}\int_0^t\lambda_k(
Z^{r_1}(s))ds(\nu'_{ik}-\nu_{ik}),\]
if $\alpha_i=0$,
\[Z^{r_1}_i(t)=Z_i(0)+\sum_{k\in\Gamma_i^{r_1}}Y_k(\int_0^t\lambda_
k(Z^{r_1}(s))ds)(\nu'_{ik}-\nu_{ik}),\]
and 
\[\tau_{\infty}=\lim_{c\rightarrow\infty}\tau_c\equiv\inf\{t:\sup_{
s\leq t}|Z^{r_1}(s)|\geq c\}.\]
\end{theorem}

\begin{remark}
By $Z^{N,r_1}\Rightarrow Z^{r_1}$ on $[0,\tau_{\infty})$,  we mean that there exist $
\tau^{N,n}$ 
and $\tau^n$ such that $(Z^{N,r_1}(\cdot\wedge\tau^{N,n}),\tau^{N
,n})\Rightarrow (Z^{r_1}(\cdot\wedge\tau^n),\tau^n)$ 
and $\lim_{n\rightarrow\infty}\tau^n=\tau_{\infty}$. 
\end{remark}

\begin{proof}
Let $\tau_{N,c}=\inf\{t:\sup_{s\leq t}|Z^{N,r_1}(s)|\geq c\}$.  The relative 
compactness of $\{Z^{N,r_1}(\cdot\wedge\tau_{N,c})\}$ follows from the 
uniform boundedness of $\lambda_k(Z^{N,r_1}(\cdot\wedge\tau_{N,c}
))$.  Then 
$(Z^{N,r_1}(\cdot\wedge\tau_{N,c}),\tau_{N,c})\Rightarrow (Z^{r_1}
(\cdot\wedge\tau_c),\tau_c)$ at least for all 
but countably many $c$.
\end{proof}

Note that $\gamma_{\theta}\geq\min_{i:\theta_i>0}\gamma_i$, so $r_
1=\min_{\theta\in [0,\infty )^{s_0}}\gamma_{\theta}$, and
Condition \ref{cndz} always holds for $\gamma =r_1$.
If Condition \ref{cndz}  holds for some $\gamma >r_1$, then the 
balance equality (\ref{cnd1}) must hold for all 
$\theta\in [0,\infty )^{s_0}$ with $\gamma_{\theta}=r_1$.

Let $\alpha_{\theta}=\max_{i:\theta_i>0}\alpha_i$, and define
\[Z_{\theta}^{N,\gamma}(t)=N^{-\alpha_{\theta}}\theta\cdot\Lambda_
N^{-1}Z^{N,\gamma}(t)=N^{-\alpha_{\theta}}\sum_{i=1}^{s_0}\theta_
iX^N_i(N^{\gamma}t).\]
As noted earlier, the natural time scale for $Z_{\theta}^{N,\gamma}$ is
\[\gamma_{\theta}=\alpha_{\theta}-\max_{k\in\Gamma^{+}_{\theta}\cup
\Gamma^{-}_{\theta}}\rho_k.\]

Let  ${\Bbb L}_1$ be the space spanned 
by ${\Bbb S}_1=\{e_i:\Gamma_i^{r_1}\neq\emptyset \}$, and let
$\Pi_1$ be the projection of ${\Bbb R}^{s_0}$ onto ${\Bbb L}_1$.
Let 
\[{\Bbb K}_2=\{\theta\in [0,\infty )^{s_0}:\theta\cdot\Pi_1\zeta_
k=0\forall k\in\cup_i\Gamma_i^{r_1}\}\]
and ${\Bbb L}_2=\mbox{\rm span}\,{\Bbb K}_2$, and let $\Pi_2$ be the projection onto $
{\Bbb L}_2$.  
Of course, ${\Bbb K}_2$ contains ${\Bbb S}_2=\{e_i:\Gamma_i^{r_1}
=\emptyset \}$,  but as in the 
example of Section \ref{sectaux}, it may be larger.  The 
projections
 $\Pi_1$ and $\Pi_2$ are not necessarily orthogonal, but 
for any $x\in {\Bbb R}^{s_0}$, $x-\Pi_2x\in {\Bbb L}_1$.

\begin{lemma}
For each $x\in {\Bbb R}^{s_0}$, $x-\Pi_2x\in {\Bbb L}_1$.
\end{lemma}

\begin{proof}
Note that ${\Bbb L}_1=\{x\in {\Bbb R}^{s_0}:e_i\cdot x=0,\forall 
e_i\in {\Bbb S}_2\}$ and that for 
$e_i\in {\Bbb S}_2$, $e_i\cdot\Pi_2x=e_i\cdot x$.  Consequently, for $
e_i\in {\Bbb S}_2$, 
$e_i\cdot (x-\Pi_2x)=0$ and $x-\Pi_2x\in {\Bbb L}_1$.
\end{proof}

\begin{lemma}$ $ 
If $\theta\in {\Bbb K}_2$ and $\theta\cdot\zeta_l\neq 0$ 
for some $l\in\Gamma_i^{r_1}$, then 
$\alpha_{\theta}>\alpha_i$. $ $

Let
\[r_2=\min_{\theta\in {\Bbb K}_2}\gamma_{\theta}=\min_{\theta\in 
{\Bbb K}_2}\{\alpha_{\theta}-\max_{k\in\Gamma^{+}_{\theta}\cup\Gamma_{
\theta}^{-}}\rho_k\}.\]
Then $r_2>r_1$.
\end{lemma}

\begin{proof}
For $\theta\in {\Bbb K}_2$, if $\theta\cdot\zeta_l\neq 0$, 
then
$\theta_j>0$ for some $j$ such that $e_j\notin {\Bbb L}_1$.  Since
\[r_1\leq\alpha_j-\max_{k\in\Gamma^{+}_j\cup\Gamma_j^{-}}\rho_k\leq
\alpha_j-\rho_l\]
and $\Gamma_j^{r_1}=\{k:\alpha_j-\rho_k=r_1\}=\emptyset$, 
\[r_1=\alpha_i-\rho_l<\alpha_j-\rho_l,\]
and $\alpha_{\theta}\geq\alpha_j>\alpha_i$.

Let $\theta\in {\Bbb K}_2$ satisfy $\gamma_{\theta}=r_2$.  Then there exists 
\[l\in\Gamma^{+}_{\theta}\cup\Gamma_{\theta}^{-}\subset\cup_{j:\theta_
j>0}\Gamma^{+}_j\cup\Gamma_j^{-}\]
and $\theta_j>0$ such that $l\in\Gamma^{+}_j\cup\Gamma_j^{-}$ and
\[r_2=\alpha_{\theta}-\rho_l\geq\alpha_j-\rho_l\geq\gamma_j\geq r_
1.\]
If $\alpha_j-\rho_l>\gamma_j$ or $\gamma_j>r_1$, then $r_2>r_1$.  If 
$\alpha_j-\rho_l=\gamma_j=r_1$, then $l\in\Gamma^{r_1}_j$ and since $
\theta\cdot\zeta_l\neq 0$, $\alpha_{\theta}>\alpha_j$ 
and $r_2>r_1$.  
\end{proof}

Unfortunately, while $r_2$ can naturally be viewed as the 
second time scale, we cannot guarantee a priori that the 
system will converge to a nondegenerate model
on that time scale.  For example, 
consider the network
\begin{center}
\[\emptyset\rightarrow S_1\quad\emptyset\rightarrow S_2\quad\emptyset
\rightarrow S_3\]
\[S_1+S_2\rightarrow\emptyset\qquad S_1+S_3\rightarrow\emptyset\]
\end{center}
and assume that the parameters scale so that
\begin{eqnarray*}
X_1(t)&=&X_1(0)+Y_1(\kappa_1t)-Y_2(\kappa_2\int_0^tX_1(s)X_2(s))d
s)-Y_5(\kappa_5N^{-1}\int_0^tX_1(s)X_3(s)ds)\\
X_2(t)&=&X_2(0)+Y_3(\kappa_3t)-Y_2(\kappa_2\int_0^tX_1(s)X_2(s))d
s)\\
X_3(t)&=&X_3(0)+Y_4(\kappa_4N^{-1}t)-Y_5(\kappa_5N^{-1}\int_0^tX_
1(s)X_3(s)ds).\end{eqnarray*}
Then Condition \ref{cnd1}\ is satisfied for all $\theta$, $r_1=0$, and 
$r_2=1$.  But if $\kappa_1>\kappa_3$, $X_1(Nt)\rightarrow\infty$ and $
X_3(Nt)\rightarrow 0$ for all 
$t>0$.

The problem is that even though the balance equations 
are satisfied for the fast subnetwork $(X_1,X_2)$, the 
subnetwork is not stable.  Consequently, to guarantee 
convergence on the second time scale, we need some 
additional condition to ensure stability for the fast 
subnetwork so that the influence of the fast components 
can be averaged in the system on the second time scale.

Of course, with reference to (\ref{chtave}) and 
(\ref{intave}), it is frequently possible to
verify convergence without any special techniques, but 
we will outline a more systematic approach.  

We assume the following condition on the scaling.

\begin{condition}\label{scl}
For each $N$, $\Lambda_N\Pi_2\Lambda_N^{-1}=\Pi_2$.
\end{condition}

Let ${\Bbb L}^{\alpha}$ be the span of $\{e_i:\alpha_i=\alpha \}$.  Then Condition 
\ref{scl}\  is equivalent to the requirement that 
$\Pi_2:{\Bbb L}^{\alpha}\rightarrow {\Bbb L}^{\alpha}$.

Define the 
occupation measure on ${\Bbb L}_1\times [0,\infty )$ by
\[V^{N,r_2}_1(C\times [0,t])=\int_0^t{\bf 1}_C((I-\Pi_2)Z^{N,r_2}
(s))ds.\]
Assume that 
\begin{equation}V^{N,r_2}_1\Rightarrow V_1\label{occnv}\end{equation}
in the sense that 
\[\int_{{\Bbb L}_1\times [0,t]}f(x)V_1^{N,r_2}(dx\times ds)\Rightarrow
\int_{{\Bbb L}_1\times [0,t]}f(x)V_1(dx\times ds)\]
for all $f\in C_b({\Bbb L}_1)$ and all $t>0$.  This requirement is 
essentially an ergodicity assumption on the fast 
subsystem.

Define $\tau_q^N=\inf\{t:|\Pi_2Z^{N,r_2}(t)|\geq q\}$.
For $\theta\in {\Bbb K}_2$, define $\Gamma^{r_2}_{\theta}=\{k:r_2
+\rho_k=\alpha_{\theta}\}$ and
\[h_{q,\theta}(y)=\sup_{x\in {\Bbb L}_2,|x|\leq q}\sum_{k\in\Gamma_{
\theta}^{r_2}}|\theta\cdot\zeta_k|\lambda_k(x+y),\]
and assume that 
for $q>0$, 
$\psi_{q,\theta}:[0,\infty )\rightarrow [0,\infty )$ satisfies $\lim_{
r\rightarrow\infty}r^{-1}\psi_{q,\theta}(r)=\infty$ and
\begin{equation}\{\int_{{\Bbb L}_1\times [0,t\wedge\tau_q^N]}\psi_{
q,\theta}(h_{q,\theta}(y))V_1^{N,r_2}(dy\times ds)\}\label{uint}\end{equation}
is stochastically bounded.  In addition, assume
\[\sum_{k:r_2+\rho_k<\alpha_{\theta}}|\theta\cdot\zeta_k|N^{r_2+\rho_
k-\alpha_{\theta}}\int_{{\Bbb L}_1\times [0,t\wedge\tau_q^N]}\lambda_
k(\Pi_2Z^{N,r_2}(s)+y)V^{N,r_2}_1(dy\times ds)\rightarrow 0.\]
Then for all but 
countably many $q$,
 at least along a subsequence, 
$\Pi_2Z^{N,r_2}(\cdot\wedge\tau^N_q)$  converges in distribution to a process 
$Z^{r_2}_2(\cdot\wedge\tau_q)$, and if $\rho_k+r_2=\alpha_{\theta}$, by Lemma \ref{wkconv}, 
\begin{equation}\int_0^{t\wedge\tau^N_q}\lambda_k(Z^{N,r_2}(s))ds
\Rightarrow\int_{{\Bbb L}_1\times [0,t\wedge\tau_q]}\lambda_k(Z^{
r_2}_2(s)+y)V_1(dy\times ds).\label{tcconv}\end{equation}

\begin{theorem}\label{avecnv} 
Define $D^{\alpha}=\mbox{\rm diag}(\ldots {\bf 1}_{\{\alpha_i=\alpha 
\}}\ldots )$.  Under the above 
assumptions, there exists a ${\Bbb L}_2$-valued process $Z_2^{r_2}$ and a 
random variable $\tau_{\infty}>0$ such  $\Pi_2Z^{N,r_2}$ converges in 
distribution to $Z_2^{r_2}$ on $[0,\tau_{\infty})$.  For $\theta\in 
{\Bbb K}_2$ with $\alpha_{\theta}=0$, 
\[\theta\cdot Z_2^{r_2}(t)=\theta\cdot Z_2^{r_2}(0)+\sum_{k\in\Gamma^{
r_2}_{\theta}}(\theta\cdot\zeta_k)Y_k(\int_{{\Bbb L}_1\times [0,t
]}\lambda_k(Z_2^{r_2}(s)+y)V_1(dy\times ds))\]
and for $\theta\in {\Bbb K}_2$ with $\alpha_{\theta}>0$,
\[\theta\cdot D^{\alpha_{\theta}}Z_2^{r_2}(t)=\theta\cdot D^{\alpha_{
\theta}}Z_2^{r_2}(0)+\sum_{k\in\Gamma^{r_2}_{\theta}}(\theta\cdot 
D^{\alpha_{\theta}}\zeta_k)\int_{{\Bbb L}_1\times [0,t]}\lambda_k
(Z_2^{r_2}(s)+y)V_1(dy\times ds),\]
for $t\in [0,\tau_{\infty})$.

In particular, $Z_{\theta}^{N,r_2}\Rightarrow\theta\cdot D^{\alpha_{
\theta}}Z_2^{r_2}$.
\end{theorem}

\begin{remark}
The statement of this theorem is somewhat misleading.  
Given $V_1$, $Z_2^{r_2}$ is uniquely determined.  However, as we 
will see in the next section, typically $V_1$ depends on $Z_2^{r_
2}$.  
There we will give conditions under which the sequence 
of pairs $\{(V_1^{N,r_2},Z^{N,r_2})\}$ is relatively compact.  Then any 
limit point $(V_1,Z^{r_2}_2)$ will satisfy the equations given by 
the present theorem, but it will still be necessary to 
show that the pair is uniquely determined.  
\end{remark}
\[\]
\begin{proof}
As for the first time-scale, stopping the process at 
$ $
\[\tau_q^N=\inf\{t:|\Pi_2Z^{N,r_2}(t)|\geq q\}\]
ensures that $\{\Pi_2Z^{N,r_2}(\cdot\wedge\tau_q^N)\}$ is relatively compact, and 
(\ref{tcconv})
ensures that any limit 
process satisfies the stochastic equations.   Uniqueness 
for the limiting system then follows by the smoothness 
of the $\lambda_k$.
\end{proof}

\setcounter{equation}{0}

\section{Averaging}
\label{sectave}
Stochastic averaging methods go back at least to 
\cite{Kha66a,Kha66b}.  In this section we summarize the 
approach taken in \cite{Kur92}.  See that article for 
additional detail and references.

Recall that $\Lambda_N=\mbox{\rm diag}(N^{-\alpha_1},\ldots ,N^{-
\alpha_{s_0}})$, $\rho_k=\beta_k+\nu_k\cdot\alpha$, and 
$\zeta_k=\nu_k'-\nu_k$.  The generator for $Z^{N,0}$ is 
\[{\Bbb B}_Nf(z)=\sum_kN^{\rho_k}\lambda_k(z)(f(z+\Lambda_N\zeta_
k)-f(z)).\]
Another way of characterizing $r_1$ is as the largest $\gamma$ 
(possibly negative) such that $\lim_{N\rightarrow\infty}N^{\gamma}
{\Bbb B}_Nf(z)$ exists 
for each $f\in C^2_c({\Bbb R}^m)$ and $z\in {\Bbb R}^m$.  Define 
$\Gamma_{\alpha}^{r_1}=\{k:r_1+\rho_k=\alpha \}$  
and set $D^{\alpha}=\mbox{\rm diag}(\ldots {\bf 1}_{\{\alpha_i=\alpha 
\}}\ldots )$. 
Then
\begin{eqnarray*}
{\Bbb C}_0f(z)&\equiv&\lim_{N\rightarrow\infty}N^{r_1}{\Bbb B}_Nf
(z)\\
&=&\sum_{k\in\Gamma^{r_1}_0}\lambda_k(z)(f(z+\Lambda^0\zeta_k)-f(
z))+(\sum_{\alpha >0}\sum_{k\in\Gamma^{r_1}_{\alpha}}\lambda_k(z)
D^{\alpha}\zeta_k)\cdot\nabla f(z),\end{eqnarray*}
which is the generator for the limit of the system on 
the first time scale.  The state space for the limit 
process is ${\Bbb E}=\prod_{i=1}^{s_0}{\Bbb E}_i$, where ${\Bbb E}_
i={\Bbb N}$ if $\alpha_i=0$ and 
${\Bbb E}_i=[0,\infty )$ if $\alpha_i>0$. 

 Note that if $k\in\Gamma_{\alpha}^{r_1}$, 
then $ $$D^{\alpha}\zeta_k=\Pi_1\zeta_k$, and by the definition of $
{\Bbb L}_2$, $\Pi_2\Pi_1\zeta_k=0$.  
Consequently, for $z\in\Pi_2{\Bbb E}$ and 
\[{\Bbb E}_z=\{y\in {\Bbb L}_1:y=(I-\Pi_2)x,\Pi_2x=z,x\in {\Bbb E}
\},\]
\[{\Bbb C}^zf(y)\equiv {\Bbb C}_0f(z+y)\]
defines a generator with state space ${\Bbb E}_z$.

As before, define
\[V^{N,r_2}_1(C\times [0,t])=\int_0^t{\bf 1}_C((I-\Pi_2)Z^{N,r_2}
(s))ds,\]
and observe that
\begin{eqnarray*}
M_f^N(t)&=&f(Z^{N,r_2}(t))-f(Z^{N,r_2}(0))-\int_0^tN^{r_2}{\Bbb B}_
Nf(Z^{N,r_2}(s))ds\\
&=&f(Z^{N,r_2}(t))-f(Z^{N,r_2}(0))-\int_{{\Bbb L}_1\times [0,t]}N^{
r_2}{\Bbb B}_Nf(\Pi_2Z^{N,r_2}(s)+y)V_1^{N,r_2}(dy\times ds)\end{eqnarray*}
is a martingale.  Since $f$ and $N^{r_1}{\Bbb B}_Nf$ are bounded by 
constants, $N^{r_1-r_2}M_f^N$ is bounded by a constant on any 
bounded time interval.  It follows that $\{N^{r_1-r_2}M_f^N\}$ is 
relatively compact, any limit point is a martingale with 
initial value zero, and any limit point is Lipschitz 
continuous with Lipschitz constant $\sup_z|{\Bbb C}_0f(z)|$.  Since 
any continuous martingale with finite variation paths is 
constant, it follows that the limit must be zero.  
Combining these observations with those of the previous 
section, we have the following theorem.  

\begin{theorem}
Assume Condition \ref{scl}.  Suppose  (\ref{occnv}) holds 
and that  (\ref{uint}) is stochastically bounded. Let 
$Z^{r_2}_2$ and $\tau_{\infty}$ be as in the conclusion of Theorem 
\ref{avecnv}.  Then for all $f\in C^2_c({\Bbb R}^{s_0})$,
\[\int_{{\Bbb L}_1\times [0,\tau_{\infty})}{\Bbb C}_0f(Z_2^{r_2}(
s)+y)V_1(dy\times ds)=\int_{{\Bbb L}_1\times [0,\tau_{\infty})}{\Bbb C}^{
Z_2^{r_2}(s)}f(y)V_1(dy\times ds)=0.\]

If for each $z\in\Pi_2{\Bbb E}$, $\pi_z$ is the unique stationary 
distribution for ${\Bbb C}^z$, then
\[V_1(dy\times ds)=\pi_{Z_2^{r_2}(s)}(dy)ds,\]
and the system of equations in Theorem \ref{avecnv} 
becomes
\[\theta\cdot Z_2^{r_2}(t)=\theta\cdot Z_2^{r_2}(0)+\sum_{k\in\Gamma^{
r_2}_{\theta}}(\theta\cdot\zeta_k)Y_k(\int_0^t\int_{{\Bbb L}_1}\lambda_
k(Z_2^{r_2}(s)+y)\pi^{Z_2^{r_2}(s)}(dy)ds),\]
for $\theta\in {\Bbb K}_2$ with $\alpha_{\theta}=0$, and 
\[\theta\cdot D^{\alpha_{\theta}}Z_2^{r_2}(t)=\theta\cdot D^{\alpha_{
\theta}}Z_2^{r_2}(0)+\sum_{k\in\Gamma^{r_2}_{\theta}}(\theta\cdot 
D^{\alpha_{\theta}}\zeta_k)\int_0^t\int_{{\Bbb L}_1}\lambda_k(Z_2^{
r_2}(s)+y)\pi^{Z_2^{r_2}(s)}(dy)ds,\]
for $\theta\in {\Bbb K}_2$ with $\alpha_{\theta}>0$, with the equations holding for 
$t\in [0,\tau_{\infty})$.
\end{theorem}

\begin{remark}
Assuming uniqueness, the system determines a piecewise 
deterministic Markov process in the sense of 
\cite{Dav93}.
If one defines 
\[\beta_k(z)=\int_{{\Bbb L}_1}\lambda_k(z+y)\pi^z(dy),\quad z\in\Pi_
2{\Bbb E},\]
the description of the system will simplify.
\end{remark}

We still need to address conditions for the relative 
compactness of the sequence of occupation measures.  
If $(I-\Pi_2){\Bbb E}$ is compact, relative compactness is 
immediate.  Otherwise, it is natural to look for some 
kind of Lyapunov function.  Note that if 
$ $$\gamma_c^N=\inf\{t:|Z^{N,r_2}(t)|\geq c\}$, then
\[f(Z^{N,r_2}(t\wedge\gamma_c^N)-f(Z^{N,r_2}(0))-\int_0^{t\wedge\gamma_
c^N}N^{r_2}{\Bbb B}_Nf(Z^{N,r_2}(s))ds\]
is a martingale for all locally bounded $f$.

\begin{lemma}
Let $h_{q,\theta}$ and $\psi_{q,\theta}$  be as in (\ref{uint}).  Suppose  
$f^{_N}_{q,\theta}$ are nonnegative functions and there exist 
positive
constants $c_1$, $c_2$ such that 
\[\sup_NN^{r_2}{\Bbb B}_Nf^N_{q,\theta}(z)<c_1-c_2\psi_{q,\theta}
(h_{q,\theta}((I-\Pi_2)z))\]
for all $z$ satisfying $|\Pi_2z|\leq q$ and for each $c\in {\Bbb R}$,
\[\sup\{|(I-\Pi_2)z|:|\Pi_2z|\mbox{\rm \ and }\sup_N{}_{}N^{r_2}{\Bbb B}_
Nf^N_{q,\theta}(z)\geq c\}<\infty .\]
Then for each $t>0$, $\{V_1^{N,r_2}\}$  
is relatively compact and (\ref{uint}) is stochastically 
bounded.
\end{lemma}

\newpage
                                             
\setcounter{equation}{0}

\section{Examples}\label{sectexamp}
We give some additional examples that demonstrate how 
identifying exponents satisfying the balance condition 
leads to reasonable approximations to the original model.  
For a ``production level'' example, see the analysis of an 
E-coli heat shock model in \cite{Kan10}.

\subsection{Goutsias's model of regulated transcription}

We consider the following model of transcription 
regulation introduced in \cite{Gou05} 
and studied further in \citet*{MBS07}.  The model involves 
six species
\[\begin{array}{rcll}
X_1&=&\mbox{\rm \ \# of M}&\quad\mbox{\rm Protein monomer}\\
X_2&=&\mbox{\rm \ \# of D}&\quad\mbox{\rm Transcription factor}\\
X_3&=&\mbox{\rm \ \# of RNA}&\quad\mbox{\rm mRNA}\\
X_4&=&\mbox{\rm \ \# of DNA}&\quad\mbox{\rm Unbound DNA}\\
X_5&=&\mbox{\rm \ \# of DNA$\cdot$D}&\quad\mbox{\rm DNA bound at one site}\\
X_6&=&\mbox{\rm \ \# of DNA$\cdot$2D}&\quad\mbox{\rm DNA bound at two sites}\end{array}
\]
and ten reactions:
\[\begin{array}{rcl}
RNA&\rightarrow&RNA+M\\
M&\rightarrow&\emptyset\\
DNA\cdot D&\rightarrow&RNA+DNA\cdot D\\
RNA&\rightarrow&\emptyset\\
DNA+D&\rightarrow&DNA\cdot D\\
DNA\cdot D&\rightarrow&DNA+D\\
DNA\cdot D+D&\rightarrow&DNA\cdot 2D\\
DNA\cdot 2D&\rightarrow&DNA\cdot D+D\\
M+M&\rightarrow&D\\
D&\rightarrow&2M\;.\end{array}
\]
Taking the volume $V=1$,
the corresponding system of equations becomes
\begin{eqnarray*}
X_1(t)&=&X_1(0)+Y_1(\kappa'_1\int_0^tX_3(s)ds))+2Y_{10}(\kappa'_{
10}\int_0^tX_2(s)ds)\\
&&\qquad -Y_2(\kappa_2'\int_0^tX_1(s)ds)-2Y_9(\kappa'_9\int_0^tX_
1(s)(X_1(s)-1)ds\\
X_2(t)&=&X_2(0)+Y_6(\kappa'_6\int_0^tX_5(s)ds)+Y_8(\kappa'_8\int_
0^tX_6(s)ds)+Y_9(\kappa'_9\int_0^tX_1(s)(X_1(s)-1)ds\\
&&\qquad -Y_5(\kappa_5'\int_0^tX_2(s)X_4(s)ds)-Y_7(\kappa'_7\int_
0^tX_2(s)X_5(s)ds)-Y_{10}(\kappa'_{10}\int_0^tX_2(s)ds)\\
X_3(t)&=&X_3(0)+Y_3(\kappa'_3\int_0^tX_5(s)ds)-Y_4(\kappa'_4\int_
0^tX_3(s)ds)\\
X_4(t)&=&X_4(0)+Y_6(\kappa'_6\int_0^tX_5(s)ds)-Y_5(\kappa'_5\int_
0^tX_2(s)X_4(s)ds)\\
X_5(t)&=&X_5(0)+Y_5(\kappa'_5\int_0^tX_2(s)X_4(s)ds)+Y_8(\kappa'_
8\int_0^tX_6(s)ds)\\
&&\qquad -Y_6(\kappa'_6\int_0^tX_5(s)ds)-Y_7(\kappa'_7\int_0^tX_2
(s)X_5(s)ds)\\
X_6(t)&=&X_6(0)+Y_7(\kappa'_7\int_0^tX_2(s)X_5(s)ds)-Y_8(\kappa'_
8\int_0^tX_6(s)ds)\;.\end{eqnarray*}

\subsection{A scaling with two fast reactions} In his 
analysis of the model, Goutsias assumes two time-scales 
and identifies reactions 9 and 10 as ``fast'' reactions.  In 
our approach, that is the same as assuming 
$\beta_9=\beta_{10}>\beta_1=\cdots =\beta_8$, so we take $N_0=100$, $
\beta_9=\beta_{10}=0$ 
and $\beta_1=\cdots =\beta_8=-1$.  Recall the relationships $\kappa_
k'=\kappa_kN_0^{\beta_k}$ 
(we are assuming the volume $V=1$)
and $\rho_k=\beta_k+\nu_k\cdot\alpha$.  Employing the rate constants from 
\cite{Gou05}, and taking $\alpha_i=0$ for all $i$, we have 

\begin{center}
\begin{longtable}{ll|ll|lr}
\caption{Scaling exponents for reaction rates}\\

\hline
Rates & & Scaled Rates & & $\rho$ &\\
\hline
$\kappa'_1$ & $4.30\times 10^{-2}$  & $\kappa_1$ & $4.30$  & $\rho_
1$ & $-1$\\
$\kappa'_2$ & $7.00\times 10^{-4}$  & $\kappa_2$ & $0.07$ & $\rho_
2$ & $-1$\\
$\kappa'_3$ & $7.15\times 10^{-2}$ & $\kappa_3$ & $7.15$ & $\rho_
3$ & $-1$\\
$\kappa'_4$ & $3.90\times 10^{-3}$  & $\kappa_4$ & $0.390$  & $\rho_
4$ & $-1$\\
$\kappa'_5$ & $1.99\times 10^{-2}$  & $\kappa_5$ & $1.99$ & $\rho_
5$ & $-1$\\
$\kappa'_6$ & $4.79\times 10^{-1}$  & $\kappa_6$ & $47.9$  & $\rho_
6$ & $-1$\\
$\kappa'_7$ & $1.99\times 10^{-4}$ & $\kappa_7$ & $0.0199$  & $\rho_
7$ & $-1$\\
$\kappa'_8$ & $8.77\times 10^{-12}$ & $\kappa_8$ & $8.77\times 10^{
-10}$  & $\rho_8$ & $-1$\\
$\kappa'_9$ & $8.30\times 10^{-2}$& $\kappa_9$ & $0.0830$ &  $\rho_9$ & $0$\\
$\kappa'_{10}$ & $5.00\times 10^{-1}$& $\kappa_{10}$ & $0.500$  & $
\rho_{10}$ & $0$\\
\hline
\end{longtable}
\end{center}

Then, for $\gamma =0$, $(Z_1^{N,0},Z_2^{N,0})$ converges to the solution of
\begin{eqnarray*}
Z_1^0(t)&=&X_1(0)+2Y_{10}(\kappa_{10}\int_0^tZ^0_2(s)ds)-2Y_9(\kappa_
9\int_0^tZ^0_1(s)(Z^0_1(s)-1)ds\\
Z^0_2(t)&=&X_2(0)+Y_9(\kappa_9\int_0^tZ^0_1(s)(Z^0_1(s)-1)ds-Y_{1
0}(\kappa_{10}\int_0^tZ^0_2(s)ds),\end{eqnarray*}
and for $k>2$, $Z_k^{N,0}$ converges to $X_k(0)$.  

For $\gamma =1$, the kind of argument employed in (\ref{intave}) 
implies 
\begin{equation}\kappa_9\int_0^tZ_1^{N,1}(s)(Z_1^{N,1}(s)-1)ds-\int_
0^t\kappa_{10}Z_2^{N,1}(s)ds\rightarrow 0\label{intave2}\end{equation}
but does not lead to a closed system for the limit of 
$(Z_3^{N,1},\ldots ,Z_6^{N,1})$.  To obtain a closed limiting system,
we introduce the following auxiliary variable:
\begin{eqnarray*}
Z^{N,1}_{12}(t)&=&Z^{N,1}_1(t)+2Z^{N,1}_2(t)\\
&=&Z^{N,1}_{12}(0)+Y_1(\kappa^{}_1\int_0^tZ^{N,1}_3(s)ds))+2Y_6(\kappa^{}_
6\int_0^tZ^{N,1}_5(s)ds)+2Y_8(\kappa^{}_8\int_0^tZ^{N,1}_6(s)ds)\\
&&\qquad -2Y_5(\kappa_5^{}\int_0^tZ^{N,1}_2(s)Z^{N,1}_4(s)ds)-2Y_
7(\kappa^{}_7\int_0^tZ^{N,1}_2(s)Z^{N,1}_5(s)ds)-Y_2(\kappa_2^{}\int_
0^tZ^{N,1}_1(s)ds),\end{eqnarray*}
and observe that the conditional equilibrium distribution 
satisfies 
\begin{eqnarray*}
\kappa_9(z_1+2)(z_1+1)\mu_s(z_1+2,z_2-1)+\kappa_{10}(z_2+1)\mu_s(
z_1-2,z_2+1)\\
=(\kappa_9z_1(z_1-1)+\kappa_{10}z_2)\mu_s(z_1,z_2)\end{eqnarray*}
and is uniquely determined by the requirement that 
\[z_1+2z_2=Z_{12}^1(s),\]
where $Z_{12}^1$ is the limit of $Z_{12}^{N,1}$.  For $m=z_1+2z_2$, the 
conditional equilibrium distribution is
\begin{equation}\mu_m(z_1,z_2)=M_m\frac {(\kappa_{10}/\kappa_9)^{
z_1+z_2}}{z_1!z_2!},\label{lavd}\end{equation}
where $M_m$ is a normalizing constant making $\mu_m$ a probability 
distribution on the collection of $(z_1,z_2)$ such that $z_1$ and 
$z_2$ are nonnegative integers satisfying $z_1+2z_2=m$.  
Define
\begin{equation}\alpha (m)=\int z_2\mu_m(dz_1,dz_2)=M_m\sum_{1\leq 
z_2\leq m/2}\frac {(\kappa_{10}/\kappa_9)^{(m-z_2)}}{(m-2z_2)!(z_
2-1)!},\label{alph}\end{equation}
and observe that $m-2\alpha (m)=\int z_1\mu_m(dz_1,dz_2)$.
Then $(Z_{12}^{N,1},Z_3^{N,1},\ldots ,Z_6^{N,1})$ converges to the solution of
\begin{eqnarray*}
Z^1_{12}(t)&=&Z^1_{12}(0)+Y_1(\kappa^{}_1\int_0^tZ^1_3(s)ds))+2Y_
6(\kappa^{}_6\int_0^tZ^1_5(s)ds)+2Y_8(\kappa^{}_8\int_0^tZ^1_6(s)
ds)\\
&&\qquad -2Y_5(\kappa_5^{}\int_0^t\alpha (Z^1_{12}(s))Z^1_4(s)ds)
-2Y_7(\kappa^{}_7\int_0^t\alpha (Z^1_{12}(s))Z^1_5(s)ds)\\
&&\qquad -Y_2(\kappa_2^{}\int_0^t(Z^1_{12}(s)-2\alpha (Z^1_{12}(s
)))ds)\\
Z^1_3(t)&=&Z^1_3(0)+Y_3(\kappa_3\int_0^tZ^1_5(s)ds)-Y_4(\kappa_4\int_
0^tZ^1_3(s)ds)\\
Z^1_4(t)&=&Z^1_4(0)+Y_6(\kappa_6\int_0^tZ^1_5(s)ds)-Y_5(\kappa_5\int_
0^t\alpha (Z^1_{12}(s))Z^1_4(s)ds)\\
Z^1_5(t)&=&Z^1_5(0)+Y_5(\kappa_5\int_0^t\alpha (Z^1_{12}(s))Z^1_4
(s)ds)+Y_8(\kappa_8\int_0^tZ^1_6(s)ds)\\
&&\qquad -Y_6(\kappa_6\int_0^tZ^1_5(s)ds)-Y_7(\kappa_7\int_0^t\alpha 
(Z^1_{12}(s))Z^1_5(s)ds)\\
Z^1_6(t)&=&Z^1_6(0)+Y_7(\kappa_7\int_0^t\alpha (Z^1_{12}(s))Z^1_5
(s)ds)-Y_8(\kappa_8\int_0^tZ^1_6(s)ds)\;,\end{eqnarray*}
which is essentially the approximation obtained by 
Goutsias.  Note that the ``fast'' reactions, reactions 9 and 
10, have been eliminated from the model.

This system is not entirely satisfactory as $\alpha (m)$ is not 
computable analytically.  For simulations, values of $\alpha (m)$ 
could be precomputed using (\ref{alph}).  \citet*{ELV07} 
suggest a Monte Carlo approach for computing $\alpha (m)$ as 
needed.  Goutsias suggests a way of approximating the 
transition rates which is equivalent to the following:  
The limit in (\ref{intave2}) implies 
\begin{equation}\kappa_{10}\alpha (m)=\kappa_9\int z_1(z_1-1)\mu_
m(dz_1,dz_2)\label{momgou}\end{equation}
as can be verified directly from the definition of $\mu_m$.  
A moment 
closure argument suggests replacing (\ref{momgou}) by
\begin{eqnarray*}
\kappa_{10}\alpha (m)&=&\kappa_9\int z_1\mu_m(dz_1,dz_2)\int (z_1
-1)\mu_m(dz_1,dz_2)\\
&=&\kappa_9(m-2\alpha (m))(m-2\alpha (m)-1),\end{eqnarray*}
which gives a quadratic equation for the approximation 
for $\alpha (m)$.  

\subsection{Alternative scaling}
 Observe that $\kappa_9'<\kappa_6'$, so reaction 6 is actually ``faster'' 
than reaction 9.  
Consequently, it is reasonable to 
look for a different solution of the balance conditions 
with $\beta_{10}=\beta_6>\beta_9$.  Drop the assumption that $\alpha_
i=0$, 
and consider a subset of the balance equations.  Recall 
that $\rho_k=\beta_k+\nu_k\cdot\alpha$.

\begin{center}
\begin{longtable}{ll}
\caption{Balance equations}\\
\hline
Variable & Balance equation \\
\hline
$X_1$ & $\rho_1\vee\rho_{10}=\rho_2\vee\rho_9$ \\
$X_2$ & $\rho_6\vee\rho_8\vee\rho_9=\rho_5\vee\rho_7\vee\rho_{10}$ \\
$X_3$ & $\rho_3=\rho_4$ \\
$X_4$ & $\rho_5=\rho_6$ \\
$X_5$ & $\rho_5\vee\rho_8=\rho_6\vee\rho_7$ \\
$X_6$ & $\rho_7=\rho_8$ \\
$X_1+2X_2+2X_5+4X_6$ & $\rho_1=\rho_2$ \\
$X_2+X_5+2X_6$ & $\rho_9=\rho_{10}$ \\
$X_5+X_6$ & $\rho_5=\rho_6$ \\
$X_4+X_5+X_6$ & $0=0$ \\
$X_4+X_5$ & $\rho_8=\rho_7$ \\
\hline
\end{longtable}
\end{center}

We take $N_0=100$, $\alpha_1=\alpha_2=1$,  and $\alpha_i=0$ for $
3\leq i\leq 6$.
We see that the following exponents satisfy the 
balance conditions and the additional requirement that 
$\kappa_k'\geq\kappa_l'$ imply $\beta_k\geq\beta_l$, except for $
\beta_8$, the exponent 
associated with the extremely small rate constant $\kappa'_8$.  
Recall that $\kappa_k$ is determined by 
the requirement $\kappa_k'=\kappa_kN_0^{\beta_k}$.

\begin{center}
\begin{longtable}{ll|lr|lr|lr}
\caption{Scaling exponents for reaction rates}\\

\hline
Rates & & Exponents & & Scaled Rates & & $\rho$ &\\
\hline
$\kappa'_1$ & $4.30\times 10^{-2}$ & $\beta_1$ & $-1$  & $\kappa_
1$ & $4.30$  & $\rho_1$ & $-1$\\
$\kappa'_2$ & $7.00\times 10^{-4}$ & $\beta_2$ & $-2$ & $\kappa_2$ & $
7.00$ & $\rho_2$ & $-1$\\
$\kappa'_3$ & $7.15\times 10^{-2}$ & $\beta_3$ & $-1$ & $\kappa_3$ & $
7.15$ & $\rho_3$ & $-1$\\
$\kappa'_4$ & $3.90\times 10^{-3}$ & $\beta_4$ & $-1$ & $\kappa_4$ & $
0.390$  & $\rho_4$ & $-1$\\
$\kappa'_5$ & $1.99\times 10^{-2}$ & $\beta_5$ & $-1$ & $\kappa_5$ & $
1.99$ & $\rho_5$ & $0$\\
$\kappa'_6$ & $4.79\times 10^{-1}$ & $\beta_6$ & $0$  & $\kappa_
6$ & $.479$  & $\rho_6$ & $0$\\
$\kappa'_7$ & $1.99\times 10^{-4}$ & $\beta_7$ & $-3$ & $\kappa_7$ & $
199$  & $\rho_7$ & $-2$\\
$\kappa'_8$ & $8.77\times 10^{-12}$ & $\beta_8$ & $-2$ & $\kappa_
8$ & $8.77\times 10^{-8}$  & $\rho_8$ & $-2$\\
$\kappa'_9$ & $8.30\times 10^{-2}$ & $\beta_9$ & $-1$ & $\kappa_9$ & $
8.30$ &  $\rho_9$ & $1$\\
$\kappa'_{10}$ & $5.00\times 10^{-1}$ & $\beta_{10}$ & $0$ & $\kappa_{
10}$ & $0.500$  & $\rho_{10}$ & $1$\\
\hline
\end{longtable}
\end{center}

Defining $Z^{N,\gamma}_i(t)=N^{-\alpha_i}X_i^N(N^{\gamma}t)$ and $
\kappa_k=N^{-\beta_k}\kappa_k'$, 
\begin{eqnarray*}
Z_1^{N,\gamma}(t)&=&Z_1^{N,\gamma}(0)+N^{-1}Y_1(\int_0^t\kappa_1N^{
\gamma -1}Z_3^{N,\gamma}(s)\,ds)+2N^{-1}Y_{10}(\int_0^t\kappa_{10}
N^{\gamma +1}Z_2^{N,\gamma}(s)\,ds)\\
&&-N^{-1}Y_2(\int_0^t\kappa_2N^{\gamma -1}Z_1^{N,\gamma}(s)\,ds)-
2N^{-1}Y_9(\int_0^t\kappa_9N^{\gamma +1}Z_1^{N,\gamma}(s)(Z_1^{N,
\gamma}(s)-N^{-1})\,ds)\\
Z_2^{N,\gamma}(t)&=&Z_2^{N,\gamma}(0)+N^{-1}Y_6(\int_0^t\kappa_6N^{
\gamma}Z_5^{N,\gamma}(s)\,ds)+N^{-1}Y_8(\int_0^t\kappa_8N^{\gamma 
-2}Z_6^{N,\gamma}(s)\,ds)\\
&&+N^{-1}Y_9(\int_0^t\kappa_9N^{\gamma +1}Z_1^{N,\gamma}(s)(Z_1^{
N,\gamma}(s)-N^{-1})\,ds)-N^{-1}Y_5(\int_0^t\kappa_5N^{\gamma}Z_2^{
N,\gamma}(s)Z_4^{N,\gamma}(s)\,ds)\\
&&-N^{-1}Y_7(\int_0^t\kappa_7N^{\gamma -2}Z_2^{N,\gamma}(s)Z_5^{N
,\gamma}(s)\,ds)-N^{-1}Y_{10}(\int_0^t\kappa_{10}N^{\gamma +1}Z_2^{
N,\gamma}(s)\,ds)\\
Z_3^{N,\gamma}(t)&=&Z_3^{N,\gamma}(0)+Y_3(\int_0^t\kappa_3N^{\gamma 
-1}Z_5^{N,\gamma}(s)\,ds)-Y_4(\int_0^t\kappa_4N^{\gamma -1}Z_3^{N
,\gamma}(s)\,ds)\\
Z_4^{N,\gamma}(t)&=&Z_4^{N,\gamma}(0)+Y_6(\int_0^t\kappa_6N^{\gamma}
Z_5^{N,\gamma}(s)\,ds)-Y_5(\int_0^t\kappa_5N^{\gamma}Z_2^{N,\gamma}
(s)Z_4^{N,\gamma}(s)\,ds)\\
Z_5^{N,\gamma}(t)&=&Z_5^{N,\gamma}(0)+Y_5(\int_0^t\kappa_5N^{\gamma}
Z_2^{N,\gamma}(s)Z_4^{N,\gamma}(s)\,ds)+Y_8(\int_0^t\kappa_8N^{\gamma 
-2}Z_6^{N,\gamma}(s)\,ds)\\
&&-Y_6(\int_0^t\kappa_6N^{\gamma}Z_5^{N,\gamma}(s)\,ds)-Y_7(\int_
0^t\kappa_7N^{\gamma -2}Z_2^{N,\gamma}(s)Z_5^{N,\gamma}(s)\,ds)\\
Z_6^{N,\gamma}(t)&=&Z_6^{N,\gamma}(0)+Y_7(\int_0^t\kappa_7N^{\gamma 
-2}Z_2^{N,\gamma}(s)Z_5^{N,\gamma}(s)\,ds)-Y_8(\int_0^t\kappa_8N^{
\gamma -2}Z_6^{N,\gamma}(s)\,ds).\end{eqnarray*}

Useful auxiliary variables include
\begin{eqnarray*}
&&NZ_1^{N,\gamma}(t)+2NZ_2^{N,\gamma}(t)+2Z_5^{N,\gamma}(t)+4Z_6^{
N,\gamma}(t)=NZ_1^N(0)+2NZ_2^N(0)+2Z_5^N(0)+4Z_6^N(0)\\
&&\qquad +Y_1(\int_0^t\kappa_1N^{\gamma -1}Z_3^{N,\gamma}(s)\,ds)
-Y_2(\int_0^t\kappa_2N^{\gamma -1}Z_1^{N,\gamma}(s)\,ds)\\
&&NZ_2^{N,\gamma}(t)+Z_5^{N,\gamma}(t)+2Z_6^{N,\gamma}(t)=NZ_2^N(
0)+Z_5^N(0)+2Z_6^N(0)\\
&&\qquad +Y_9(\int_0^t\kappa_9N^{\gamma +1}Z_1^{N,\gamma}(s)(Z_1^{
N,\gamma}(s)-N^{-1})\,ds)-Y_{10}(\int_0^t\kappa_{10}N^{\gamma +1}
Z_2^{N,\gamma}(s)\,ds)\\
&&Z_5^{N,\gamma}(t)+Z_6^{N,\gamma}(t)=Z_5^N(0)+Z_6^N(0)+Y_5(\int_
0^t\kappa_5N^{\gamma}Z_2^{N,\gamma}(s)Z_4^{N,\gamma}(s)\,ds)-Y_6(
\int_0^t\kappa_6N^{\gamma}Z_5^{N,\gamma}(s)\,ds)\\
&&Z_4^{N,\gamma}(t)+Z_5^{N,\gamma}(t)+Z_6^{N,\gamma}(t)=Z_4^N(0)+
Z_5^N(0)+Z_6^N(0)\\
&&Z_4^{N,\gamma}(t)+Z_5^{N,\gamma}(t)=Z_4^N(0)+Z_5^N(0)+Y_8(\int_
0^t\kappa_8N^{\gamma -2}Z_6^{N,\gamma}(s)\,ds)-Y_7(\int_0^t\kappa_
7N^{\gamma -2}Z_2^N(s)Z_5^{N,\gamma}(s)\,ds)\end{eqnarray*}

For $\gamma =0$, the limiting system is the piecewise 
deterministic model
\begin{eqnarray}
Z^0_1(t)&=&Z_1(0)+\int_0^t\big(2\kappa_{10}Z^0_2(s)-2\kappa_9Z^0_
1(s)^2\big)\,ds\nonumber\\
Z^0_2(t)&=&Z_2(0)+\int_0^t\big(\kappa_9Z^0_1(s)^2-\kappa_{10}Z^0_
2(s)\big)\,ds\label{altgam0}\\
Z^0_4(t)&=&Z_4(0)+Y_6(\int_0^t\kappa_6Z^0_5(s)\,ds)-Y_5(\int_0^t\kappa_
5Z^0_2(s)Z^0_4(s)\,ds)\nonumber\\
Z^0_5(t)&=&Z_5(0)+Y_5(\int_0^t\kappa_5Z^0_2(s)Z^0_4(s)\,ds)-Y_6(\int_
0^t\kappa_6Z^0_5(s)\,ds),\nonumber\end{eqnarray}
with $Z_3^0(t)\equiv Z_3(0)$ and $Z_6^0(t)\equiv Z_6(0)$.

For $\gamma =1$, we introduce the auxiliary variables

\begin{eqnarray*}
Z_{12}^{N,1}(t)&\equiv&Z_1^{N,1}(t)+2Z_2^{N,1}(t)\\
Z_{45}^{N,1}(t)&\equiv&Z_4^{N,1}(t)+Z_5^{N,1}(t)\\
&=&Z_4^N(0)+Z_5^N(0)+Y_8(\int_0^t\kappa_8N^{-1}Z_6^{N,1}(s)\,ds)-
Y_7(\int_0^t\kappa_7N^{-1}Z_2^{N,1}(s)Z_5^{N,1}(s)\,ds).\end{eqnarray*}
Observing that $Z_{12}^{N,1}$ is asymptotically the same as 
$Z_1^{N,1}+2Z_2^{N,1}+2N^{-1}Z_5^{N,1}+4N^{-1}Z_6^{N,1}$,
 $Z_{12}^{N,1}$ 
converges to $Z_{12}^1(t)\equiv Z_{12}(0)=\lim_{N\rightarrow\infty}
(Z_1^N(0)+2Z_2^N(0))$.  In 
particular, $Z_{12}^1$ is constant in time.  We also have 
$Z_{45}^1(t)\equiv Z_{45}(0)=\lim_{N\rightarrow\infty}(Z_4^N(0)+Z_
5^N(0))$.

Let $V^{N,1}$ denote the occupation measure for 
$(Z_1^{N,1},Z_2^{N,1},Z_4^{N,1},Z_5^{N,1})$.  The stochastic boundedness of 
$Z_{12}^{N,1}$ and $Z_{45}^{N,1}$ ensures the relative compactness of $
\{V^{N,1}\}$, 
and as in Section \ref{sectave}, $V^{N,1}$ converges to 
$V^1(dz,ds)=v_s(dz)ds$  where $v_s$ satisfies
\[\int {\Bbb C}fv_s(dz)=0\]
and 
\begin{eqnarray*}
{\Bbb C}f(z_1,z_2,z_4,z_5)&=&(\kappa_{10}z_2-\kappa_9z_1^2)(2\partial_{
z_1}f(z)-\partial_{z_2}f(z))\\
&&\qquad +\kappa_6z_5(f(z+e_4-e_5)-f(z))\\
&&\qquad +\kappa_5z_2z_4(f(z-e_4+e_5)-f(z)).\end{eqnarray*}
Consequently, $v_s$ is uniquely determined for each $s$ by the 
requirement that $z_1+2z_2=Z^1_{12}(s)=Z_{12}(0)$  and 
$z_4+z_5=Z_{45}^1(s)=Z_{45}(0)$, and 
hence,
\[v_s(dz)=\delta_{\varphi_1(Z_{12}(0))}(dz_1)\delta_{\varphi_2(Z_{
12}(0))}(dz_2)B(Z_{45}(0),\frac {\kappa_6}{\kappa_6+\kappa_5\varphi_
2(Z_{12}(0))},dz_4,dz_5),\]
where 
\begin{eqnarray*}
\varphi_1(y)&=&\frac {\sqrt {\kappa_{10}^2+8\kappa_9\kappa_{10}y}
-\kappa_{10}}{4\kappa_9}\\
\varphi_2(y)&=&\frac {4\kappa_9y+\kappa_{10}-\sqrt {\kappa_{10}^2
+8\kappa_9\kappa_{10}y}}{8\kappa_9}\end{eqnarray*}
and 
$B(n,p,dz_4,dz_5)$ is given by the binomial distribution
\[P\{Z_4=k,Z_5=n-k\}={n\choose k}p^k(1-p)^{n-k}.\]
Averaging gives
\begin{equation}Z^1_3(t)=Z_3(0)+Y_3(\int_0^t\frac {\kappa_3\kappa_
5\varphi_2(Z_{12}(0))}{\kappa_6+\kappa_5\varphi_2(Z_{12}(0))}Z_{4
5}(0)\,ds)-Y_4(\int_0^t\kappa_4Z^1_3(s)\,ds).\label{altgam1}\end{equation}

Finally, for $\gamma =2$, dividing the equation for $Z_3^{N,2}$ by $
N$, 
we see that  
\[\int_0^tZ_3^{N,2}(s)ds\approx\frac {\kappa_3}{\kappa_4}\int_0^t
Z_5^{N,2}(s)ds,\]
and $(Z_{12}^{N,2},Z^{N,2}_{45},Z_6^{N,2})$ converges to the 
solution of
\begin{eqnarray*}
Z_{12}^2(t)&=&Z_{12}(0)+\int_0^t\left(\frac {\kappa_1\kappa_3}{\kappa_
4}\bar Z_5^2(s)-\kappa_2\varphi_1(Z_{12}^2(s))\right)ds\label{altgam2}\\
Z_{45}^2(t)&=&Z_{45}(0)+Y_8(\int_0^t\kappa_8Z_6^2(s)\,ds)-Y_7(\int_
0^t\kappa_7\varphi_2(Z_{12}^2(s))\bar {Z}_5^2(s)\,ds)\nonumber\\
Z^2_6(t)&=&Z_6(0)+Y_7(\int_0^t\kappa_7\varphi_2(Z_{12}^2(s))\bar {
Z}_5^2(s)\,ds)-Y_8(\int_0^t\kappa_8Z_6^2(s)\,ds)\nonumber\\
\bar {Z}^2_5(t)&=&\frac {\kappa_5\varphi_2(Z_{12}^2(t))}{\kappa_6
+\kappa_5\varphi_2(Z_{12}^2(t))}Z_{45}^2(t)\nonumber\end{eqnarray*}

\subsubsection{Simulation results} We compare simulation 
results for the full model with the approximations given 
by the limiting systems.  The mean and standard 
deviations of the number of molecules for each species 
or for the auxiliary variables of interest are given from 
$100$ simulations of the full model and from $1000$ 
simulations of the limiting systems.  The evolution of 
the processes in the full model is approximated by the 
evolution of the processes in the limiting system using 
the relationship 
\begin{eqnarray*}
X_i(t)\equiv X_i^{N_0}(t)&\approx&N_0^{\alpha_i}Z_i^{\gamma}(tN_0^{
-\gamma}).\end{eqnarray*}
Following \cite{Gou05}, initial values are taken as 
$X_1(0)=2$, $X_2(0)=6$, $X_5(0)=2$, and all other values equal 
to zero.  

\begin{figure}[htp]
\centering
\suppressfloats 
\hspace{-1cm}\includegraphics[totalheight=0.18\textheight,clip]{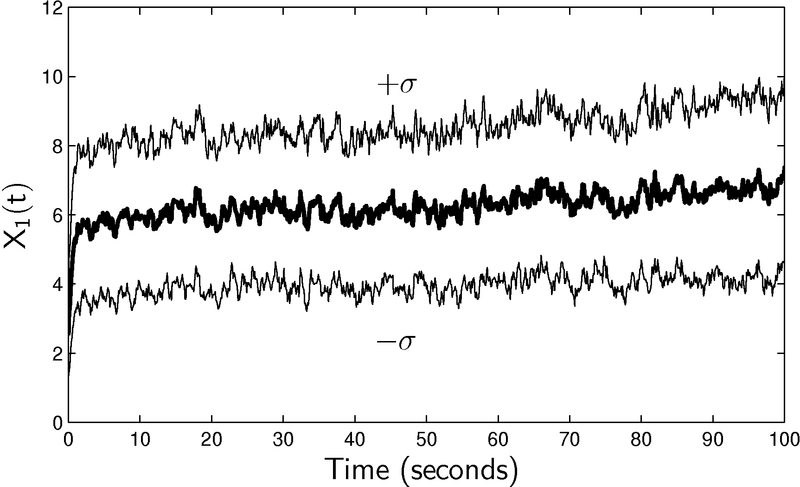}
\includegraphics[totalheight=0.18\textheight,clip]{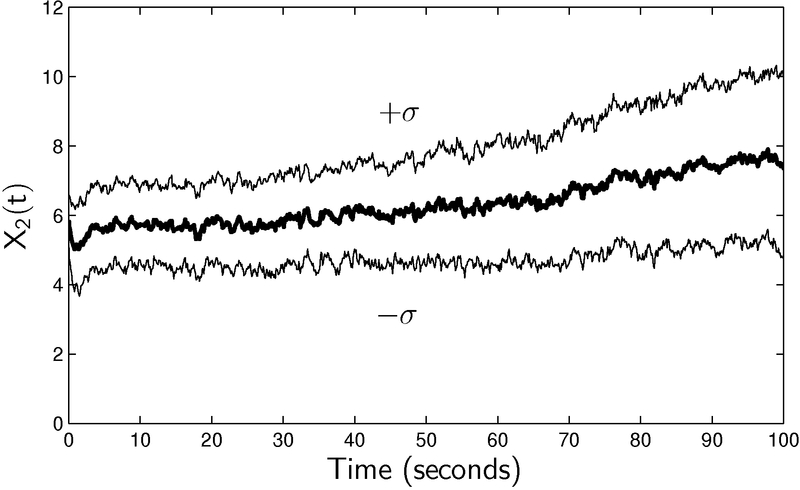}\\
\hspace{-1cm}\includegraphics[totalheight=0.18\textheight,clip]{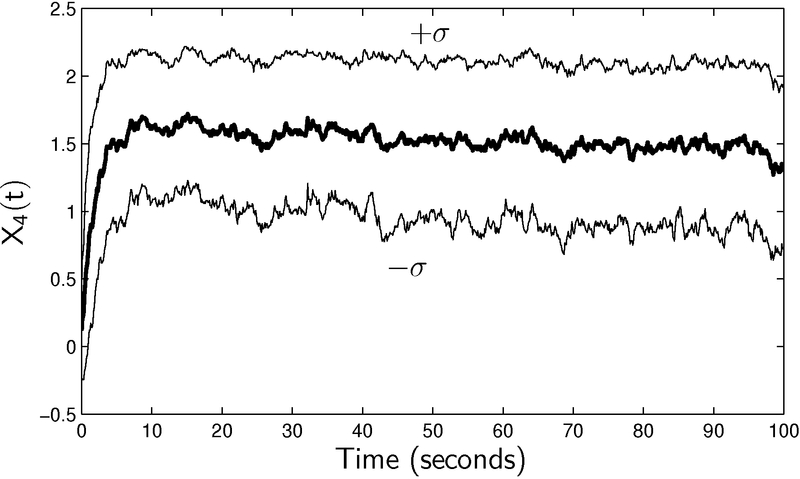}
\includegraphics[totalheight=0.18\textheight,clip]{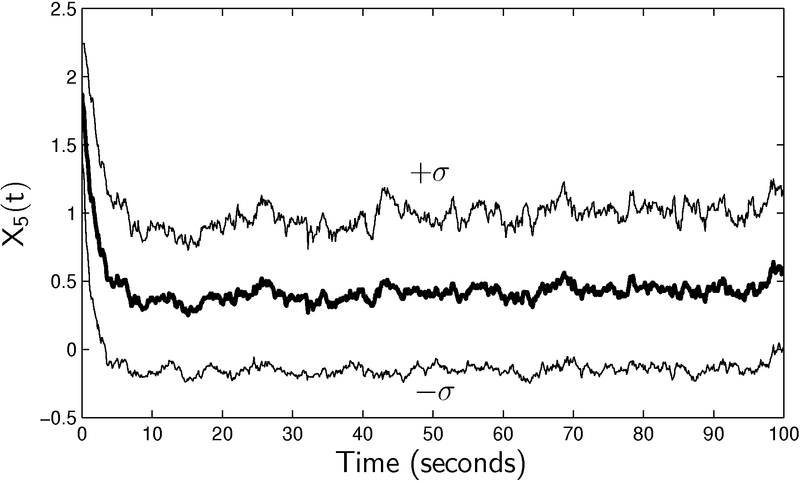}\\
\caption{Simulation of the full model during $t=0$ to $t=100$}
\label{fig: goutsias0}
\hspace{-1cm}\includegraphics[totalheight=0.18\textheight,clip]{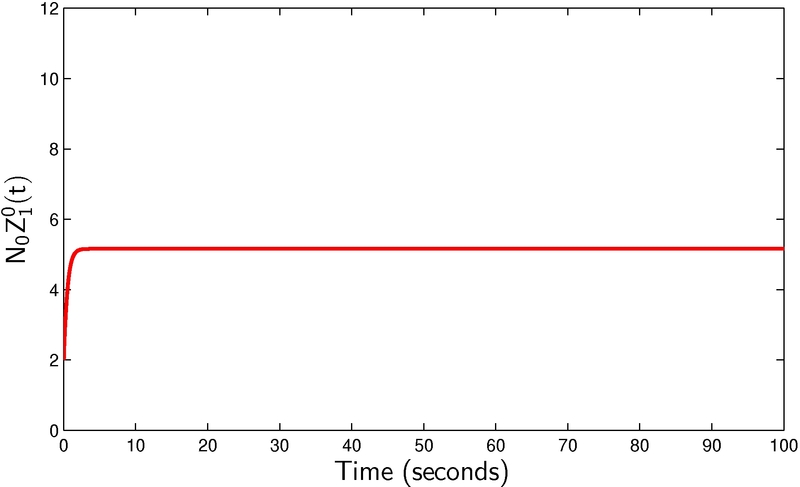}
\includegraphics[totalheight=0.18\textheight,clip]{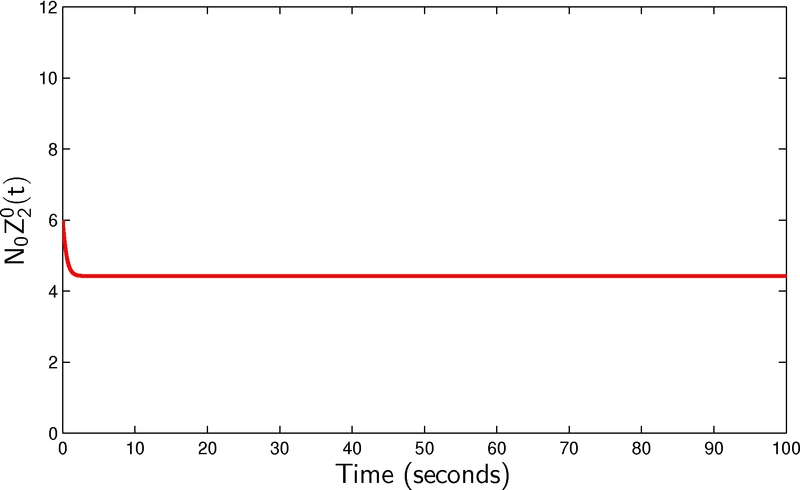}\\
\hspace{-1cm}\includegraphics[totalheight=0.18\textheight,clip]{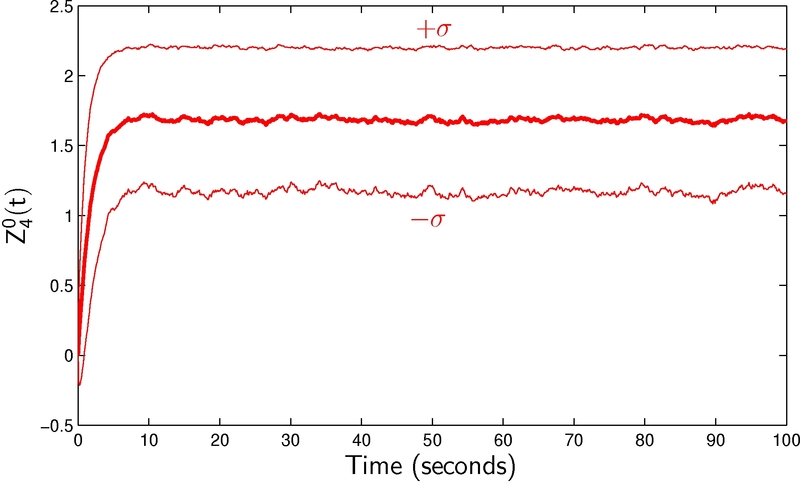}
\includegraphics[totalheight=0.18\textheight,clip]{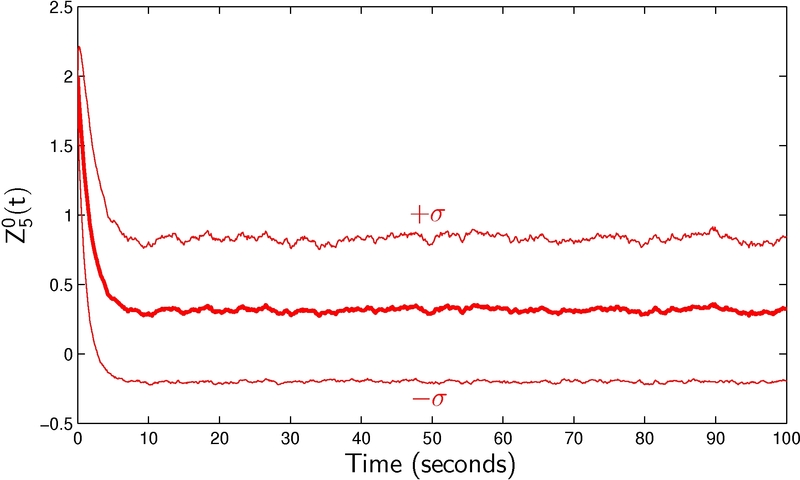}
\caption{Approximation using the limiting model for $\gamma=0$ in the 
alternative scaling}
\label{fig: alternativescaling0}
\end{figure}

For $\gamma =0$, we observe the evolution of the processes during 
the time interval $[0,100]$. 
The full model is reduced to the 
$4$-dimensional hybrid model (\ref{altgam0}) 
in which
 $Z_1^0$ and $Z_2^0$ 
are the solution of a pair of ordinary differential equations and 
$Z_4^0$ and $Z_5^0$ are discrete with transition intensities 
depending on $Z_2^0$.
The evolution of 
$X_1$, $X_2$, $X_4$, and $X_5$ in the full model is given in Figure 
\ref{fig:  goutsias0} and the evolution of the 
approximation is given in Figure \ref{fig:  
alternativescaling0}.  The exact simulations of the full 
model are done using  Gillespie's stochastic simulation 
algorithm (SSA) from \cite{Gill77}.  For the 
approximation, $Z_1^0$ and $Z_2^0$ are solved by the Matlab ODE 
solver, and $Z_4^0$ and $Z_5^0$ are computed by Gillespie's 
SSA taking $Z_2^0$ from the solution of ODE. The evolution of $X_1$ and $X_2$ are 
well captured by $Z_1^0$ and $Z_2^0$ in Figure \ref{fig:  
alternativescaling0}.  These deterministic values 
approximate the evolution of the mean of $X_1$ and $X_2$ 
given in Figure \ref{fig:  goutsias0} except for a slight 
increase over time in the simulation of the full model.  
Note that in the approximate model $Z_1^0(t)+2Z_2^0(t)$ is 
constant, but that is not the case in the full model.

\begin{figure}[htp]
\centering
\suppressfloats
\hspace{-1cm}\includegraphics[totalheight=0.18\textheight,clip]{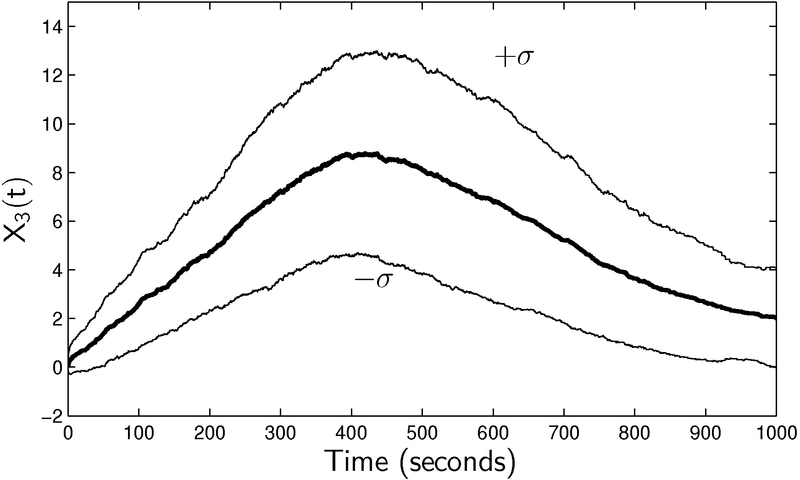}
\includegraphics[totalheight=0.18\textheight,clip]{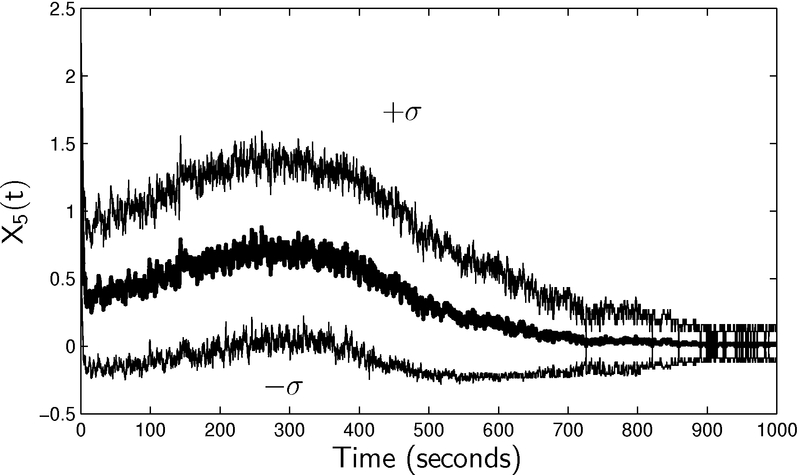}\\
\caption{Simulation of the full model during $t=0$ to $t=1000$}
\label{fig: goutsias1}
\hspace{-1cm}\includegraphics[totalheight=0.18\textheight,clip]{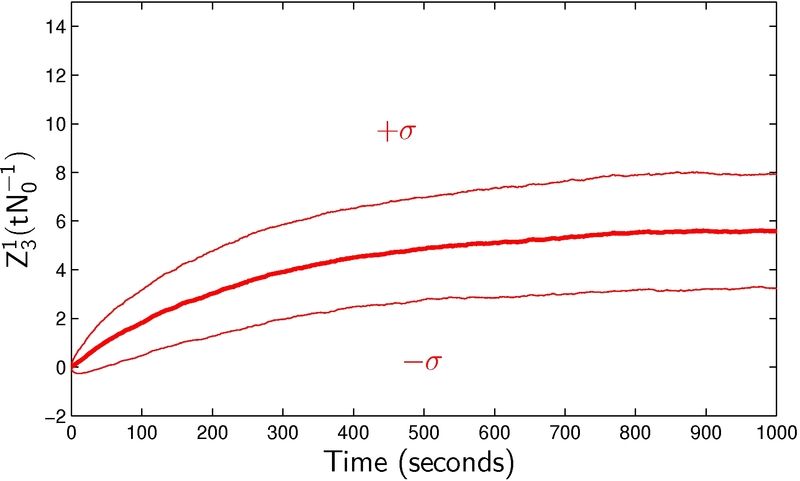}
\includegraphics[totalheight=0.18\textheight,clip]{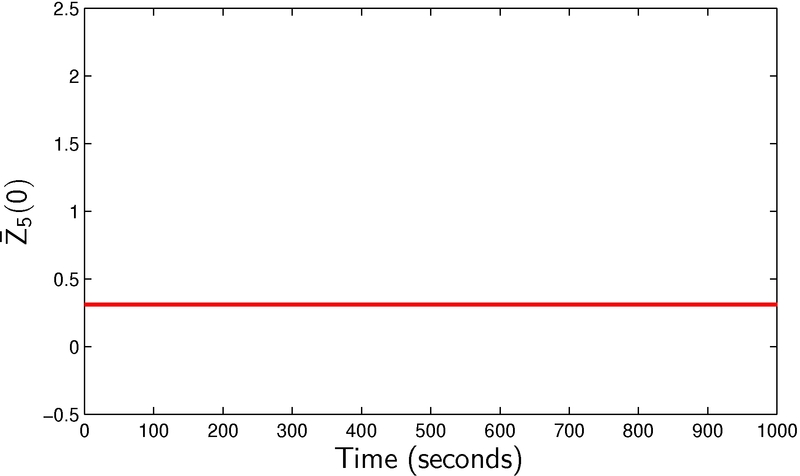}
\caption{Approximation using the limiting model for $\gamma =1$ in the 
alternative scaling}
\label{fig: alternativescaling1}
\end{figure}

For $\gamma =1$, we consider the evolution of the processes 
on the time interval $[0,1000]$.
The full model is reduced to the 
$1$-dimensional limiting system (\ref{altgam1})
with a single jump process $Z_3^1$.  
Comparing the governing equations for $Z_3^{N,1}$ and $Z_3^1$, the 
different behavior of the evolution of the two processes 
comes from the difference between $Z_5^{N,1}$ and  
$\bar{Z}_5^{1}(t)=\frac{\kappa_5\varphi_2(Z_{12}(0))}
{\kappa_6+\kappa_5\varphi_2(Z_{12}(0))}Z_{45}(0)$.  Therefore, plots of the 
evolution of both $X_3$ and $X_5$ in the exact simulation are 
given in Figure \ref{fig:  goutsias1}.  In Figure \ref{fig:  
alternativescaling1}, the evolution of $Z_3^1$ and  
of $\bar{Z}_5^1$ is given.  For both exact and approximate 
simulations, we use Gillespie's SSA.  In Figure 
\ref{fig:  goutsias1}, $Z_5^{N,1}$ increases slightly and then 
decreases to zero.  Since $\bar{Z}_5^1$ is 
approximated as constant in Figure \ref{fig:  
alternativescaling1}, the increase during the early time 
and the decrease to zero of $X_3$ is not captured by the 
approximation.

\begin{figure}[htp]
\centering
\suppressfloats
\hspace{-1cm}\includegraphics[totalheight=0.18\textheight,clip]{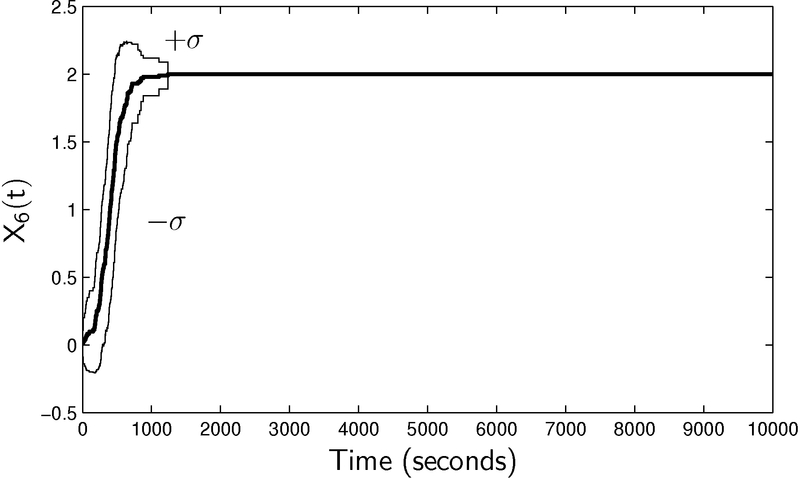}
\includegraphics[totalheight=0.18\textheight,clip]{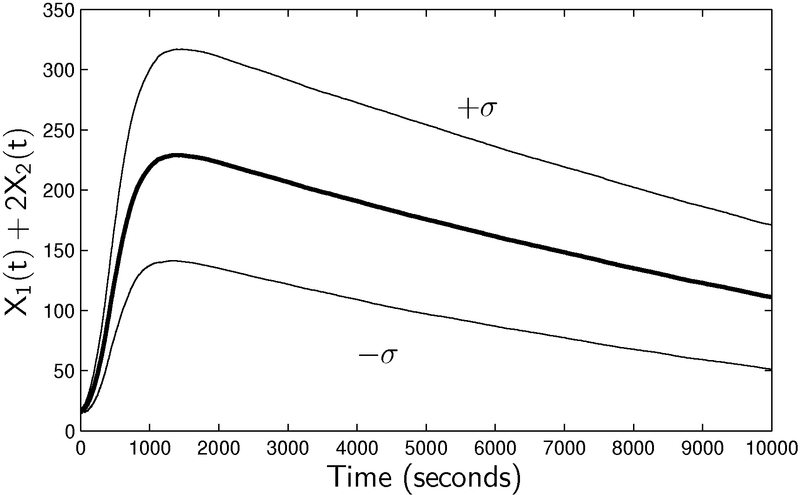}\\
\hspace{-1cm}\includegraphics[totalheight=0.18\textheight,clip]{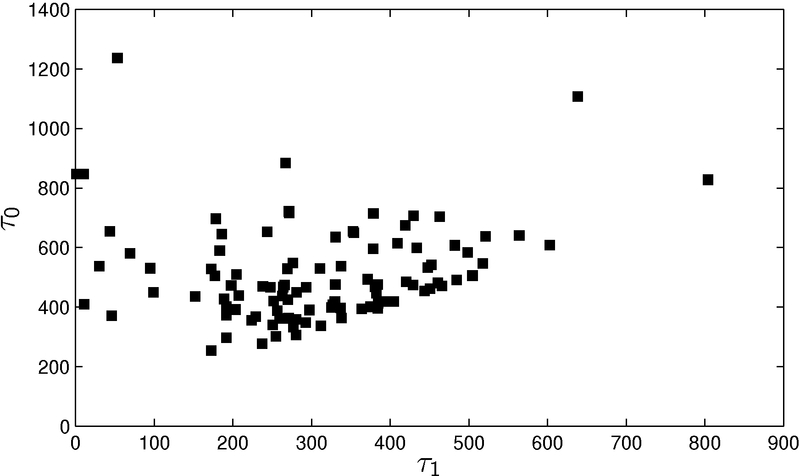}\\
\caption{Simulation of the full model during $t=0$ to $t=10000$}
\label{fig: goutsias2}
\hspace{-1cm}\includegraphics[totalheight=0.18\textheight,clip]{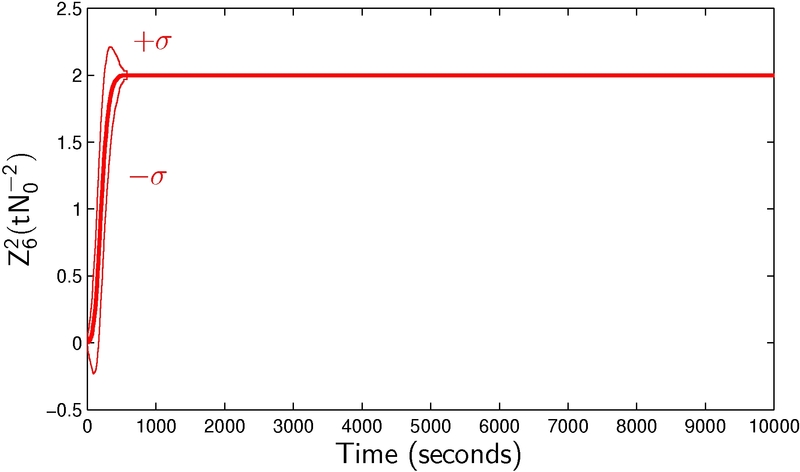}
\includegraphics[totalheight=0.18\textheight,clip]{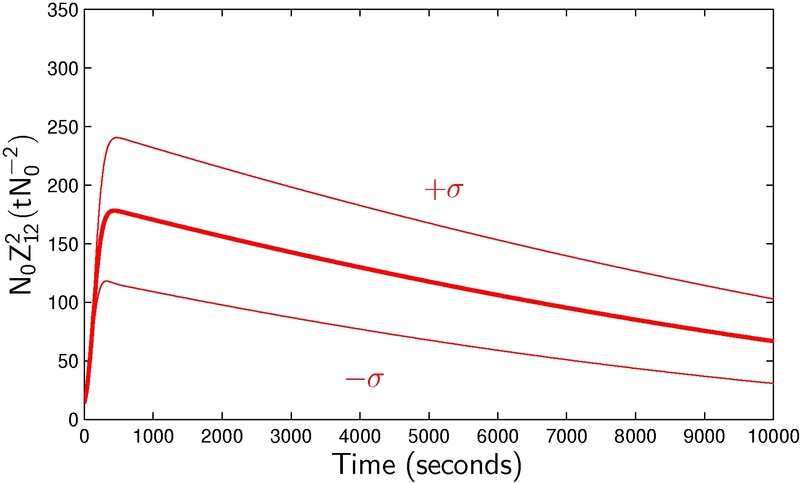}\\
\hspace{-1cm}\includegraphics[totalheight=0.18\textheight,clip]{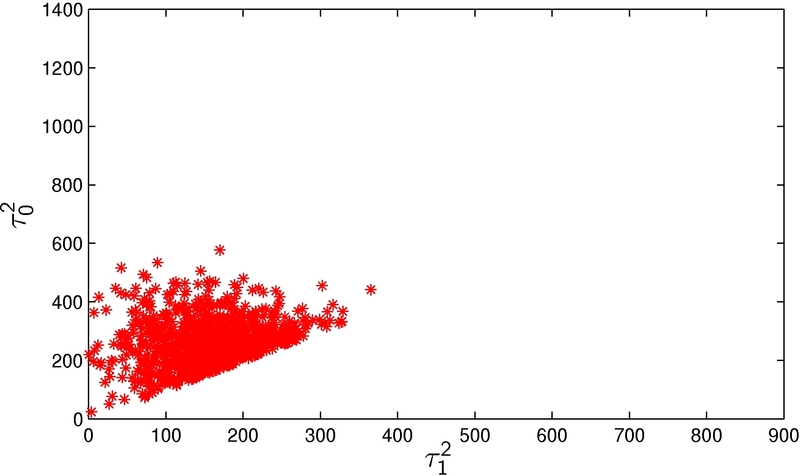}\\
\caption{Approximation using the limiting model for $\gamma =2$ in the 
alternative scaling}
\label{fig: alternativescaling2}
\end{figure}

For $\gamma =2$, the simulation is carried out on the time 
interval $[0,10000]$.  The $3$-dimensional limiting model 
(\ref{altgam2}) is piecewise deterministic  and 
includes the auxiliary variables 
$Z_{12}^2$, $Z_{45}^2$, and the species abundance $Z_6^2$.
$Z_{12}^2$ is governed by a random
differential equation driven by a component of the jump 
process, $Z_{45}^2$.  $Z_{45}^2$ and $Z_6^2$ are discrete with transition 
intensities that depend on $Z_{12}^2$.  Since there is mutual 
dependence between the
continuous and discrete components, we 
modify Gillespie's SSA to simulate the limiting 
system.  Here is a brief description of the simulation 
method for the limiting system.  
\begin{enumerate}
\item Assume that the process has been simulated up to 
$t_i$, the $i$th jump time of the jump process.  Simulate a 
unit exponential random variable $\Delta$ by simulating a
uniform $[0,1]$ random number $r_1$  and setting 
$\Delta =\log\frac 1{r_1}$.

\item Solve the differential equation for $Z_{12}^2$ starting at
$Z_{12}^2(t_i)$ holding $Z_{45}^2(t)=Z_{45}^2(t_i)$ and $Z_6^2(t)
=Z_6^2(t_i)$ until 
time $t_{i+1}$ satisfying
\begin{eqnarray*}
&&\int_{t_i}^{t_{i+1}}\big(\kappa_7\varphi_2(Z_{12}^2(s))\bar {Z}_
5^2(s)+\kappa_8Z_6^2(s)\big)\,ds\\
&&\qquad\qquad =Z_{45}^2(t_i)\int_{t_i}^{t_{i+1}}\frac {\kappa_5\kappa_
7\varphi_2(Z_{12}^2(s))^2}{\kappa_6+\kappa_5\varphi_2(Z_{12}^2(s))}
ds+\kappa_8Z_6^2(t_i)(t_{i+1}-t_i)\\
&&\qquad\qquad =\Delta .\end{eqnarray*}
(We compute the integral by the trapezoid rule using 
the grid points from the ODE solver.)

\item Simulate a uniform $[0,1]$ random number $r_2$.  If 
\begin{eqnarray}
r_2&\leq&\frac {\kappa_7\varphi_2(Z_{12}^2(t_{i+1}))\bar {Z}_5^2(
t_{i+1}-)}{\kappa_7\varphi_2(Z_{12}^2(t_{i+1}))\bar {Z}_5^2(t_{i+
1}-)+\kappa_8Z_6^2(t_{i+1}-))}\label{jumpdet}\\
&=&\frac {\kappa_5\kappa_7\varphi_2(Z_{12}^2(t_{i+1}))^2Z_{45}^2(
t_i)}{\kappa_5\kappa_7\varphi_2(Z_{12}^2(t_{i+1}))^2Z^2_{45}(t_i)
+\kappa_8Z_6^2(t_i)(\kappa_6+\kappa_5\varphi_2(Z_{12}^2(t_{i+1})))},\nonumber\end{eqnarray}
set
\[\left(\begin{array}{c}
Z_{45}^2(t_{i+1})\\
Z^2_6(t_{i+1})\end{array}
\right)=\left(\begin{array}{c}
Z_{45}^2(t_i)\\
Z^2_6(t_i)\end{array}
\right)+\left(\begin{array}{c}
-1\\
1\end{array}
\right),\]
and if the reverse inequality holds in (\ref{jumpdet}), set

\[\left(\begin{array}{c}
Z_{45}^2(t_{i+1})\\
Z^2_6(t_{i+1})\end{array}
\right)=\left(\begin{array}{c}
Z_{45}^2(t_i)\\
Z^2_6(t_i)\end{array}
\right)+\left(\begin{array}{c}
1\\
-1\end{array}
\right).\]
\item Go back to step $1$. 
\end{enumerate}

Comparing plots for $X_1(t)+2X_2(t)$ in Figure \ref{fig:  
goutsias2} and for $N_0Z_{12}^2(tN_0^{-2})$ in Figure \ref{fig:  
alternativescaling2}, the plot in the approximation 
increases more rapidly at early times and starts to drop 
earlier than the plot in the exact simulation.  Also, the 
peak level in the approximation is much lower than the 
peak level in the exact simulation.  

Since $\kappa_8=8.77\times 10^{-8}$ is small compared to the time 
interval, 
Reaction 8 will rarely occur on the time scales we are 
considering.  We retained this reaction in the limiting 
model only to emphasize that a long time after the 
model appears to equilibrate, action may restart after 
the dissociation
\[DNA\cdot 2D\rightharpoonup DNA\cdot D+D.\]

If Reaction 8 does not occur, the stochastic behavior of 
the limiting model just depends on the two jump times
\[\tau_1^2=\inf\{t:Z_{45}^2(t)=1\},\quad\tau_0^2=\inf\{t:Z_{45}^2
(t)=0\},\]
so we compare these random variable to the corresponding 
variables
\[\tau_1=\inf\{t:X_4(t)+X_5(t)=1\},\quad\tau_0=\inf\{t:X_4(t)+X_5
(t)=0\},\]
from the original model or more precisely, because of the 
change of time scale, we compare $(N_0^2\tau_1^2,N_0^2\tau_2^2)$ to $
(\tau_1,\tau_2)$.

In Figure \ref{fig:  goutsias2}, plots for $\tau_1$ and $\tau_0$ for 
$100$ exact simulations are given.  Taking the average, the 
mean of first hitting time of $X_4(t)+X_5(t)$ to $1$ is $305.44$ 
and the mean of the first hitting time of $X_4(t)+X_5(t)$ to 
$0$ is $512.45$.  In Figure \ref{fig:  alternativescaling2}, 
plots for 1000 simulations of $\tau_1^2$ and $\tau_0^2$ are given.  The mean 
of the first hitting time of $Z_{45}^2(tN_0^{-2})$ to $1$ is $155
.95$ and 
the mean of the first hitting time of $Z_{45}^2(tN_0^{-2})$ to $0$ is 
$261.01$.  Comparing the two stopping times in the simulations 
of the full model and of the approximation, the mean 
hitting time to $1$ and $0$ in the approximation is much 
faster than in the full model.  Consequently, the
quicker decrease of $Z_{45}^2$ to $0$ gives a discrepancy in the 
peak levels and the peak times in the full model and in 
the approximation.


\subsection{Derivation of Michaelis-Menten equation}
\cite{Dard79, Dard82} derives the Michaelis-Menten 
equation from a stochastic reaction network model.  His 
result can be obtained as a special case of the methods 
developed here.

Consider the reaction system
\[S_1+S_2\mathop{\rightleftharpoons}^{\kappa_1'}_{\kappa_2'}S_3\mathop{
\rightharpoonup}^{\kappa_3'}S_4+S_2,\]
where $S_1$  is the substrate, $S_2$ the enzyme, $S_3$ the 
enzyme-substrate complex, and $S_4$ the product.  
Assume that the parameters scale so that
\begin{eqnarray*}
Z^N_1(t)&=&Z^N_1(0)-N^{-1}Y_1(N\int_0^t\kappa_1Z^N_1(s)Z^N_2(s)ds
)+N^{-1}Y_2(N\int_0^t\kappa_2Z^N_3(s)ds)\\
Z^N_2(t)&=&Z^N_2(0)-Y_1(N\int_0^t\kappa_1Z^N_1(s)Z^N_2(s)ds)+Y_2(
N\int_0^t\kappa_2Z^N_3(s)ds)\\
&&\qquad +Y_3(N\int_0^t\kappa_3Z^N_3(s)ds)\\
Z_3^N(t)&=&Z^N_2(0)+Y_1(N\int_0^t\kappa_1Z^N_1(s)Z^N_2(s)ds)-Y_2(
N\int_0^t\kappa_2Z^N_3(s)ds)\\
&&\qquad -Y_3(N\int_0^t\kappa_3Z^N_3(s)ds)\\
Z^N_4(t)&=&N^{-1}Y_3(N\int_0^t\kappa_3Z^N_3(s)ds),\end{eqnarray*}
that is, $\alpha_1=\alpha_4=1$, $\alpha_2=\alpha_3=0$, $\beta_1=0$, and $
\beta_2=\beta_3=1$.

Note that $M=Z^N_3(t)+Z^N_2(t)$ is constant, and define
\[V_2^N(t)=\int_0^tZ_2^N(s)ds.\]

\begin{theorem}\label{dard}
Assume that $Z^N_1(0)\rightarrow x_A(0)$.
Then $(Z^N_1,V^N_2)$ converges to $(x_1(t),v_2(t))$ satisfying
\begin{eqnarray}
x_1(t)&=&x_1(0)-\int_0^t\kappa_1x_1(s)\dot {v}_2(s)ds+\int_0^t\kappa_
2(M-\dot {v}_2(s))ds\label{limsys}\\
0&=&-\int_0^t\kappa{}_1x_1(s)\dot {v}_2(s)ds+\int_0^t(\kappa_2+\kappa_
3)(M-\dot {v}_2(s))ds,\nonumber\end{eqnarray}
and hence $\dot {v}_2(s)=\frac {(\kappa_2+\kappa_3)M}{\kappa_2+\kappa_
3+\kappa_1x_1(s)}$ and
\[\dot {x}_1(t)=-\frac {M\kappa_1\kappa_3x_1(t)}{\kappa_2+\kappa_
3+\kappa_1x_1(s)}.\]
\end{theorem}

\begin{proof}
Relative compactness of the sequence $(Z_1^N,V_2^N)$ is 
straightforward.  Dividing the second equation by $N$ and 
passing to the limit, we see that any limit point $(x_1,v_2)$ of 
$(Z_1^N,V_2^N)$ must satisfy
\begin{equation}0=-\int_0^t\kappa{}_1x_1(s)dv_2(s)+(\kappa_2+\kappa_
3)Mt-\int_0^t(\kappa_2+\kappa_3)dv_2(s).\label{limm}\end{equation}
Since $v_2$ is Lipschitz, it is absolutely continuous, and 
rewriting (\ref{limm}) in terms of the derivative gives 
the second equation in (\ref{limsys}).  The first equation 
follows by a similar argument.
\end{proof}

\subsection{Limiting models when the balance conditions 
fail}\label{rawl}
The balance condition, Condition \ref{cndz}, has as its goal 
ensuring that the normalized species numbers remain 
positive, at least on average, and bounded. \citet*{MHR07} 
consider examples in which model reduction is achieved 
by eliminating species whose numbers are zero most of 
the time.  We translate some of their examples into our 
notation and see how one can obtain reduced models 
even though the balance conditions fail.

Consider  
\[S_1\mathop{\rightleftharpoons}^{\kappa'_1}_{\kappa'_2}2S_2,\quad 
S_2\mathop{\rightharpoonup}^{\kappa'_3}S_3,\]
where we assume $\kappa_2',\kappa_3'>>\kappa_1'$.  We take the scaled 
system to be
\begin{eqnarray*}
Z_1^N(t)&=&Z_1(0)-Y_1(\int_0^t\kappa_1Z_1^N(s)ds)+Y_2(N\int_0^t\kappa_
2Z_2^N(s)(Z_2^N(s)-1)ds)\\
Z_2^N(t)&=&Z_2(0)+2Y_1(\int_0^t\kappa_1Z_1^N(s)ds)-2Y_2(N\int_0^t
\kappa_2Z_2^N(s)(Z_2^N(s)-1)ds)-Y_3(N\int_0^t\kappa_3Z_2^N(s)ds)\\
Z_3^N(t)&=&Z_3(0)+Y_3(N\int_0^t\kappa_3Z_2^N(s)ds).\end{eqnarray*}
Consequently, assuming $Z_2^N(0)=0$, for most $t>0$, 
$Z^N_2(t)=0$ and
\[2Y_1(\int_0^t\kappa_1Z_1^N(s)ds)=Y_3(N\int_0^t\kappa_3Z_2^N(s)d
s)+2Y_2(N\int_0^t\kappa_2Z_2^N(s)(Z_2^N(s)-1)ds).\]
To be 
precise, letting $\varLambda$ denote Lebesgue measure and defining
\[\hat {R}^N_2(t)=\int_0^t{\bf 1}_{\{Z_2^N(r-)=2\}}dR_2^N(r),\quad
\hat {R}_3^N(t)=\int_0^t{\bf 1}_{\{Z_2^N(r-)=2\}}dR_3^N(r),\]
for each 
$t>0$,
\begin{eqnarray*}
\lim_{N\rightarrow\infty}\varLambda \{0\leq s\leq t:Z^N_2(s)\neq 
0\}\leq\lim_{N\rightarrow\infty}\int_0^tZ^N_2(s)ds&=&0\\
\limsup_{N\rightarrow\infty}\sup_{s\leq t}Z_2^N(s)&\leq &2\\
\lim_{N\rightarrow\infty}\int_0^t|R_2^N(s)-\hat {R}_2^N(s)|ds&=&0\\
\lim_{N\rightarrow\infty}\int_0^t|R_3^N(s)-2\hat {R}_3^N(s)|ds&=&
0\end{eqnarray*}
so
\[\lim_{N\rightarrow\infty}\int_0^t|R_1^N(s)-\hat {R}_2^N(s)-\hat {
R}_3^N(s)|ds=0.\]
Setting 
 $Q^N(t)={\bf 1}_{\{Z_2^N(t)=2\}}$, 
\[\hat {R}_2^N(t)-\int_0^tNQ^N(s)\kappa_22ds\mbox{\rm \ and }\hat {
R}_3^N(t)-\int_0^tNQ^N(s)\kappa_32ds\]
are martingales.  Working first with a subsequence 
satisfying (\ref{cvreq}), by Lemma \ref{alloc}, $(\hat {R}_2^N,\hat {
R}_3^N)$ 
converges to counting processes $(\hat {R}_2,\hat {R}_3)$ with intensities
\[\hat{\lambda}_2(t)=\frac {\kappa_1\kappa_2}{\kappa_2+\kappa_3}Z_
1(t),\quad\hat{\lambda}_3(t)=\frac {\kappa_1\kappa_3}{\kappa_2+\kappa_
3}Z_1(t),\]
where $Z_1(t)=Z_1(0)-\hat {R}_3(t)$.  It follows that the finite 
dimensional distributions of $(Z_1^N,Z_3^N)$ converge to those 
of a solution to
\begin{eqnarray*}
Z_1(t)&=&Z_1(0)-Y(\int_0^t\frac {\kappa_1\kappa_3}{\kappa_2+\kappa_
3}Z_1(s)ds)\\
Z_3(t)&=&Z_3(0)+2Y(\int_0^t\frac {\kappa_1\kappa_3}{\kappa_2+\kappa_
3}Z_1(s)ds),\end{eqnarray*}
which is the reduced model obtained in \cite{MHR07}.  
More precisely, $(Z_1^N,Z_3^N)$ converges in the Jakubowski 
topology as described in Remark \ref{jakrem}

(Note the relationship between our rate constants 
and those of \cite{MHR07}:  $\kappa_1=k_1$, $\kappa_2=\frac 12k_{
-1}$, and 
$\kappa_3=k_2$.)

\[\]

\renewcommand {\theequation}{A.\arabic{equation}}
\appendix

\setcounter{equation}{0}

\section{Appendix}

\subsection{Convergence of random measures} \label{rmsect}
The material in this section is taken from \cite{Kur92}.  
Proofs of the results can be found there.

Let  $({\Bbb L},d)$  be a complete, separable 
metric space, and let  ${\cal M}({\Bbb L})$  
be the space of finite measures on  ${\Bbb L}$  with 
the weak topology.  The Prohorov metric on  ${\cal M}({\Bbb L})$  is 
defined by 
\begin{equation}\rho (\mu ,\nu )=\inf\{\epsilon >0:\mu (B)\leq\nu 
(B^{\epsilon})+\epsilon ,\nu (B)\leq\mu (B^{\epsilon})+\epsilon ,
B\in {\cal B}({\Bbb L})\},\label{pro}\end{equation}
where  $B^{\epsilon}=\{x\in {\Bbb L}:\inf_{y\in B}d(x,y)<\epsilon 
\}$.
The following lemma is a simple consequence 
of Prohorov's theorem.

\begin{lemma}\ \label{rmcnv}
Let $\{\Gamma_n\}$ be a sequence of ${\cal M}({\Bbb L})$-valued random variables.  
Then $\Gamma_n$ is relatively compact if and only if $\{\Gamma_n(
{\Bbb L})\}$ is 
relatively compact as a family of ${\Bbb R}$-valued random 
variables and for each $\epsilon >0$, there exists a compact 
$K\subset {\Bbb L}$
such that $\sup_nP\{\Gamma_n(K^c)>\epsilon \}<\epsilon$.
\end{lemma}

\begin{corollary}Let 
$\{\Gamma_n\}$ be a sequence of ${\cal M}({\Bbb L})$-valued random variables.  
Suppose that $\sup_nE[\Gamma_n({\Bbb L})]<\infty$ and that for each $
\epsilon >0$, there 
exists a compact $K\subset {\Bbb L}$ such that 
\[\limsup_{n\rightarrow\infty}E[\Gamma_n(K^c)]\leq\epsilon .\]
Then  $\{\Gamma_n\}$  is relatively compact.
\end{corollary}

Let  ${\cal L}({\Bbb L})$  be the space of measures on  ${\Bbb L}
\times [0,\infty )$
  such that  $\mu$$({\Bbb L}\times [0,t])<\infty$    for 
each  $t>0$, and let  ${\cal L}_m({\Bbb L})\subset {\cal L}({\Bbb L}
)$  
be the subspace on which  $\mu ({\Bbb L}\times [0,t])=t$.    
For    $\mu\in {\cal L}({\Bbb L})$, let   $\mu^t$   denote the restriction of $
\mu$   
to  ${\Bbb L}\times [0,t]$.  Let   $\rho_t$
denote the Prohorov metric on  ${\cal M}({\Bbb L}\times [0,t])$, and define   
$\hat{\rho}$ on ${\cal L}({\Bbb L})$ by
\[\hat{\rho }(\mu ,\nu )=\int_0^{\infty}e^{-t}1\wedge\rho_t(\mu^t
,\nu^t)dt,\]
that is, $\{\mu_n\}$ converges in $\hat{\rho}$ if and only if $\{
\mu_n^t\}$ converges 
weakly for almost every $t$.
In particular, if $\hat{\rho }(\mu_n,\mu )\rightarrow 0$, then  $
\rho_t(\mu_n^t,\mu^t)\rightarrow 0$ if and only 
if 
$\mu_n({\Bbb L}\times [0,t])\rightarrow\mu ({\Bbb L}\times [0,t])$.  The following lemma 
is an immediate consequence of Lemma \ref{rmcnv}.

\begin{lemma}\ \label{titmeas}
A sequence of $({\cal L}_m({\Bbb L}),\hat{\rho })$-valued random variables $
\{\Gamma_n\}$ is 
relatively compact if and only if for each $\epsilon >0$ and each 
$t>0$,
there exists a compact $K\subset {\Bbb L}$ such that 
$\inf_nE[\Gamma_n(K\times [0,t])]\geq (1-\epsilon )t$.
\end{lemma}

\begin{lemma}\label{lem1a4}
  Let  $\Gamma$  be an $({\cal L}({\Bbb L}),\hat{\rho })$-valued random variable adapted to a 
complete filtration  $\{{\cal F}_t\}$  in the sense that for each  
$t\geq 0$ and $H\in {\cal B}({\Bbb L})$,  $\Gamma (H\times [0,t])$  
is  ${\cal F}_t$-measurable.  Let  $\lambda (G)=\Gamma ({\Bbb L}\times 
G)$.  Then there exists 
an $\{{\cal F}_t\}$-optional, ${\cal P}({\Bbb L})$-valued process   $
\gamma$ such that  
\begin{equation}\int_{{\Bbb L}\times [0,t]}h(y,s)\Gamma (dy\times 
ds)=\int_0^t\int_{{\Bbb L}}h(y,s)\gamma_s(dy)\lambda (ds)\label{disint}\end{equation}
for all  $h\in B({\Bbb L}\times [0,\infty ))$  with probability one.  
If  $\lambda ([0,t])$  is continuous, 
then  $\gamma$  can be taken to be $\{{\cal F}_t\}$-predictable.
\end{lemma}

\begin{lemma}\label{contlem}
Let $\{(x_n,\mu_n)\}\subset D_{{\Bbb E}}[0,\infty )\times {\cal L}
({\Bbb L})$, and  $(x_n,\mu_n)\rightarrow (x,\mu )$.  Let 
$h\in C({\Bbb E}\times {\Bbb L})$ and $\psi$
be a nonnegative function on $[0,\infty )$ satisfying 
$\lim_{r\rightarrow\infty}\psi (r)/r=\infty$ 
such that 
\begin{equation}\sup_n\int_{{\Bbb L}\times [0,t]}\psi (|h(x_n(s),
y)|)\mu_n(dy\times ds)<\infty\label{uint2}\end{equation}
for each  $t>0$.

  Define
\[u_n(t)=\int_{{\Bbb L}\times [0,t]}h(x_n(s),y)\mu_n(dy\times ds)
,\quad u(t)=\int_{{\Bbb L}\times [0,t]}h(x(s),y)\mu (dy\times ds)\]
$z_n(t)=\mu_n({\Bbb L}\times [0,t])$, and $z(t)=\mu ({\Bbb L}\times 
[0,t])$.

\begin{itemize}
\item[a)] If  $x$  is continuous on  $[0,t]$  and  
$\lim_{n\rightarrow\infty}z_n(t)=z(t)$, then  $\lim_{n\rightarrow
\infty}u_n(t)=u(t)$.

\item[b)]  If  $(x_n,z_n,\mu_n)\rightarrow (x,z,\mu )$ in $D_{{\Bbb E}
\times {\Bbb R}}[0,\infty )\times {\cal L}({\Bbb L})$,
then $(x_n,z_n,u_n,\mu_n)\rightarrow (x,z,u,\mu )$ in $D_{{\Bbb E}
\times {\Bbb R}\times {\Bbb R}}[0,\infty )\times {\cal L}({\Bbb L}
)$.
In particular, $\lim_{n\rightarrow\infty}u_n(t)=u(t)$  at 
all points of continuity of  $z$.

\item[c)] The continuity assumption on $h$ can be replaced 
by the assumption that $h$ is continuous a.e.  $\nu_t$ for each 
$t$, where $\nu_t\in {\cal M}({\Bbb E}\times {\Bbb L})$ is the measure determined by 
$\nu_t(A\times B)=\mu \{(y,s):x(s)\in A,s\leq t,y\in B\}$. 

\end{itemize}
\end{lemma}

The Lemma \ref{contlem}\ and the continuous mapping 
theorem give the following.

\begin{lemma}\label{wkconv}
Suppose $(Z^N,V^N)\Rightarrow (Z,V)$ in $D_{{\Bbb E}}[0,\infty )\times 
{\cal L}_m({\Bbb L})$. Let 
$h\in C({\Bbb E}\times {\Bbb L})$ and $\psi$ be as in Lemma \ref{contlem}. 
If
$\{\int_0^t\psi (|h(Z^N(s),y)|)V^N(dy\times ds)\}$is stochastically bounded for 
all $t>0$, then
\[\int_{{\Bbb L}\times [0,\cdot ]}h(Z^N(s),y)V^N(dy\times ds)\Rightarrow
\int_{{\Bbb L}\times [0,\cdot ]}h(Z(s),y))V(dy\times ds).\]
\end{lemma}

\subsection{Martingale properties of counting processes}
A cadlag stochastic process $R$ is a {\em counting process\/} if 
$R(0)=0$ and $R$ is constant except for jumps of plus one.  
If $R$ is adapted to a filtration $\{{\cal F}_t\}$, then a nonnegative 
$\{{\cal F}_t\}$-adapted process $\lambda$ is an $\{{\cal F}_t\}$-intensity for $
R$ if
\[M(t)=R(t)-\int_0^t\lambda (s)ds\]
is an $\{{\cal F}_t\}$-local martingale.  Specifically, letting $
\tau_l$ denote 
the $l$th jump time of $R$,
\[M^{\tau_l}(t)\equiv M(t\wedge\tau_l)=R(t\wedge\tau_l)-\int_0^{t
\wedge\tau_l}\lambda (s)ds\]
is an $\{{\cal F}_t\}$-martingale for each $l$.

For simplicity, we assume that $\lambda$ is cadlag.

\begin{remark}
For $R_k$ defined in (\ref{cntr}) and 
$\{{\cal F}_t\}=\sigma (R_l(s):s\leq t,l=1,\ldots ,r_0)$, the intensity for $
R_k$ 
is 
$t\rightarrow\lambda_k(X(t))$.
\end{remark}

\begin{lemma}
For each $t\geq 0$ and each $l$,
\begin{equation}l\geq E[R(t\wedge\tau_l)]=E[\int_0^{t\wedge\tau_l}
\lambda (s)ds]\label{mnint}\end{equation}
and
\[E[R(t)]=E[\int_0^t\lambda (s)ds],\]
where we allow $\infty =\infty$.
If $E[R(t)]<\infty$ for all $t>0$, then 
\[R(t)-\int_0^t\lambda (s)ds\]
is an  $\{{\cal F}_t\}$-martingale. 
\end{lemma}

Two counting processes, $R_1$, $R_2$, are {\em orthogonal\/} if they have no 
simultaneous jumps.  

\begin{lemma}
Let $R_1,\ldots ,R_m$ be pairwise orthogonal $\{{\cal F}_t\}$-adapted 
counting processes with $\{{\cal F}_t\}$-intensities $\lambda_k$.  
Then, perhaps on a larger probability space,
 there exist independent unit Poisson processes 
$Y_1,\ldots ,Y_m$ such 
that 
\[R_k(t)=Y_k(\int_0^t\lambda_k(s)ds),\]
and $R=\sum_{k=1}^mR_k$ is a counting process with intensity 
$\lambda =\sum_{k=1}^m\lambda_k$.  

If $\tau_l$ is the $l$th jump time of $R$, then 
\begin{equation}P\{R_k(\tau_l)-R_k(\tau_l-)=1|{\cal F}_{\tau_l}\}
=\frac {\lambda_k(\tau_l-)}{\lambda (\tau_l-)}.\label{jmpprb}\end{equation}
\end{lemma}

\begin{remark}
Note that the right side of (\ref{jmpprb}) involves the 
left limits of the intensities.  If the intensities are not 
cadlag, then $\lambda_k(\tau_l-)$ should be replaced by
\[\limsup_{h\rightarrow 0+}h^{-1}\int_{\tau_l-h}^{\tau_l}\lambda_
k(s)ds.\]

The intensity of a counting process does not necessarily 
uniquely determined its distribution.  For example, 
consider the system
\begin{eqnarray*}
R_1(t)&=&Y_1(\int_0^t\lambda (R_1(s))ds)\\
R_2(t)&=&Y_2(\int_0^t\lambda (R_1(s))ds).\end{eqnarray*}
The intensity for each component is $\lambda (R_1(t))$, but the two 
components will not have the same distribution.
\end{remark}

\begin{proof}
See  \cite{Mey71} and \cite{Kur80a}.
\end{proof}

\begin{lemma}\label{skconv}
Suppose that $R_1^N,\ldots ,R_m^N$ are pairwise orthogonal
counting processes adapted to a filtration $\{{\cal F}_t^N\}$ with 
$\{{\cal F}_t^N\}$-intensities $\lambda_1^N,\ldots ,\lambda_m^N$.  Let $
\Lambda_k^N(t)=\int_0^t\lambda_k^N(s)ds$, and 
suppose that $(\Lambda_1^N,\ldots ,\Lambda_m^N)\Rightarrow (\Lambda_
1,\ldots ,\Lambda_m)$ in the Skorohod 
topology on $D_{{\Bbb R}^m}[0,\infty )$.  Then 
$\{(R_1^N,\ldots ,R_m^N)\}$ is relatively compact in the Skorohod 
topology  and any limit point $(R_1,\ldots ,R_m)$ consists of 
pairwise orthogonal counting processes.  

At least along a further subsequence,  
\[(\Lambda_1^N,\ldots ,\Lambda_m^N,R_1^N,\ldots ,R_m^N)\Rightarrow 
(\Lambda_1,\ldots ,\Lambda_m,R_1,\ldots ,R_m),\]
and letting $\{{\cal F}_t^{\Lambda ,R}\}$ be the filtration generated by 
$(\Lambda_1,\ldots ,\Lambda_m,R_1,\ldots ,R_m)$,  $R_k-\Lambda_k$ are $
\{{\cal F}_t^{\Lambda ,R}\}$-local martingales and 
there exist independent unit Poisson processes $(Y_1,\ldots ,Y_m)$ 
such that
\begin{equation}R_k(t)=Y_k(\Lambda_k(t)),\quad k=1,\ldots ,m.\label{cpsys}\end{equation}
\end{lemma}

\begin{remark}
If the $\Lambda_k$ are adapted to $\{{\cal F}_t^R\}$, then $R$ will be the 
unique solution of (\ref{cpsys}) and $R^N\Rightarrow R$ in the 
Skorohod topology.
\end{remark}

\begin{proof}
See \citet*{KLS84}.
\end{proof}

In Section \ref{rawl}, we consider an example for which 
the integrated intensities did not have a continuous 
limit.  The next lemma covers that situation.

\begin{lemma}\label{alloc}
Suppose that $R_0^N,R_1^N,\ldots ,R_m^N$ are counting processes 
adapted to a filtration $\{{\cal F}_t^N\}$, and 
$R_1^N,\ldots ,R_m^N$ are pairwise orthogonal.
Suppose 
$R_0^N$ has $\{{\cal F}_t^N\}$-intensity $\lambda_0^N$, and $R_1^
N,\ldots ,R_m^N$ have $\{{\cal F}_t^N\}$-intensities 
$\lambda_k^N=NQ^N\mu_k^N$, where $Q^N\geq 0$.
Suppose 
\begin{equation}(\lambda_0^N,\mu_1^N,\ldots ,\mu_m^N)\Rightarrow 
(\lambda_0,\mu_1,\ldots ,\mu_m),\label{cvreq}\end{equation}
and 
\begin{equation}\int_0^t|R_0^N(s)-\sum_{k=1}^mR_k^N(s)|ds\rightarrow 
0,\label{ctasym}\end{equation}
for each $t>0$.  Then $\{(R_0^N,R_1^N,\ldots ,R_m^N)\}$ is relatively 
compact in the Jakubowski topology and for any limit point 
$(R_0,R_1,\ldots ,R_m)$,
\[R_0=\sum_{k=1}^mR_k,\]
and
$R_1,\ldots ,R_m$ are pairwise 
orthogonal counting processes with intensities
\[\lambda_k(t)=\frac {\mu_k(t)}{\sum_{l=1}^m\mu_l(t)}\lambda_0(t)
.\]
\end{lemma}

\begin{remark}\label{jakrem}
The sequence may not be relatively compact in the Skorohod 
topology since we have not ruled out the possibility that
the sequence has discontinuities that 
coalesce.  See the example in Section \ref{rawl}.

The Meyer-Zheng conditions (\cite{MZ84}) imply relative 
compactness in the Jakubowski topology (\cite{Jak97}).
A sequence of cadlag functions $\{x_n\}$ 
converges to a cadlag function $x$ in the Jakubowski topology 
 if and only if  there exists a sequence of 
time changes $\{\gamma_n\}$ such that $(x_n\circ\gamma_n,\gamma_n
)\rightarrow (x\circ\gamma ,\gamma )$ in the 
Skorohod topology.  (See 
\cite{Kur91}.)  The 
time-changes are continuous, nondecreasing mappings from 
$[0,\infty )$ onto $[0,\infty )$ but are not necessarily strictly 
increasing.  Convergence implies
$\int_0^t|x_n(s)-x(s)|\wedge 1ds\rightarrow 0$.  In contrast to the Skorohod 
topology, if $x_n\rightarrow x$ and $y_n\rightarrow y$ in the Jakubowski topology, 
then $(x_n,y_n)\rightarrow (x,y)$ in the Jakubowski topology on cadlag 
functions in the product space.
\end{remark}

\begin{proof}
By Lemma \ref{skconv}, $\{R_0^N\}$ is relatively compact in 
the Skorohod topology and hence in the Jakubowski 
topology.  Let 
\[\hat {R}_0^N=\sum_{k=1}^mR_k^N.\]
The stochastic boundedness of $\{R_0^N(t)\}$ for each $t>0$
 and (\ref{ctasym}) 
imply the stochastic boundedness of $\{\hat {R}_0^N(t)\}$ for each 
$t>0$ which by (\ref{mnint}) implies the stochastic 
boundedness of 
\[\{\int_0^tNQ^N(s)\sum_{k=1}^m\mu_k^N(s)ds\}.\]
Let $\gamma_N$ be defined by 
\[\int_0^{\gamma_N(t)}(1+NQ^N(s)\sum_{k=1}^m\mu_k^N(s))ds=t.\]
Since $|\gamma_N(s)-\gamma_N(t)|\leq |s-t|$,
 $\{\gamma_N\}$ is relatively compact.  
Define 
\[\Lambda_k^N(t)=\int_0^t\lambda_k^N(s)ds,\]
and observe that 
\[\Lambda_l^N\circ\gamma_N(t)=\int_0^t\frac {NQ^N\circ\gamma_N(s)
\mu_l^N\circ\gamma_N(s)}{1+NQ^N\circ\gamma_N(s)\sum_k\mu_k^N\circ
\gamma_N(s)}ds\]
is also Lipschitz with Lipschitz constant $1$.  Since 
$\{\gamma_N(t),t\geq 0\}$ are stopping times, 
\[R_l^N-\Lambda_l^N\circ\gamma_N\]
are martingales with respect to the filtration $\{{\cal F}_{\gamma_
N(t)}^N\}$.  

The Lipschitz properties imply the relative compactness 
of $\{(\Lambda_1^N\circ\gamma_N,\ldots ,\Lambda_m^N\circ\gamma_N,
\gamma_N)\}$ in the Skorohod topology
 which in turn, by Lemma \ref{skconv}, implies the 
relative compactness of 
\[\{(\Lambda_1^N\circ\gamma_N,\ldots ,\Lambda_m^N\circ\gamma_N,\gamma_
N,R_1^N\circ\gamma_N,\ldots ,R_m^N\circ\gamma_N)\}.\]
Relative compactness  of this sequence in the Skorohod 
topology ensures relative compactness of $\{(R_1^N,\ldots ,R_m^N)
\}$ in 
the Jakubowski topology which in turn implies relative 
compactness of $\{(R_0^N,R_1^N,\ldots ,R_m^N)\}$ in the Jakubowski 
topology.

Along an appropriate subsequence, we have convergence 
of $\gamma_N$ to a limit $\gamma$,
\[\int_0^t\frac {NQ^N\circ\gamma_N(s)\sum_k\mu_k^N\circ\gamma_N(s
)}{1+NQ^N\circ\gamma_N(s)\sum_k\mu_k^N\circ\gamma_N(s)}ds\Rightarrow
\hat{\Lambda },\]
convergence of $\Lambda_k^N\circ\gamma_N$ to
\[\hat{\Lambda}_k(t)=\int_0^t\frac {\mu_k\circ\gamma (s)}{\sum_l\mu_
l\circ\gamma (s)}d\hat{\Lambda }(ds),\]
and convergence of $(R^N_0,R_1^N,\ldots ,R_m^N)$ in the Jakubowski topology
to a process satisfying 
\[R_0=\sum_{k=1}^mR_k.\]
Since $R_0\circ\gamma (t)-\int_0^{\gamma (t)}\lambda_0(s)ds$ is a martingale, we must
have 
\[\int_0^{\gamma (t)}\lambda_0(s)ds=\hat{\Lambda }(t)\]
and 
\[\hat{\Lambda}_k(t)=\int_0^t\frac {\mu_k\circ\gamma (s)}{\sum_l\mu_
l\circ\gamma (s)}\lambda_0\circ\gamma (s)\gamma'(s)ds.\]
Since $R_0$ is a counting process, the $R_k$, $k=1,\ldots ,m$, must 
be orthogonal, and $R_k$ must have intensity 
$\frac {\mu_k}{\sum_l\mu_l}\lambda_0$.
 \end{proof}

\bibliography{chem2}

\begin{thebibliography}{32}
\providecommand{\natexlab}[1]{#1}
\providecommand{\url}[1]{\texttt{#1}}
\expandafter\ifx\csname urlstyle\endcsname\relax
  \providecommand{\doi}[1]{doi: #1}\else
  \providecommand{\doi}{doi: \begingroup \urlstyle{rm}\Url}\fi

\bibitem[Ball et~al.(2006)Ball, Kurtz, Popovic, and Rempala]{BKPR06}
Karen Ball, Thomas~G. Kurtz, Lea Popovic, and Greg Rempala.
\newblock Asymptotic analysis of multiscale approximations to reaction
  networks.
\newblock \emph{Ann. Appl. Probab.}, 16\penalty0 (4):\penalty0 1925--1961,
  2006.
\newblock ISSN 1050-5164.

\bibitem[Cao et~al.(2005)Cao, Gillespie, and Petzold]{CGP05}
Yang Cao, Daniel~T. Gillespie, and Linda~R. Petzold.
\newblock The slow-scale stochastic simulation algorithm.
\newblock \emph{The Journal of Chemical Physics}, 122\penalty0 (1):\penalty0
  014116, 2005.
\newblock URL \url{http://link.aip.org/link/?JCP/122/014116/1}.

\bibitem[Crudu et~al.(2009)Crudu, Debussche, and Radulescu]{CDR09}
Alina Crudu, Arnaud Debussche, and Ovidiu Radulescu.
\newblock Hybrid stochastic simplifications for multiscale gene networks.
\newblock \emph{BMC Systems Biology}, 3:89, 2009.
\newblock \doi{10.1186/1752-0509-3-89}.

\bibitem[Darden(1979)]{Dard79}
Thomas Darden.
\newblock A pseudo-steady state approximation for stochastic chemical kinetics.
\newblock \emph{Rocky Mountain J. Math.}, 9\penalty0 (1):\penalty0 51--71,
  1979.
\newblock ISSN 0035-7596.
\newblock Conference on Deterministic Differential Equations and Stochastic
  Processes Models for Biological Systems (San Cristobal, N.M., 1977).

\bibitem[Darden(1982)]{Dard82}
Thomas~A. Darden.
\newblock Enzyme kinetics: stochastic vs. deterministic models.
\newblock In \emph{Instabilities, bifurcations, and fluctuations in chemical
  systems (Austin, Tex., 1980)}, pages 248--272. Univ. Texas Press, Austin, TX,
  1982.

\bibitem[Davis(1993)]{Dav93}
M.~H.~A. Davis.
\newblock \emph{Markov models and optimization}, volume~49 of \emph{Monographs
  on Statistics and Applied Probability}.
\newblock Chapman \& Hall, London, 1993.
\newblock ISBN 0-412-31410-X.

\bibitem[E et~al.(2005)E, Liu, and Vanden-Eijnden]{ELV05}
Weinan E, Di~Liu, and Eric Vanden-Eijnden.
\newblock Nested stochastic simulation algorithm for chemical kinetic systems
  with disparate rates.
\newblock \emph{The Journal of Chemical Physics}, 123\penalty0 (19):\penalty0
  194107, 2005.
\newblock URL \url{http://link.aip.org/link/?JCP/123/194107/1}.

\bibitem[E et~al.(2007)E, Liu, and Vanden-Eijnden]{ELV07}
Weinan E, Di~Liu, and Eric Vanden-Eijnden.
\newblock Nested stochastic simulation algorithms for chemical kinetic systems
  with multiple time scales.
\newblock \emph{J. Comput. Phys.}, 221\penalty0 (1):\penalty0 158--180, 2007.
\newblock ISSN 0021-9991.

\bibitem[Ethier and Kurtz(1986)]{EK86}
Stewart~N. Ethier and Thomas~G. Kurtz.
\newblock \emph{Markov processes}.
\newblock Wiley Series in Probability and Mathematical Statistics: Probability
  and Mathematical Statistics. John Wiley \& Sons Inc., New York, 1986.
\newblock ISBN 0-471-08186-8.
\newblock Characterization and convergence.

\bibitem[Gillespie(1977)]{Gill77}
Daniel~T. Gillespie.
\newblock Exact stochastic simulation of coupled chemical reactions.
\newblock \emph{J. Phys. Chem.}, 81:\penalty0 2340--61, 1977.

\bibitem[Goutsias(2005)]{Gou05}
John Goutsias.
\newblock Quasiequilibrium approximation of fast reaction kinetics in
  stochastic biochemical systems.
\newblock \emph{The Journal of Chemical Physics}, 122\penalty0 (18):\penalty0
  184102, 2005.
\newblock \doi{10.1063/1.1889434}.
\newblock URL \url{http://link.aip.org/link/?JCP/122/184102/1}.

\bibitem[Haseltine and Rawlings(2002)]{HR02}
Eric~L. Haseltine and James~B. Rawlings.
\newblock Approximate simulation of coupled fast and slow reactions for
  stochastic chemical kinetics.
\newblock \emph{J. Chem. Phys.}, 117\penalty0 (15):\penalty0 6959--6969, 2002.

\bibitem[Hensel et~al.(2009)Hensel, Rawlings, and Yin]{HRY09}
Sebastian~C. Hensel, James~B. Rawlings, and John Yin.
\newblock Stochastic kinetic modeling of vesicular stomatitis virus
  intracellular growth.
\newblock \emph{Bull. Math. Biol.}, 71\penalty0 (7):\penalty0 1671--1692, 2009.
\newblock ISSN 0092-8240.
\newblock \doi{10.1007/s11538-009-9419-5}.
\newblock URL
  \url{http://dx.doi.org.ezproxy.library.wisc.edu/10.1007/s11538-009-9419-5}.

\bibitem[Jakubowski(1997)]{Jak97}
Adam Jakubowski.
\newblock A non-{S}korohod topology on the {S}korohod space.
\newblock \emph{Electron. J. Probab.}, 2:\penalty0 no.\ 4, 21 pp.\
  (electronic), 1997.
\newblock ISSN 1083-6489.
\newblock URL
  \url{http://www.math.washington.edu/~ejpecp/EjpVol2/paper4.abs.html}.

\bibitem[Kabanov et~al.(1984)Kabanov, Liptser, and Shiryaev]{KLS84}
Yu.~M. Kabanov, R.~Sh. Liptser, and A.~N. Shiryaev.
\newblock Weak and strong convergence of the distributions of counting
  processes.
\newblock \emph{Theory of Probability and its Applications}, 28\penalty0
  (2):\penalty0 303--336, 1984.
\newblock \doi{10.1137/1128026}.
\newblock URL \url{http://link.aip.org/link/?TPR/28/303/1}.

\bibitem[Kang(2009)]{Kan10}
Hye-Won Kang.
\newblock The multiple scaling approximation in the heat shock model of e.
  coli.
\newblock In Preparation, 2009.

\bibitem[Kang et~al.(2010)Kang, Kurtz, and Popovic]{KKP10}
Hye-Won Kang, Thomas~G. Kurtz, and Lea Popovic.
\newblock Diffusion approximations for multiscale chemical reaction models.
\newblock in preparation, 2010.

\bibitem[Khas{\cprime}minski{\u\i}(1966{\natexlab{a}})]{Kha66a}
R.~Z. Khas{\cprime}minski{\u\i}.
\newblock On stochastic processes defined by differential equations with a
  small parameter.
\newblock \emph{Theory Probab. Appl.}, 11:\penalty0 211--228,
  1966{\natexlab{a}}.

\bibitem[Khas{\cprime}minski{\u\i}(1966{\natexlab{b}})]{Kha66b}
R.~Z. Khas{\cprime}minski{\u\i}.
\newblock A limit theorem for the solutions of differential equations with
  random right-hand sides.
\newblock \emph{Theory Probab. Appl.}, 11:\penalty0 390--406,
  1966{\natexlab{b}}.

\bibitem[Kurtz(1972)]{Kur72}
Thomas~G. Kurtz.
\newblock The relationship between stochastic and deterministic models for
  chemical reactions.
\newblock \emph{J. Chem. Phys.}, 57\penalty0 (7):\penalty0 2976--2978, 1972.

\bibitem[Kurtz(1977/78)]{Kur77}
Thomas~G. Kurtz.
\newblock Strong approximation theorems for density dependent {M}arkov chains.
\newblock \emph{Stochastic Processes Appl.}, 6\penalty0 (3):\penalty0 223--240,
  1977/78.

\bibitem[Kurtz(1980)]{Kur80a}
Thomas~G. Kurtz.
\newblock Representations of {M}arkov processes as multiparameter time changes.
\newblock \emph{Ann. Probab.}, 8\penalty0 (4):\penalty0 682--715, 1980.
\newblock ISSN 0091-1798.
\newblock URL
  \url{http://links.jstor.org/sici?sici=0091-1798(198008)8:4<682:ROMPAM>2.0.CO%
;2-W&origin=MSN}.

\bibitem[Kurtz(1991)]{Kur91}
Thomas~G. Kurtz.
\newblock Random time changes and convergence in distribution under the
  {M}eyer-{Z}heng conditions.
\newblock \emph{Ann. Probab.}, 19\penalty0 (3):\penalty0 1010--1034, 1991.
\newblock ISSN 0091-1798.
\newblock URL
  \url{http://links.jstor.org/sici?sici=0091-1798(199107)19:3<1010:RTCACI>2.0.%
CO;2-2&origin=MSN}.

\bibitem[Kurtz(1992)]{Kur92}
Thomas~G. Kurtz.
\newblock Averaging for martingale problems and stochastic approximation.
\newblock In \emph{Applied stochastic analysis (New Brunswick, NJ, 1991)},
  volume 177 of \emph{Lecture Notes in Control and Inform. Sci.}, pages
  186--209. Springer, Berlin, 1992.

\bibitem[Macnamara et~al.(2007)Macnamara, Burrage, and Sidje]{MBS07}
Shev Macnamara, Kevin Burrage, and Roger~B. Sidje.
\newblock Multiscale modeling of chemical kinetics via the master equation.
\newblock \emph{Multiscale Model. Simul.}, 6\penalty0 (4):\penalty0 1146--1168,
  2007.
\newblock ISSN 1540-3459.

\bibitem[Mastny et~al.(2007)Mastny, Haseltine, and Rawlings]{MHR07}
Ethan~A. Mastny, Eric~L. Haseltine, and James~B. Rawlings.
\newblock Two classes of quasi-steady-state model reductions for stochastic
  kinetics.
\newblock \emph{The Journal of Chemical Physics}, 127\penalty0 (9):\penalty0
  094106, 2007.
\newblock \doi{10.1063/1.2764480}.
\newblock URL \url{http://link.aip.org/link/?JCP/127/094106/1}.

\bibitem[Meyer(1971)]{Mey71}
P.~A. Meyer.
\newblock D\'emonstration simplifi\'ee d'un th\'eor\`eme de {K}night.
\newblock In \emph{S\'eminaire de {P}robabilit\'es, {V} ({U}niv. {S}trasbourg,
  ann\'ee universitaire 1969--1970)}, pages 191--195. Lecture Notes in Math.,
  Vol. 191. Springer, Berlin, 1971.

\bibitem[Meyer and Zheng(1984)]{MZ84}
P.-A. Meyer and W.~A. Zheng.
\newblock Tightness criteria for laws of semimartingales.
\newblock \emph{Ann. Inst. H. Poincar\'e Probab. Statist.}, 20\penalty0
  (4):\penalty0 353--372, 1984.
\newblock ISSN 0246-0203.
\newblock URL \url{http://www.numdam.org/item?id=AIHPB_1984__20_4_353_0}.

\bibitem[Rao and Arkin(2003)]{RA03}
Christopher~V. Rao and Adam~P. Arkin.
\newblock Stochastic chemical kinetics and the quasi-steady-state assumption:
  application to the gillespie algorithm.
\newblock \emph{J. Chem. Phys.}, 118\penalty0 (11):\penalty0 4999--5010, 2003.

\bibitem[Segel and Slemrod(1989)]{SS89}
Lee~A. Segel and Marshall Slemrod.
\newblock The quasi-steady-state assumption: a case study in perturbation.
\newblock \emph{SIAM Rev.}, 31\penalty0 (3):\penalty0 446--477, 1989.
\newblock ISSN 0036-1445.
\newblock \doi{10.1137/1031091}.
\newblock URL \url{http://dx.doi.org.ezproxy.library.wisc.edu/10.1137/1031091}.

\bibitem[Srivastava et~al.(2001)Srivastava, Peterson, and Bentley]{SPB01}
R.~Srivastava, M.~S. Peterson, and W.~E. Bentley.
\newblock Stochastic kinetic analysis of escherichia coli stress circuit using
  sigma(32)-targeted antisense.
\newblock \emph{Biotechnol. Bioeng.}, 75:\penalty0 120--129, 2001.

\bibitem[Zeiser et~al.(2010)Zeiser, Franz, and Liebscher]{ZFL10}
Stefan Zeiser, Uwe Franz, and Volkmar Liebscher.
\newblock Autocatalytic genetic networks modeled by piecewise-deterministic
  {M}arkov processes.
\newblock \emph{J. Math. Biol.}, 60\penalty0 (2):\penalty0 207--246, 2010.
\newblock ISSN 0303-6812.
\newblock \doi{10.1007/s00285-009-0264-9}.
\newblock URL \url{http://dx.doi.org/10.1007/s00285-009-0264-9}.

\end{thebibliography}
\end{document}